\documentstyle[epsf]{article}

\newcommand{\ee}{\end{equation}}
\newcommand{\be}{\begin{equation}}

\begin{document}

\leftline{Running head: {\bf Elliptic Representation}}
\vskip 0.5 cm

\centerline{\bf\huge  The Elliptic Representation}
\vskip 0.5  cm 
\centerline{\bf \huge of the General  Painlev\'e 6 Equation   }
\vskip 0.3 cm
\centerline{\bf   Davide Guzzetti}
\vskip 0.3 cm
\centerline{\bf  Research Institute for Mathematical Sciences (RIMS)}
\centerline{\bf  Kyoto University}
\centerline{\bf Kitashirakawa, Sakyo-ku, Kyoto 606-8502, Japan}
\vskip 0.2 cm 
\centerline{\bf E-mail: guzzetti@kurims.kyoto-u.ac.jp}
\vskip 0.2 cm 
\vskip 0.5 cm 
%
%
%
%
%
%


\vskip 1 cm 
AMS CLASSIFICATION: 34M55

\section{Introduction}\label{Introduction}

In this paper we study the elliptic representation 
of the  sixth Painlev\'e equation  
$$
{d^2y \over dx^2}={1\over 2}\left[ 
{1\over y}+{1\over y-1}+{1\over y-x}
\right]
           \left({dy\over dx}\right)^2
-\left[
{1\over x}+{1\over x-1}+{1\over y-x}
\right]{dy \over dx}
$$
$$
+
{y(y-1)(y-x)\over x^2 (x-1)^2}
\left[
\alpha+\beta {x\over y^2} + \gamma {x-1\over (y-1)^2} +\delta
{x(x-1)\over (y-x)^2}
\right]
,~~~~~\hbox{(PVI)}.
$$
 Though the elliptic representation of PVI has been known since R.Fuchs \cite{fuchs}, in the 
literature  
there is no general study of its analytic implications.  To fill this 
gap, we study here 
 the analytic properties of the solutions in elliptic representation  
for all values of $\alpha, \beta,\gamma,\delta$
and 
we derive their critical behavior close to the singular points $x=0,1,\infty$. Moreover,   we solve 
the connection problem  for generic values 
of $\alpha,\beta,\gamma,\delta$ and for the special (non-generic) case $\beta=\gamma=1-2\delta=0$, which is important in 2-D topological field theory.  
\vskip 0.2 cm

The six classical Painlev\'e  equations were discovered by Painlev\'e \cite{pain} 
 and 
Gambier \cite{gamb}, 
who classified all the second order ordinary differential equations
of the type 
$$
         {d^2 y \over dx^2}= {\cal R}\left(x,y,{dy\over dx}\right)
$$
where ${\cal R}$ is rational in ${dy\over dx}$,  $x$ and
 $y$. The Painlev\'e  equations 
satisfy the {\it Painlev\'e property} of  absence of movable branch points and essential singularities. These singularities will be called {\it critical points}; for PVI  
they are 0,1,$\infty$. The behavior of a solution close to a critical point is called {\it critical 
behavior}.    
 A solution of the sixth Painlev\'e equation  
 can be analytically continued to a meromorphic function on the universal
   covering of ${\bf P}^1\backslash \{ 0,1 ,\infty \}$. 
 For generic values of the integration constants and of the parameters in
   the equation, it can not be expressed via elementary or classical
   transcendental functions. For this reason, it  is called a {\it
   Painlev\'e transcendent}.

 The Painlev\'e equations turned out to be important in physical models, since the works  \cite{Mccoy1} \cite{Mccoy2} \cite{Mccoy3}. The reader may find a 
review of the physical applications in \cite{IN}. More recently, 
the Painlev\'e 6 equation  found applications in topological field theory and Frobenius
  manifolds \cite{Dub1}, \cite{Dub2} 

\vskip 0.2 cm
  The first analytical problem with Painlev\'e equations is to determine 
the critical behavior of the transcendents 
at the critical points. Such a behavior must depend on 
two parameters, which are integration constants. The second problem, 
called {\it connection problem}, is 
to find the relation between the couples  
of parameters at different critical points. 
The method of isomonodromic deformations developed in \cite{JMU} \cite{JM1} 
was applied  to the Painlev\'e 6  equation in \cite{Jimbo},  to solve such problems for a 
class of solutions of PVI with  generic values of the 
parameters.  The non-generic case $\beta=\gamma=1-2\delta=0$ 
 is studied in \cite{DM} \cite{M} \cite{guz1} 
 for its applications to topological field theory. Studies on the critical behavior can be also found in \cite{Sh}. 

In the present paper we show that the elliptic representation is a valuable tool to study the critical behavior of the Painlev\'e 6 transcendents. We 
obtain results which include  the results of \cite{Jimbo} 
\cite{DM}  and extend the class of solutions to which they apply. On the other hand, we 
needed to use  the isomonodromic deformation theory to solve the connection problem, to be 
formulated below, 
 for the  elliptic representation.  

\vskip 0.2 cm

The elliptic representation was introduced by  R. Fuchs in \cite{fuchs}. Let    
$${\cal L}
  :=
 x(1-x) {d^2\over dx^2}+(1-2x) {d\over dx} -{1\over 4}.  
$$ 
be a linear differential operator and 
 let $\wp(z;\omega_1,\omega_2)$ be 
 the Weierstrass elliptic function of the independent variable $z\in 
{\bf P}^1$, with 
 {\it half-periods } $\omega_1$, $\omega_2$. Let us consider the following independent solutions of the  
{\it hyper-geometric equation} $
{\cal L} \omega=0$:
$$
\omega_1(x):={\pi \over 2} F\left({1\over 2},{1\over 2},1;x\right), ~~~~
\omega_2(x):=i{ \pi \over 2}F\left({1\over 2},{1\over 2},1;1-x\right), 
$$
where $F\left({1\over 2},{1\over 2},1;x\right)$ is the standard notation for 
the hyper-geometric function. R.Fuchs proved that the Painlev\'e 6 equation is equivalent to the 
following differential equation for a new function $u(x)$:
$$
 {\cal L}(u)= {1\over 2 x (1-x)} ~
{\partial \over \partial u} \left\{ 2\alpha\left[
\wp\left({u\over 2};\omega_1,\omega_2\right)+{1+x\over 3}\right]
-2\beta {x\over 
\wp\left({u\over 2};\omega_1,\omega_2\right)+{1+x\over 3} }
+\right.
$$
\be
 \left. +2\gamma{1-x \over \wp\left({u\over 2};\omega_1,\omega_2\right)+{x-2\over 3}}
+(1-2\delta){x(1-x)\over \wp\left({u\over 2};\omega_1,\omega_2\right)+
{1-2x\over 3}}
\right\}
\label{introELLITTICA}
\ee
 The connection to Painlev\'e 6 is given by the following representation of the transcendents: 
$$
   y(x)= \wp \left({u(x)\over 2};\omega_1(x),\omega_2(x) \right)
+{1+x\over 3}. 
$$ 
 We review these facts in section \ref{The Elliptic Representation}.

\vskip 0.2 cm 
 The algebraic-geometrical properties of the elliptic representations  
where studied   
in \cite{Manin2}. 
Nevertheless, the  analytic properties of the function $u(x)$ have
not been 
 studied so far, except for some special cases. The most simple case is
  $\alpha=\beta = \gamma=1-2\delta=0$. The function $u(x)$ is a 
linear combination of $\omega_1$ and $\omega_2$. This case was well known 
to Picard 
 \cite{Picard}, and the critical behavior  was 
studied in \cite{M}.   
A more general case was studied in \cite{guz1}, for  
$\beta = \gamma=1-2\delta=0$ and $\alpha$ any complex number. 
The motivation of \cite{guz1}  was that this case 
 is equivalent to
the WDVV equations of associativity in 2-D topological field theory introduced
by Witten \cite{Witten}, Dijkgraaf, Verlinde E., Verelinde H. \cite{DVV} and it has applications to Frobenius manifolds \cite{Dub1} and quantum cohomology \cite{Manin}.


In this paper, we study the analytic properties of $u(x)$  and 
 we compute the critical behavior of the transcendents  for {\it any} value of 
$\alpha,\beta,\gamma,\delta$; moreover,  we solve 
the connection problem  in elliptic representation for generic 
$\alpha,\beta,\gamma,\delta$ and for $\beta = \gamma=1-2\delta=0$.

\subsection{Our results}


\subsubsection{ Local Representation}\label{theorem1}

The equation  ${\cal L}(u)=0$ has a general solution $ 
   u_0= 2\nu_1 \omega_1 + 2 \nu_2 \omega_2$, $\nu_1,\nu_2\in{\bf C}
$. 
 We look for a solution of (\ref{introELLITTICA}) of the 
form $u(x) = 2\nu_1 \omega_1 +2\nu_2 \omega_2(x) + 2v(x)$, 
where $v(x)$ is a perturbation of $u_0$. 
 Let ${\bf C}_0:=
{\bf C}\backslash \{0\}$, $\widetilde{{\bf C}_0}$ the universal covering
 and let $0<r<1$. We define the domains 
 \be
{\cal D}(r;\nu_1,\nu_2):= 
\left\{ x\in \widetilde{{\bf C}_0}~ \hbox{ such that }
  |x|<r, \left|{e^{-i\pi \nu_1}\over 16^{1-\nu_2}} x^{1-\nu_2} 
\right|<r,
\left| 
{e^{i\pi \nu_1} \over 16^{\nu_2}} x^{\nu_2}\right|<r \right\}
\label{DOMINOO1talk}
\ee
\be
{\cal D}_0(r):= 
\left\{ x\in \widetilde{{\bf C}_0}~ \hbox{ such that }
  |x|<r \right\}
\label{DOMINOO3}
\ee
We observe that the translations  $\nu_i \mapsto \nu_i+2N_i$, $i=1,2$, $N_i \in {\bf Z}$ do not change a transcendent in the elliptic representation
$$
 y(x) = \wp \left(\nu_1\omega_1(x)+\nu_2\omega_2(x) +v(x); \omega_1(x),\omega_2(x)
 \right)+{1+x\over 3}.
$$
This is a consequence of  the periodicity of 
the $\wp$-function. Therefore, one  can  take  
$0\leq \Re  \nu_i <2$, $i=1,2$. Nevertheless,
 we don't need to suppose such a 
range explicitly. Only in the case 
 $\Im \nu_2=0$ we need to  suppose that $0\leq \nu_2<2$.  Finally, let us 
introduce the following expansion:
\be 
v(x;\nu_1,\nu_2):= \sum_{n\geq 1} a_n x^n + \sum_{n\geq 0, m\geq 1} b_{nm} x^n 
\left[ e^{-i\pi \nu_1} \left({x\over 16}\right)^{1-\nu_2}\right]^m 
+\sum_{n\geq 0,m\geq 1} c_{nm} x^n 
\left[ e^{i\pi \nu_1} \left({x\over 16}\right)^{\nu_2}\right]^m 
\label{expatalk}
\ee

\vskip 0.3 cm 
\noindent
{\bf Theorem 1:} { \it  Let $\nu_1$, $\nu_2$ be two complex numbers. 

\vskip 0.2 cm 
{\bf I)} For any 
complex  $\nu_1$, $\nu_2$ such that $ 
\Im \nu_2\neq 0
$ 
there exist a  positive number $r<1$  and a  transcendent 
$$
 y(x) = \wp \Bigl(\nu_1\omega_1(x)+\nu_2\omega_2(x) +v(x;\nu_1,\nu_2);~ \omega_1(x),\omega_2(x)
 \Bigr)+{1+x\over 3}
$$
such that $v(x;\nu_1,\nu_2)$ is holomorphic in the domain $ 
{\cal D}(r;\nu_1,\nu_2)  
$ and it is given by  the expansion (\ref{expatalk}) which is  convergent in 
 $ {\cal D}(r;\nu_1,\nu_2)$. 
The coefficients $a_n$, 
$b_{nm}$, $c_{nm}$, $i=1,2$, are certain 
rational functions of $\nu_2$. Moreover, 
there exists a positive constant $M(\nu_2)$ such that 
  \be
  |v(x;\nu_1,\nu_2) | \leq M(\nu_2) \left(|x|+\left| e^{-i\pi \nu_1} \left({x\over 16}
\right)^{1-\nu_2}\right| + \left|  e^{i\pi \nu_1} 
\left({x\over 16}\right)^{\nu_2} \right| 
 \right)~~~ \hbox{ in } {\cal D}(r;\nu_1,\nu_2)
\label{bundooo1}
\ee

\vskip 0.2 cm  

{\bf II)} For any complex $\nu_1$ and  real $\nu_2$, with the constraint  $0<\nu_2<1$ or $1<\nu_2<2$, 
there exists a positive  $r<1$ and a  transcendent 
$$
 y(x) = \wp \Bigl(\nu_1\omega_1(x)+\nu_2\omega_2(x) +v(x;\nu_1,\nu_2);~ \omega_1(x),\omega_2(x)
 \Bigr)+{1+x\over 3}, ~~\hbox{ if } 0<\nu_2<1 
$$
or 
$$
 y(x) = \wp \Bigl(\nu_1\omega_1(x)+\nu_2\omega_2(x) +v(x;-\nu_1,2-\nu_2);~ \omega_1(x),\omega_2(x)
 \Bigr)+{1+x\over 3}, ~~\hbox{ if } 1<\nu_2<2 
$$
such that $v(x;\nu_1,\nu_2)$ and $v(x;-\nu_1,2-\nu_2)$ are 
 holomorphic in ${\cal D}_0(r)$,  
with  convergent expansion  (\ref{expatalk}) and bound (\ref{bundooo1}) (for $1<
\nu_2<2$ substitute $\nu_1 \mapsto -\nu_1$, $\nu_2 \mapsto 2-\nu_2$). 

}

\vskip 0.3 cm
 Note that in the theorem 
$$ 
\nu_2\neq 0,1
$$ 
We stress  that in case {\bf II)}, if $\nu_2$ 
 is greater that 2 or less then 0, we can always make a 
translation $\nu_2\mapsto \nu_2+2N$ to obtain $0< 
\nu_2 <2$ (on the other hand, 
if $-2N < \nu_2< 2-2N$, the formulae of case {\bf II)} hold 
with the substitution $\nu_2 \mapsto \nu_2+2N$). 

\vskip 0.3 cm 

\noindent
{\it Observation 1:} 
 As a consequence of the theorem, for any $N\in {\bf Z}$ and  for any complex 
$\nu_1,\nu_2$ such that  $\Im \nu_2\neq 0$, there exists $r_N<1$ and a transcendent 
 $y(x) = \wp \Bigl(\nu_1\omega_1(x)+[\nu_2+2N]\omega_2(x) +v(x;\nu_1,\nu_2+2N);~ \omega_1(x),\omega_2(x)
 \Bigr)+{1+x\over 3}$ in  ${\cal D}(r;\nu_1,\nu_2+2N)$. By periodicity of the 
$\wp$-function we re-write the transcendent as follows: 
$$
y(x) = \wp \Bigl(\nu_1\omega_1(x)+\nu_2\omega_2(x) +v(x;\nu_1,\nu_2+2N);~ \omega_1(x),\omega_2(x)
 \Bigr)+{1+x\over 3}~~\hbox{ in } {\cal D}(r;\nu_1,\nu_2+2N).
$$
Moreover, we will show in section \ref{twoquestions} that if a transcendent has the  
 elliptic representation 
$$
y(x)= 
 \wp \Bigl(\nu_1\omega_1(x)+\nu_2\omega_2(x) +v(x;\nu_1,\nu_2);~ \omega_1(x),\omega_2(x)
 \Bigr)+{1+x\over 3}
$$ 
in ${\cal D}(r, \nu_1,\nu_2)$ for some $\nu_1, \nu_2$, $\Im \nu_2 \neq 0$, then for any integer 
$N$ there exists $\nu_1^{\prime}$ (depending on $\nu_1$, $\nu_2$ and $N$) such that the 
transcendent has also the representation 
$$y(x)= \wp \Bigl(\nu_1^{\prime}\omega_1(x)+\nu_2\omega_2(x) +v(x;\nu_1,\nu_2+2N);~ \omega_1(x),\omega_2(x)
 \Bigr)+{1+x\over 3}
$$
 in ${\cal D}(r, \nu_1^{\prime} ,\nu_2+2N)$.

\vskip 0.3 cm 

\noindent
{\it Observation 2:} 
 Another consequence of the theorem is that for  any complex  $\nu_1,\nu_2$ 
such that 
 $\Im \nu_2\neq 0$ there exists $ y(x) =  
\wp \Bigl(-\nu_1\omega_1(x)+[2-\nu_2]~\omega_2(x) +v(x;-\nu_1,2-\nu_2);~ \omega_1(x),\omega_2(x)
 \Bigr)+{1+x\over 3}$. Again we use the fact that  the 
$\wp$-function  is periodic w.r.t. $2\omega_2$ and it is an even function. Therefore the transcendent becomes 
$$
 y(x) =  \wp \Bigl(\nu_1\omega_1(x)+\nu_2\omega_2(x) -v(x;-\nu_1,2-\nu_2);~ \omega_1(x),\omega_2(x)
 \Bigr)+{1+x\over 3},~~\hbox{ in } {\cal D}(r;-\nu_1,2-\nu_2)
$$
 Note that the series $-v(x;-\nu_1,2-\nu_2)$ is of the form 
$$
\sum_{n\geq 1} a_n x^n + \sum_{n\geq 0, m\geq 1} b_{nm} x^n 
\left[ e^{-i\pi \nu_1} \left({x\over 16}\right)^{2-\nu_2}\right]^m 
+\sum_{n\geq 0,m\geq 1} c_{nm} x^n 
\left[ e^{i\pi \nu_1} \left({x\over 16}\right)^{\nu_2-1}\right]^m 
$$
where we have re-named the constants $a_n, b_{nm}, c_{nm}$.

\begin{figure}
\epsfxsize=12cm
\centerline{\epsffile{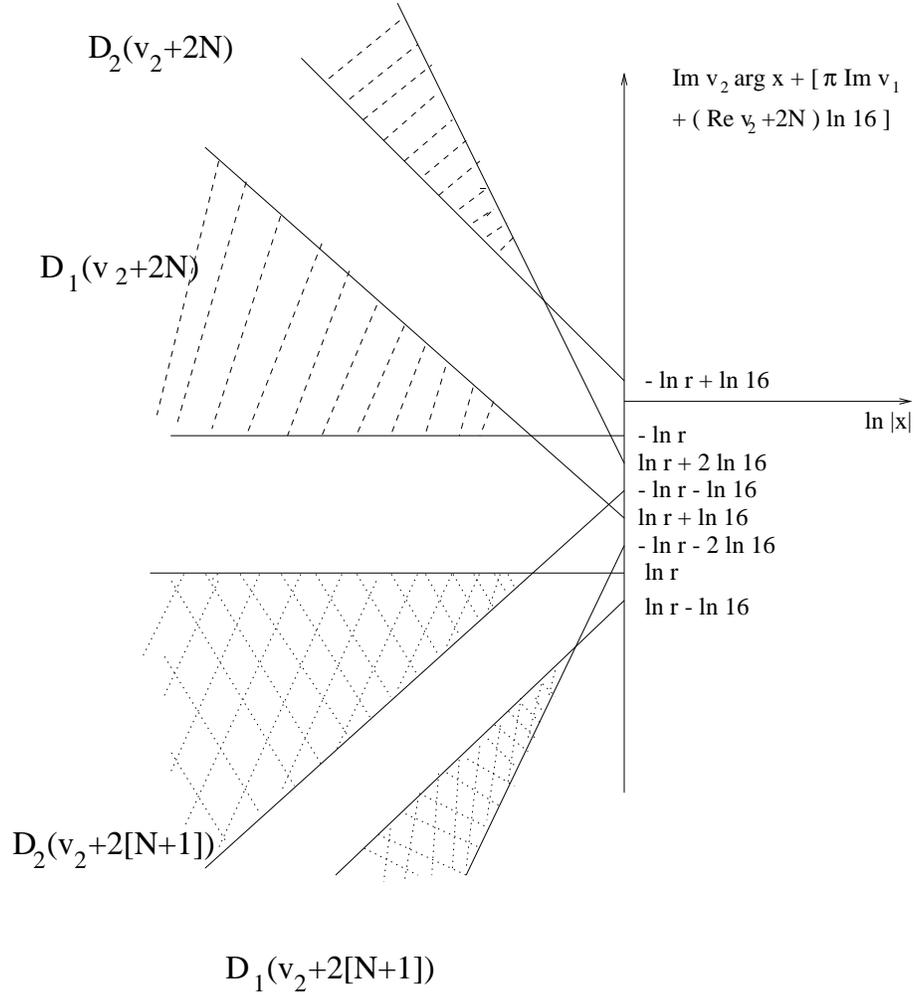}}
\caption{The domains ${\cal D}_1(r;\nu_1,\nu_2+2N):={\cal D}(r;\nu_1,\nu_2+2N)$, 
${\cal D}_2(r;\nu_1,\nu_2+2N):={\cal D}(r;-\nu_1,2-\nu_2-2N)$ 
and  ${\cal D}_1(r;\nu_1,\nu_2+2[N+1])$, 
${\cal D}_2(r;\nu_1,\nu_2+2[N+1])$ for  arbitrarily 
 fixed values of $\nu_1$, $\nu_2$, $N$. They are represented in the plane  
$(\ln|x|,~\Im \nu_2 \arg x +[\pi \Im \nu_1 +(\Re \nu_2 +2N)\ln 16])$. }
\label{figur1}
\end{figure}

\vskip 0.3 cm
The domain  ${\cal D}(r_N;\nu_1,\nu_2+2N)$ can be written as 
follows: 
$$ 
  (\Re \nu_2+2N) \ln{|x|\over 16} -\pi \Im \nu_1
-\ln r_N < \Im \nu_2
   \arg x <
$$
$$ < (\Re \nu_2-1 +2N)\ln{|x|\over 16}-\pi \Im \nu_1
+\ln r_N,~~~~|x|<r_N
$$
 Therefore 
the domain  ${\cal D}(r_N,-\nu_1,2-\nu_2-2N)$ is 
$$   (\Re \nu_2-1+2N) \ln{|x|\over 16}-\pi \Im \nu_1
-\ln r_N < \Im \nu_2
   \arg x < $$
$$ 
< (\Re \nu_2 -2+2N)\ln{|x|\over 16}-\pi \Im \nu_1
+\ln r_N, ~~~~|x|<r_N
$$
We can draw their picture in the $(\ln|x|,\Im\nu_2 \arg x)$-plane. See figure 
\ref{figur1}.

\subsubsection{Critical Behavior}

 It is possible to compute the critical behavior for $x\to 0$ of a transcendent of Theorem 1. For simplicity, we consider $x\to 0$  along  the  paths defined below. 
Let $\Im \nu_2 \neq 0$ and ${\cal V}\in {\bf C}$.  We define the following 
family of 
paths joining a point $x_0\in {\cal D}(r;\nu_1,\nu_2)$  to $x=0$
\be
  \arg x = \arg x_0 + {\Re \nu_2 - {\cal V}  \over \Im \nu_2} \ln{|x|\over 
|x_0|},~~~~0\leq {\cal V}\leq1
\label{patheli}
\ee
 The paths are contained 
in
   ${\cal D}(r;\nu_1,\nu_2)$.  
If $\Im \nu_2=0$ any regular path  contained in 
${\cal D}_0(r)$ can be considered. 

\vskip 0.3 cm 
\noindent
{\bf Theorem 2:} {\it  Let $\nu_1$, $\nu_2$ be given.  

\vskip 0.2 cm 
If $\Im \nu_2 \neq 0$, the critical behavior of the transcendent 
$y(x)=\wp(\nu_1\omega_1+
\nu_2\omega_2+v(x;\nu_1,\nu_2);\omega_1,\omega_2)+(1+x)/3$   when $x\to 0$ along the path (\ref{patheli}) is: 

\vskip 0.2 cm 
\noindent
For $0<{\cal V}<1$:
\be
y(x)= -{1\over 4} \left[{e^{i\pi \nu_1} \over 16^{\nu_2-1}}\right] ~
          x^{\nu_2} ~\left(1+ O(|x^{\nu_2}|+|x^{1-\nu_2}|)\right).
\label{critA}
\ee

\vskip 0.2 cm 
\noindent
For ${\cal V}=0$:
\be
y(x)=  \left[
       {x\over 2} +\sin^{-2}
       \left( 
  -i{\nu_2\over 2} \ln {x\over 16} +{\pi \nu_1 \over 2} + 
  \sum_{m\geq 1} c_{0m} 
\left[
e^{i\pi \nu_1} 
\left(
                     {x\over 16}
\right)^{\nu_2}
\right]^m
\right)           
       \right] ~\bigl(1+O(x)\bigr).
\label{critC}
\ee

\vskip 0.2 cm 
\noindent
For ${\cal V}=1$:
\be
y(x) =
x~\sin^2
                                                             \left(
i{1-\nu_2\over 2}\ln{x\over 16}+{\pi \nu_1 \over 2}+ 
\sum_{m\geq 1} b_{0m}
                                 \left[
e^{-i\pi \nu_1}
                                 \left(
{x\over 16}
                                        \right)^{1-\nu_2}
                                       \right]^m
                                                              \right)~(1+O(x)).
\label{critB}
\ee

\vskip 0.2 cm 
\noindent
For $\nu_2$  real we have two cases. 
For $0<\nu_2<1$, the transcendent $y(x)=\wp(\nu_1\omega_1+
\nu_2\omega_2+v(x;\nu_1,\nu_2);\omega_1,\omega_2)+(1+x)/3$ defined in 
${\cal D}_0(r)$  has behavior 
\be
 y(x)= -{1\over 4} \left[{e^{i\pi \nu_1} \over 16^{\nu_2-1}}\right] ~
          x^{\nu_2} ~\left(1+ O(|x^{\nu_2}|+|x^{1-\nu_2}|)\right),
~~~0<\nu_2<1
\label{critE}
\ee
For $1<\nu_2<2$, the transcendent $y(x)=\wp(\nu_1\omega_1+
\nu_2\omega_2+v(x;-\nu_1,2-\nu_2);\omega_1,\omega_2)+(1+x)/3$ defined in 
${\cal D}_0(r)$  has behavior 
\be
y(x)= -{1\over 4} \left[{e^{i\pi \nu_1} \over 16^{\nu_2-1}}\right]^{-1} ~
          x^{2-\nu_2} ~\left(1+ O(|x^{2-\nu_2}|+|x^{\nu_2-1}|)\right),
~~~ 1<\nu_2<2
\label{critEE}
\ee
}

\subsubsection{The Critical Points $x=1,\infty$}

 Theorem 1 and 2 deal with the point $x=0$. We now turn to the other critical points. 
Let us define  $\omega_1^{(0)}:=\omega_1$, $\omega_2^{(0)}:=\omega_2$; $\omega_1^{(1)}:=\omega_2$, $\omega_2^{(1)}:=\omega_1$ and 
 $\omega_1^{(\infty)}:=\omega_1+\omega_2$, $\omega_2^{(\infty)}:=\omega_2$. 
We construct solutions 
$$ 
  {u(x)\over 2} = \nu_1^{(1)}\omega_1^{(1)}(x)
+\nu_2^{(1)}\omega_2^{(1)}(x) +
v^{(1)}(x)
$$
in a neighborhood of $x=1$, and solutions 
$$ 
  {u(x)\over 2} = \nu_1^{(\infty)}\omega_1^{(\infty)}(x)
+\nu_2^{(\infty)}\omega_2^{(\infty)}(x) +
v^{(\infty)}(x)
$$
in a neighborhood of $x=\infty$. Let 
${\bf P}^1_i
:={\bf P}^1 \backslash \{i\}$, $i=1,\infty$, and let $\widetilde{{\bf P}^1_i}$ be 
the universal covering.

\vskip 0.3 cm 
\noindent
{\bf Theorem 3}: (In a Neighborhood of $x=1$). { \it For any 
complex $\nu^{(1)}_1$ $\nu_2^{(1)}$ such that 
$\Im \nu_2^{(1)} \neq 0$ 
 there exists  a  transcendent $y(x)= 
\wp \bigl( \nu_1^{(1)} \omega_1^{(1)}+ \nu_2^{(1)} \omega_2^{(1)}+
v^{(1)}(x;\nu_1^{(1)},\nu_2^{(1)});\omega_1^{(1)},\omega_2^{(1)}\bigr)+{1+x\over 3}$ 
 such that $v^{(1)}(x)$ is holomorphic in the domain
$$
{\cal D}(r;\nu_1^{(1)},\nu_2^{(1)}) := 
\left\{
x\in \widetilde{{\bf P}^1_1}~\hbox{ such that }~|1-x|<r, ~\left|e^{i\pi \nu_1^{(1)}} \left( {1-x\over 16}\right)^{1-\nu^{(1)}_2}\right|<r,\right.
$$
$$ \left. \left|e^{-i\pi \nu_1^{(1)}} 
\left( {1-x\over 16}\right)^{\nu^{(1)}_2} \right|<r 
\right\}
$$
where it  has the convergent 
 expansion: 
$$
  v^{(1)}(x;\nu_1^{(1)},\nu_2^{(1)}) = 
\sum_{n\geq 1} a_n (1-x)^n+\sum_{n\geq 0,m\geq 1} 
 b_{nm} (1-x)^n \left[e^{i\pi \nu_1^{(1)}} \left( {1-x\over 16}\right)^{1-\nu^{(1)}_2} \right]^m +$$
\be +\sum_{n\geq 0,m\geq 1} 
 c_{nm} (1-x)^n \left[e^{-i\pi \nu_1^{(1)}} \left( {1-x\over 16}\right)^{\nu^{(1)}_2} \right]^m 
\label{stufo1}
\ee

 For any complex $\nu_1^{(1)}$ and real $\nu^{(1)}_2$ with the constraint 
$
0<\nu_2^{(1)}<1$, there exists a sufficiently small $r$ and a transcendent  $y(x)= 
\wp ( \nu_1^{(1)} \omega_1^{(1)}+ \nu_2^{(1)} 
\omega_2^{(1)}+v^{(1)}(x;\nu_1^{(1)},\nu_2^{(1)}))+{1+x\over 3}$ [we omit half-periods for semplicity] such that 
$v^{(1)}
(x;\nu_1^{(1)},\nu_2^{(1)})$ is holomorphic in  
$$
{\cal D}_0(r):= \left\{
x\in \widetilde{{\bf P}^1_1}~\hbox{ such that }~|1-x|<r \right\}
$$
 where it   has convergent expansion 
(\ref{stufo1}). 

 For any complex $\nu_1^{(1)}$ and real $\nu^{(1)}_2$ with the constraint 
$
1<\nu_2^{(1)}<2$,  there exists a sufficiently small $r$ and a transcendent 
$y(x)= 
\wp ( \nu_1^{(1)} \omega_1^{(1)}+ \nu_2^{(1)} 
\omega_2^{(1)}+v^{(1)}(x;-\nu_1^{(1)},2-\nu_2^{(1)}))+{1+x\over 3}$ such that 
$v^{(1)}(x;-\nu_1^{(1)},2-\nu_2^{(1)}))$ is holomorphic in  $
{\cal D}_0(r)$,  
 where it   has convergent expansion 
(\ref{stufo1}) with the substitution $(\nu_1^{(1)},\nu_2^{(1)})\mapsto 
(-\nu_1^{(1)},2-\nu_2^{(1)})$.

}

\vskip 0.3 cm 
\noindent
 (In a Neighborhood of $x=\infty$). 
{ \it For any complex 
  $\nu^{(\infty)}_1$ $\nu_2^{(\infty)}$ such that $\Im \nu_2^{(\infty)}\neq 0$ 
     there exists a  transcendent 
$y(x)= 
\wp \bigl( \nu_1^{(\infty)} \omega_1^{(\infty)}+ \nu_2^{(\infty)} 
\omega_2^{(\infty)}+v^{(\infty)}(x;\nu_1^{(\infty)},\nu_2^{(\infty)}); \omega_1^{(\infty)}, \omega_2^{(\infty)}\bigr)+{1+x\over 3}$ 
 such that $v^{(\infty)}(x;\nu_1^{(\infty)},\nu_2^{(\infty)})$ 
is holomorphic  in 
$$
{\cal D}(r;\nu_1^{(\infty)},\nu_2^{(\infty)}) := 
\left\{
x\in \widetilde{{\bf P}^1_{\infty}}~\hbox{ such that }~|x^{-1}|<r, ~\left|
e^{-i\pi \nu_1^{(\infty)}} \left( {16\over x}\right)^{1-\nu^{(\infty)}_2
}\right|<r,\right.$$
$$
\left.
 \left|e^{i\pi \nu_1^{(\infty)}} \left( {16\over 
x}\right)^{\nu^{(\infty)}_2} \right|<r 
\right\}
$$ where it  has the 
 convergent expansion: 
$$
 x^{1\over 2} 
v^{(\infty)}(x;\nu_1^{(\infty)},\nu_2^{(\infty)}) = \sum_{n\geq 1} a_n \left({1\over x}\right)^n
+\sum_{n\geq0,m\geq 1} 
 b_{nm} \left({1\over x}\right)^n \left[e^{-i\pi \nu_1^{(\infty)}} 
\left( {16\over x}\right)^{1-\nu^{(\infty)}_2
} \right]^m +$$
\be
+\sum_{n\geq 0,m\geq 1} 
 c_{nm} \left({1\over x}\right)^n \left[e^{i\pi \nu_1^{(\infty)}} \left( {16\over x}\right)^{\nu^{(\infty)}_2} \right]^m 
\label{stufo11}
\ee

 For any $\nu^{(\infty)}_1$ and real  $\nu^{(\infty)}_2$ with  the constraint 
$0<\nu_2^{(\infty)}<1$, there exists a sufficiently small $r$ and a transcendent  $y(x)= 
\wp ( \nu_1^{(\infty)} \omega_1^{(\infty)}+ \nu_2^{(\infty)} 
\omega_2^{(\infty)}+v^{(\infty)}(x;\nu_1^{(\infty)},\nu_2^{(\infty)}))+{1+x\over 3}$ such that $
v^{(\infty)}(x;\nu_1^{(\infty)},\nu_2^{(\infty)})$ is holomorphic in  
$$
{\cal D}_0(r):= \left\{
x\in  \widetilde{{\bf P}_{\infty}}~\hbox{ such that }
~|x^{-1}|<r\right\}
$$
where  it  has convergent expansion  
(\ref{stufo11}).    
For any $\nu^{(\infty)}_1$ and real  $\nu^{(\infty)}_2$ with  the constraint 
$1<\nu_2^{(\infty)}<2$ there exists a sufficiently small $r$ and a transcendent $y(x)= 
\wp ( \nu_1^{(\infty)} \omega_1^{(\infty)}+ \nu_2^{(\infty)} 
\omega_2^{(\infty)}+v^{(\infty)}(x;-\nu_1^{(\infty)},2-\nu_2^{(\infty)}))+{1+x\over 3}$ such that 
$v^{(\infty)}(x;-\nu_1^{(\infty)},2-\nu_2^{(\infty)})$ is holomorphic in $
{\cal D}_0(r)$ 
where  it  has convergent expansion  
(\ref{stufo11}) with the substitution $(\nu_1^{(\infty)},\nu_2^{(\infty)})\mapsto (-\nu_1^{(\infty)},2-\nu_2^{(\infty)})$. 
}

\vskip 0.3 cm 
\noindent
{\it Note:} We have used the notations ${\cal D}(r;\nu_1,\nu_2)$ for the 
domains both at $x=0$, $x=1$ and $x=\infty$. We believe this will not confuse 
the reader, because it is always clear which is the critical point we are 
considering.  
Note that $\nu_1^{(1)}$ comes with sign changed w.r.t. $\nu_1$ at $x=0$; this is due to the definition of $\omega_1^{(1)}$.

\vskip 0.3 cm

In section \ref{theorem3} we also compute the critical behaviors 
for $x\to 1 $ and $x\to \infty$.

\subsubsection{Connection Problem}

 The elliptic representation allows us to obtained 
 detailed information about the critical behavior of the
Painlev\'e transcendents.  On the other hand, the local analysis does 
 not solve the
 {\it connection problem}. This is the problem of determining the 
critical behavior of a given transcendent at both $x=0$, $x=1$ and 
$x=\infty$.  In our framework, we ask if a transcendent may have, at 
the same time, three representations
$$
  y(x) =\wp (\nu_1^{(0)}\omega_1^{(0)}+\nu_2^{(0)}\omega_2^{(0)}+v^{(0)}) +
{1+x\over 3}
$$
$$ = \wp (\nu_1^{(1)}\omega_1^{(1)}+\nu_2^{(1)}\omega_2^{(1)}+v^{(1)}) +
{1+x\over 3}
$$
$$
  = \wp (\nu_1^{(\infty)}\omega_1^{(\infty)}+
\nu_2^{(\infty)}\omega_2^{(\infty)}+v^{(\infty)}) +
{1+x\over 3}
$$
in the domains of Theorems 1 and 3. 
Moreover, we look for formulae which connect the three 
couples 
of parameters $(\nu_1^{(0)},\nu_2^{(0)})$, $(\nu_1^{(1)},\nu_2^{(1)})$, 
 $(\nu_1^{(\infty)},\nu_2^{(\infty)})$.

\vskip 0.2 cm
The connection problem may be solved using the method of isomonodromic deformations. Works on this problem are \cite{Jimbo}, \cite{DM} and \cite{guz1}. 
 The PVI is the isomonodromy deformation equation of  a Fuchsian system of
 differential equations
$$
   {dY\over dz}=\left[ {A_0(x)\over z}+{A_x(x) \over z-x}+{A_1(x)\over
z-1}\right] Y$$
The  $2\times 2$ matrices  $A_i(x)$ ($i=0,x,1$ are labels) depend  on $x$ in such a way that the monodromy of a fundamental solution $Y(z,x)$ does not change for  small deformations of $x$. They depend on the parameters $\alpha,\beta,\gamma,\delta$ of PVI as follows:
$$
A_0(x)+A_1(x)+A_x(x) = -{1\over 2} \pmatrix{ \theta_{\infty} & 0 \cr 
                                                 0       & -\theta_{\infty}\cr
}  , ~~~~ 
  \hbox{  eigenvalues of } A_i(x)=\pm {1\over 2} \theta_i,~~~~i=0,1,x
$$
$$ 
    \alpha= {1\over 2} (\theta_{\infty} -1)^2,
~~\beta=-{1\over 2} \theta_0^2, 
~~ \gamma={1\over 2} \theta_1^2,
~~\delta={1\over 2} (1-\theta_x^2) 
$$

In section \ref{twoquestions} we solve the connection problem for the elliptic 
representation for generic values of $\alpha$, $\beta$, $\gamma$, $\delta$ . More precisely, by {\it generic case} we mean:  
\be
   \nu_2^{(i)}, ~\theta_0,~\theta_x,~\theta_1,~\theta_{\infty}~\not \in {\bf Z};
~~~~
{\pm1\pm \nu_2^{(i)} \pm \theta_1 \pm \theta_{\infty}\over 2},~
{\pm1\pm\nu_2^{(i)}  \pm \theta_0 \pm \theta_x \over 2} ~\not \in {\bf Z}
\label{generisspecies}
\ee
The signs $\pm$ vary independently. 
This is a technical condition which can be abandoned (except for $\nu_2^{(i)}\not \in 
{\bf Z}$) at the price of making the computations more 
complicated. For example, the non-generic case $\beta=\gamma=1-2\delta=0$ and $\alpha$ any complex number was analyzed in \cite{guz1}  for its relevant applications to Frobenius 
manifolds and quantum cohomology. We will review it in the paper. 

\vskip 0.3 cm
 To summarize the results for the generic case, we first observe that the critical behaviors provided by the elliptic representations along regular paths (except special directions for ${\cal V}=0,1$, see Theorem 2 and section \ref{HIEI-ZAN}) at $x=0$, $x=1$ and $x=\infty$ respectively are  
\be
y(x) = 
 a^{(0)} x^{\nu_2^{(0)}}(1+ \hbox{ higher orders in $x$}), ~~~x\to 0  
\label{LOc1introduzione}
\ee
\be
y(x)=
 1- a^{(1)} (1-x)^{\nu_2^{(1)}}(1+\hbox{ higher orders in $(1-x)$}),
~~~x\to 1
\label{LOc2introduzione}
\ee
\be
y(x) = 
a^{(\infty)} x^{1-\nu_2^{(\infty)}}(1+\hbox{ higher orders in $x^{-1}$}),~~~x\to 
\infty 
\label{LOc3introduzione}
\ee
and the parameters $\nu_1^{(i)}$ are given by 
$$
    e^{i\pi \nu_1^{(0)}}= -4a^{(0)}~16^{\nu_2^{(0)}-1},
~~~e^{-i\pi \nu_1^{(1)}}= -4a^{(1)}~16^{\nu_2^{(1)}-1},
~~~ e^{i\pi \nu_1^{(\infty)}}= -4a^{(\infty)}~16^{\nu_2^{(\infty)}-1}
$$
If $\nu_2^{(i)}$ is real, the behavior  is as above when $0<\nu_2^{(i)}<1$. Otherwise, when $1<\nu_2^{(i)}<2$ it is:
\be
y(x) = 
 a^{(0)} x^{2-\nu_2^{(0)}}(1+ \hbox{ higher orders in $x$}), ~~~x\to 0  
\label{LOc1introduzioNE}
\ee
\be
y(x)=
 1- a^{(1)} (1-x)^{2-\nu_2^{(1)}}(1+\hbox{ higher orders in $(1-x)$}),
~~~x\to 1
\label{LOc2introduzioNE}
\ee
\be
y(x) = 
a^{(\infty)} x^{\nu_2^{(\infty)}-1}(1+\hbox{ higher orders in $x^{-1}$}),~~~x\to 
\infty 
\label{LOc3introduzioNE}
\ee
with
$$
    e^{-i\pi \nu_1^{(0)}}= -4a^{(0)}~16^{1-\nu_2^{(0)}},
~~~e^{i\pi \nu_1^{(1)}}= -4a^{(1)}~16^{1-\nu_2^{(1)}},
~~~ e^{-i\pi \nu_1^{(\infty)}}= -4a^{(\infty)}~16^{1-\nu_2^{(\infty)}}
$$
Note that the ambiguity $\nu_1^{(i)}\mapsto \nu_1^{(i)}+2k$, $k$ integer, is 
 natural, because $v^{(i)}(x)$ does not change and the $\wp$-function is periodic. 

\vskip 0.2 cm

Let  $M_0$, $M_1$, $M_x$ be the monodromy matrices 
 at $z=0,1,x$, for a given basis in the fundamental 
group of ${\bf P}^1\backslash\{0,1,x,\infty\}$. Such basis is chosen as 
in figure \ref{figure3}.

 If 
$$
\theta_0,~\theta_x,~\theta_1,~\theta_{\infty}~\not \in {\bf Z}
$$
 there is a one to one correspondence between a given choice of   {\it monodromy data} 
$\theta_0$, $\theta_x$, $\theta_1$, $\theta_{\infty}$, $\hbox{tr}(M_0M_x)$,  $\hbox{tr}(M_0M_1)$,
 $\hbox{tr}(M_1M_x)$ and a transcendent $y(x)$. Namely:
\be
y(x)=
y\bigl(x;~ \theta_0, \theta_x,\theta_1,\theta_{\infty}, \hbox{tr}(M_0M_x),  \hbox{tr}(M_0M_1),
 \hbox{tr}(M_1M_x)~ \bigr)
\label{chineseramen}
\ee
We prove that such a transcendent has elliptic representations at $x=0,1,\infty$, provided 
that (\ref{generisspecies}) is satisfied.  The three sets of parameters $(\nu_1^{(i)}, \nu_2^{(i)})$, $i=0,1,\infty$ are functions of the monodromy 
data $\theta_0$, $\theta_x$, $\theta_1$, $\theta_{\infty}$, $\hbox{tr}(M_0M_x)$,  $\hbox{tr}(M_0M_1)$,
 $\hbox{tr}(M_1M_x)$. Namely, we show that
\be
   2\cos(\pi \nu_2^{(0)})=- \hbox{tr}(M_0M_x),~~~ 2\cos(\pi \nu_2^{(1)})=- \hbox{tr}(M_1M_x),~~~
 2\cos(\pi \nu_2^{(0)})=- \hbox{tr}(M_0M_1)
\label{pastasciutta}
\ee
$$
  a^{(i)}=a^{(i)}\bigl(\nu_2^{(i)};\theta_0, \theta_x, \theta_1, \theta_{\infty}, \hbox{tr}(M_0M_x),  \hbox{tr}(M_0M_1),
 \hbox{tr}(M_1M_x)\bigr), ~~~~i=0,1,\infty
$$
 The formulas of $a^{(i)}$ are quite long, so we do not write them here. They depend on the monodromy data through rational, trigonometric and $\Gamma$-functions. In particular, $\nu_2^{(i)}$ enters explicitly. The procedure for computing such formulae is given in the Appendix. We 
note that the condition $\nu_2^{(i)}\not \in {\bf Z}$ is equivalent to tr$(M_iM_j)\neq \pm 2$.

Conversely, we prove that a transcendent 
 $y(x)$ given by its elliptic representation close to one critical point is a transcendent (\ref{chineseramen}) for some monodromy data. This follows from the consideration 
that  the couple $(\nu_1^{(i)},
\nu_2^{(i)})$ is given at the critical point $x=i$, and    $\theta_0$, $\theta_x$, $\theta_1$, $\theta_{\infty}$ are fixed by the equation PVI we are considering.  
 From these data we 
 can  compute   $\hbox{tr}(M_0M_x)$,  $\hbox{tr}(M_1M_x)$,  $\hbox{tr}(M_0M_1)$. One of 
the traces is $- 2\cos(\pi \nu_2^{(i)})$, the others depend on 
$\nu_1^{(i)},\nu_2^{(i)}$,   $\theta_0$, $\theta_x$, $\theta_1$, $\theta_{\infty}$ through 
rational, trigonometric and $\Gamma$-functions. The formulae are rather long, so we refer 
the reader to the Appendix.  In this way the transcendent  (\ref{chineseramen}) is obtained. 
  From the monodromy data we compute the couples $(\nu_1^{(j)},
\nu_2^{(j)})$ at the other two critical points and we get the the elliptic representation 
of the initial transcendent at the other critical points. Therefore, the connection 
problem is solved.  

  Note  that if we start from the elliptic representation at one critical point, say for 
example $x=0$, then $\nu_1^{(0)},\nu_2^{(0)}$  are given. As explained above, we can 
 compute the monodromy data 
and from them we compute $\nu^{(j)}_2$ and $a^{(j)}$ (then $\nu_1^{(j)}$) at the other two 
critical 
points. As already observed, the ambiguity  $\nu_1^{(j)}\mapsto \nu_1^{(j)}+2k$  ($k$ integer) does not change the elliptic representation. On the other hand, 
the ambiguities $\nu_2^{(j)}\mapsto \nu_2^{(j)}+2N$ ($N$ integer),  $\nu_2^{(j)} \mapsto - \nu_2^{(j)}$ 
and the ambiguity in the choice $0\leq \Re \nu_2^{(j)}\leq 1$ or 
 $1\leq \Re \nu_2^{(j)}\leq 2$, which results from the cosines in (\ref{pastasciutta}), is 
due to the fact that the same transcendent has  
different elliptic representations in different domains (the choice of $\nu_2^{(j)}$ 
determines the representation and  the domain!).  

\vskip 0.2 cm
 To conclude the discussion of the generic case, some comments about our extension of previous 
 known results are in order. 
The critical behavior  for a class of solutions to  
the Painlev\'e 6 equation was found  by
Jimbo in  \cite{Jimbo}  
for  generic values of 
 $\alpha$, $\beta$, $\gamma$ $\delta$.     
 A transcendent in this class has behavior: 
\be
y(x)= a^{(0)} x^{1-\sigma^{(0)}}(1+O(|x|^{\delta})),~~~~x\to 0,
\label{loc1introduzione}
 \ee
\be
y(x)= 1-a^{(1)}(1-x)^{1-\sigma^{(1)}} (1+O(|1-x|^{\delta})),~~~~x\to 1,
\label{loc2introduzione}
 \ee
\be
y(x)= a^{(\infty)}
 x^{-\sigma^{(\infty)}}(1+O(|x|^{-\delta})),~~~~x\to \infty,
\label{loc3introduzione}
\ee
where $\delta$ is a small positive number, $a^{(i)}$ and
 $\sigma^{(i)}$ are complex numbers such that $a^{(i)}\neq 0$ and 
 \be
0\leq \Re \sigma^{(i)}<1.
\label{RESTOacasa}
\ee 
We remark that  $x$ converges
 to the critical points {\it inside a sector} with vertex
 on the corresponding critical point.   
The  connection problem, i.e.  the problem of 
finding  the relation among the three pairs $(\sigma^{(i)},a^{(i)})$,
$i=0,1,\infty$, was solved  in \cite{Jimbo} for the above class of transcendents 
  using the
isomonodromy deformations theory. Actually, a transcendent in the class above coincides 
with a transcendent (\ref{chineseramen}). 
In particular
$$ 
 2 \cos(\pi \sigma^{(0)})=\hbox{tr}(M_0M_x),~~~2 \cos(\pi \sigma^{(1)})=\hbox{tr}(M_1M_x),~~~2 \cos(\pi \sigma^{(\infty)})=\hbox{tr}(M_0M_1)
$$
and 
$$
  a^{(i)}=a^{(i)}\bigl(\sigma^{(i)};\theta_0, \theta_x, \theta_1, \theta_{\infty}, \hbox{tr}(M_0M_x),  \hbox{tr}(M_0M_1),
 \hbox{tr}(M_1M_x)\bigr), ~~~~i=0,1,\infty
$$
For the formulas of $a^{(i)}$ we refer to \cite{Jimbo}. The  monodromy data are 
restricted by the following condition, equivalent to (\ref{RESTOacasa}):
\be
|\hbox{tr}(M_iM_j)|\leq 2,~~~\Re\{\hbox{tr}(M_iM_j)\}\neq-2
\label{UFFAAAA}
\ee

 As explained above, we have shown that the transcendents (\ref{chineseramen}) have elliptic 
representation. Therefore,  Jimbo's transcendents are included in our class 
of transcendents obtained by the  elliptic representation. Observe that the behaviors 
(\ref{loc1introduzione})--(\ref{loc3introduzione}) are included in  the behaviors (\ref{LOc1introduzione})--(\ref{LOc3introduzione}) (and (\ref{LOc1introduzioNE})--(\ref{LOc3introduzioNE})), which  hold for any  $\nu_2^{(i)}\not \in (-\infty,0]\cup 
\{1\}\cup [2,+\infty)$.  Therefore the  condition (\ref{RESTOacasa}) is extended to any $
   \sigma^{(i)}\in {\bf C}$ such that $\sigma^{(i)}\not \in  (-\infty,0]\cup  [1,+\infty)$.  
 This corresponds to the fact that we have solved the connection problem for any complex 
value of 
$\hbox{tr}(M_iM_j)$ with the only constraint $\hbox{tr}(M_iM_j)\neq \pm 2$. This condition 
extends (\ref{UFFAAAA}).

 To be more precise, 
the condition  $\nu_2^{(i)}\neq 1$ is equivalent to $\hbox{tr}(M_0M_x)\neq 2$  
at $x=0$; to  
$\hbox{tr}(M_1M_x)\neq 2$ at $x=1$; to   
$\hbox{tr}(M_0M_1)\neq 2$ at $x=\infty$. Nevertheless, in the case $\hbox{tr}(M_iM_j)=2$ 
the critical behavior and the  solution of the connection problem were achieved  by Jimbo. 
 Unfortunately, the condition 
$\nu_2^{(i)}\neq 1$  which we had to impose to study the elliptic representation (except for 
 non-generic  
cases like $\beta=\gamma=1-2\delta=0$) does not allow us to know  the analytic properties and the critical behavior of the  elliptic representation in this case.  We expect that the properties of $u(x)$ are such to exactly produce the critical behavior found by Jimbo for $\hbox{tr}(M_iM_j)=2$, but we still have to cover this case.

The condition  $\nu_2^{(i)} \neq  0$ (and 2),    
implies that we can not give the critical 
behaviors (and the elliptic representation) of (\ref{chineseramen}) 
at $x=0$ for 
$\hbox{tr}(M_0M_x)=- 2$; at $x=1$ for 
$\hbox{tr}(M_1M_x)=- 2$; at $x=\infty$ for 
$\hbox{tr}(M_0M_1)=- 2$. 
 To our knowledge, these cases have not yet been studied in the literature.

To conclude, the results of \cite{Jimbo} together with our extension provide the 
critical behaviors and the solution of the connection  problem for the transcendents 
(\ref{chineseramen}) in the generic case for
$$
 \hbox{ any value of } \hbox{tr}(M_iM_j)\neq -2
$$
which corresponds to exponents 
$$ 
   \sigma^{(i)}\in {\bf C} \hbox{ such that } \sigma^{(i)}\not \in  (-\infty,0)\cup  [1,+\infty).
$$

\vskip 0.3 cm

 We turn now to the  special 
 case $\beta=\gamma=1-2\delta=0$, important  for its applications to topological filed theory, Frobenius manifolds \cite{Dub2} 
and  quantum cohomology \cite{KM}  \cite{guz2}.  
In this case is fully studied in \cite{guz1}. We can give a representation of   $u(x)$ in a 
domain which is wider than the generic case. Namely, at $x=0$, the domain is  
$$
{\cal D}(r;\nu_1,\nu_2):= \left\{ 
x\in \tilde{{\bf C}_0} ~|~ |x|<r,~\left| e^{-i\pi \nu_1} \left({x\over 16} 
\right)^{2-\nu_2}\right|<r,~\left| e^{i\pi \nu_1} \left({x\over 16} 
\right)^{\nu_2}\right|<r
\right\}
$$
In this domain $v(x)$ is holomorphic 
with convergent expansion
$$  
 v(x) = \sum_{n\geq 1} a_n x^n + \sum_{n\geq 0,m\geq 1} b_{nm}
x^n\left[e^{-i\pi \nu_1}\left({x\over 16} \right)^{2-\nu_2} \right]^m
+
 \sum_{n\geq 0,m\geq 1} c_{nm}
x^n\left[e^{i\pi \nu_1}\left({x\over 16} \right)^{\nu_2} \right]^m 
$$
If $\nu_2$ is real, the value $\nu_2=1$ is now allowed, namely, the 
constraint is $ 
   \nu_2 \not \in (-\infty,0]\cup [2,+\infty)
$. 
Therefore,  by periodicity of the $\wp$-function we can assume $
 0\leq \Re \nu_2 <2$, $\nu_2 \neq 0$.  
A similar result holds at $x=1$ and $x=\infty$ (see section \ref{so2qng}).

 According to \cite{DM}, we define  $
   2-x^2_0:=\hbox{tr} ~M_0M_x$, $2-x^2_1:=\hbox{tr}~ M_1 M_x$, 
$2-x^2_{\infty}:=\hbox{tr}~ M_0M_1$.  
 There is a one to one correspondence between triples 
 $(x_0,x_1,x_{\infty})$ (defined up to the change of two signs) 
and Painlev\'e transcendents, provided that at most one $x_i$ 
is zero and not all the $x_i$ are $\pm 2$ at the same time. Therefore we write 
$y(x)=y(x;x_0,x_1,x_{\infty})$. 
We show that one such transcendent has elliptic representations (half-periods are 
understood) 
$$
  y(x;x_0,x_1,x_{\infty}) =\wp \bigl(\nu_1^{(0)}\omega_1^{(0)}(x)+\nu_2^{(0)}\omega_2^{(0)}(x)+v^{(0)}(x;\nu_1^{(0)},\nu_2^{(0)})\bigr) +
{1+x\over 3}
$$
\be
 = \wp \bigl(\nu_1^{(1)}\omega_1^{(1)}(x)+\nu_2^{(1)}\omega_2^{(1)}(x)
+v^{(1)}(x;\nu_1^{(1)},\nu_2^{(1)})\bigr) +
{1+x\over 3}
\label{REF1}
\ee
\be
  = \wp \bigl(\nu_1^{(\infty)}\omega_1^{(\infty)}(x)+
\nu_2^{(\infty)}\omega_2^{(\infty)}(x)+v^{(\infty)}(x;\nu_1^{(\infty)},
\nu_2^{(\infty)})\bigr) +
{1+x\over 3}
\label{REF2}
\ee
 The parameters $\nu_2^{(i)}$ are obtained from 
 $$ 
\cos \pi \nu_2^{(i)}= {x_i^2\over 2} -1,~~~0\leq \Re \nu_2^{(i)}\leq 1,
~~~\nu_2^{(i)} \neq 0,  
~~~~i=0,1,\infty
$$ 
 Note  that the condition  $x_i\neq \pm 2$, 
$i=0,1,\infty$, corresponds to  $\nu_2^{(i)} \neq 0$. 
The parameter $\nu_1^{(0)}$ is obtained by the formula 
$$
  e^{i\pi \nu_1^{(0)}}= 
-{i \Gamma^4\left( 1-{\nu_2^{(0)}\over 2}\right)\over 2 \sin(\pi \nu_2^{(0)})\Gamma^2\left({3\over 2} 
-\mu - {\nu_2^{(0)}\over 2}\right) \Gamma^2\left({1\over 2} +\mu -{\nu_2^{(0)}\over 2}
\right)}
\left[2(1-e^{i\pi \nu_2^{(0)}}) -\right.
$$
$$
\left.  - f(x_0,x_1,x_{\infty})(x_{\infty}^2- e^{i\pi \nu_2^{(0)}}
x_1^2) \right]f(x_0,x_1,x_{\infty})
$$
where 
$$
f(x_0,x_1,x_{\infty}):= {4-x_0^2\over x_1^2+x_{\infty}^2-x_0x_1x_{\infty}}, ~~~
\alpha={(2\mu-1)^2\over 2}
$$  
Moreover, 
$\exp\{ -i \pi \nu_1^{(1)}\}$, $\exp\{i \pi \nu_1^{(\infty)}\}$ 
are given by an  analogous formula  with the substitutions 
$(x_0,x_1,x_{\infty}) 
\mapsto (x_1,x_0,x_0x_1-x_{\infty})$, $\nu_2^{(0)}\mapsto \nu_2^{(1)}$ and 
 $(x_0,x_1,x_{\infty}) 
\mapsto (x_{\infty}, -x_1,x_0-x_1x_{\infty})$,  $\nu_2^{(0)}\mapsto 
\nu_2^{(\infty)}$ respectively. 

 The most general choice 
 of $\nu_2$  is $0\leq \Re \nu_2 <2$. This corresponds to the fact that the 
  transcendent  $y(x;x_0,x_1,x_{\infty}) $ also  has three  representations 
$$
  y(x;x_0,x_1,x_{\infty}) =\wp \bigl(\tilde{\nu}_1^{(0)}\omega_1^{(0)}(x)+
\tilde{\nu}_2^{(0)}\omega_2^{(0)}(x)+v^{(0)}(x;\tilde{\nu}_1^{(0)},
\tilde{\nu}_2^{(0)})\bigr) +
{1+x\over 3}
$$
$$ = \wp \bigl(\tilde{\nu}_1^{(1)}\omega_1^{(1)}(x)+\tilde{\nu}_2^{(1)}\omega_2^{(1)}(x)+v^{(1)}(x;\tilde{\nu}_1^{(1)},\tilde{\nu}_2^{(1)})\bigr) +
{1+x\over 3}
$$
$$
  = \wp \bigl(\tilde{\nu}_1^{(\infty)}\omega_1^{(\infty)}(x)+
\tilde{\nu}_2^{(\infty)}\omega_2^{(\infty)}(x)+v^{(\infty)}(x;\tilde{\nu}_1^{(\infty)},\tilde{\nu}_2^{(\infty)})\bigr) +
{1+x\over 3}
$$
where
 $$ 
\cos \pi \tilde{\nu}_2^{(i)}= {x_i^2\over 2} -1,~~~1\leq \Re \nu_2^{(i)}<2,   
~~~~i=0,1,\infty
$$ 
The parameter $\tilde{\nu}_1^{(0)}$ is obtained by the formula 
$$
  e^{-i\pi \tilde{\nu}^{(0)}_1}= 
{i \Gamma^4\left( {\tilde{\nu}^{(0)}_2\over 2}\right)\over 2 \sin(\pi \tilde{\nu}_2^{(0)})
\Gamma^2\left({1\over 2} 
-\mu + {\tilde{\nu}_2^{(0)}\over 2}\right) \Gamma^2\left(-{1\over 2} +\mu +{\nu_2\over 2}
\right)}
\left[2(1-e^{-i\pi \tilde{\nu}_2^{(0)}})-\right.
$$
$$
\left. 
 - f(x_0,x_1,x_{\infty})(x_{\infty}^2- e^{-i\pi \tilde{\nu}^{(0)}_2}
x_1^2) \right]f(x_0,x_1,x_{\infty}).
$$
$\exp\{ i \pi \tilde{\nu}_1^{(1)}\}$, 
$\exp\{- i \pi \tilde{\nu}_1^{(\infty)}\}$ 
are given by an  analogous formula  with the substitutions 
$(x_0,x_1,x_{\infty}) 
\mapsto (x_1,x_0,x_0x_1-x_{\infty})$, $\tilde{\nu}_2^{(0)}\mapsto \tilde{\nu}_2^{(1)}$ and 
 $(x_0,x_1,x_{\infty}) 
\mapsto (x_{\infty}, -x_1,x_0-x_1x_{\infty})$,  $\tilde{\nu}_2^{(0)}\mapsto 
\tilde{\nu}_2^{(\infty)}$ respectively. 

\vskip 0.2 cm 
\noindent 
The formulae above have limits for $\nu_2=1,1\pm 2\mu+2m$, 
$m$ integer. They are listed in subsection \ref{so2qng}.

\vskip 0.3 cm 
 Conversely,  a transcendent  
\be
y(x)= \wp \bigl({\nu}_1\omega_1^{(0)}(x)+
{\nu}_2\omega_2^{(0)}(x)+v^{(0)}(x;{\nu}_1,
{\nu}_2)\bigr) +
{1+x\over 3},~~~ \hbox{ at } x=0
\label{startingpoint}
\ee
coincides with $y(x;~x_0,x_1,x_{\infty})$, with the following monodromy data.

\vskip 0.2 cm 
\noindent 
If $0\leq \Re \nu_2 \leq 1$:  
$$
  x_0=2\cos\left({\pi \over 2}\nu_2\right)
$$
$$
x_1=  \left[{ 4^{-\nu_2}~2~e^{i {\pi\over 2} \nu_1 } \over f(\nu_2,\mu) G(\nu_2,\mu)}+{G(\nu_2,\mu) \over  4^{-\nu_2}~2~e^{i{ \pi\over 2} \nu_1 }} \right]
$$
$$
x_{\infty} = \left[ { 4^{-\nu_2}~2~e^{i{\pi \over 2} (\nu_1-\nu_2)}\over 
  f(\nu_2,\mu) G(\nu_2,\mu)}+ {G(\nu_2,\mu)\over  4^{-\nu_2}~2~e^{i{\pi \over 2} (\nu_1-\nu_2)}} \right]
$$
where 
$$ 
  f(\nu_2,\mu)= - {2\sin^2\left({\pi \over 2} \nu_2\right) \over 
\cos(\pi \nu_2) +\cos(2\pi \mu)}
,~~~G(\nu_2,\mu)=  4^{-\nu_2}~2~ { \Gamma\left(1-{\nu_2\over 2}\right)^2 \over
\Gamma\left({3\over 2} - \mu - {\nu_2 \over 2} \right) 
\Gamma\left({1\over 2} +\mu - {\nu_2\over 2} \right) 
}
$$

\vskip 0.2 cm 
\noindent 
If $1\leq \Re \nu_2 <2$:  
$$
x_0=  2 \cos \left({\pi\over 2} \nu_2 \right)
$$
$$
x_1= \left[ {e^{-i{\pi \over 2} \nu_1} \over  4^{1-\nu_2}~ 2~f(\nu_2,\mu) G_1(\nu_2,\mu)} + {  4^{1-\nu_2}~  2~G_1(\nu_2,\mu) \over e^{-i{\pi \over 2} \nu_1}}
 \right] 
$$
$$ 
x_{\infty}
= \left[{ e^{i {\pi \over 2} (\nu_2-\nu_1)} \over  4^{1-\nu_2}~2~ f(\nu_2,\mu) 
G_1(\nu_2,\mu) } +{   4^{1-\nu_2}~2~ G_1(\nu_2,\mu)\over e^{i {\pi \over 2} (\nu_2-\nu_1)}} \right]
$$
where 
$$ 
  G_1(\nu_2,\mu)= {1\over 4^{1-\nu_2}~2} {\Gamma\left({\nu_2\over 2} \right)^2
\over 
\Gamma\left( {1\over 2}- \mu +{\nu_2\over 2}\right)\Gamma\left(
-{1\over 2} +\mu+{\nu_2\over 2}\right)}
$$
 After computing the monodromy data, we can write the  elliptic 
representations of  $y(x;x_0,x_1,x_{\infty})$ at $x=1$ and $x=\infty$, namely 
(\ref{REF1}), (\ref{REF2}).    Since they 
 are the elliptic representations at $x=1$, $x=\infty$  of 
(\ref{startingpoint}),  we 
 have solved the connection problem for (\ref{startingpoint}).

\vskip 0.2 cm
 
 We observed that there  is a one to one correspondence between   
 Painlev\'e 
transcendents and   triples
 of monodromy data $(x_0,x_1,x_{\infty})$, 
defined up to the change of two signs,  satisfying $
x_i \neq \pm 2$, $i=0,1,\infty$, i.e. $\nu_2^{(i)}\neq 0$ (and 2),   
and at most one $x_i=0$. The cases when  these conditions are not satisfied are studied in \cite{M}. 
However, if  $x_i=\pm2$ (namely the  trace is $ -2$) the problem of 
finding the critical behavior at the corresponding critical point $x=i$ is still 
open (except when  {\it all the 
three} $x_i$  are $\pm 2$: in this case  there is a one-parameter class of 
 solutions called {\it Chazy solutions } in \cite{M}). We conclude that the results of our paper 
(and \cite{guz1}), 
plus the results of \cite{M} cover all the possible transcendents, except the special case when 
 one or two 
$x_i$ are $\pm 2$. We plan to cover this last case soon.

\vskip 0.3 cm 
 Finally, we expect that  in  all  non-generic cases we can solve the connection problem and express the parameters $\nu_1$, $\nu_2$ in terms of monodromy data.  
 From the conceptual point of 
 view nothing should change with respect to  \cite{Jimbo} \cite{DM} \cite{guz1} and the 
present paper; but the technical details may require a long time for computations.


\section{The Elliptic Representation}\label{The Elliptic Representation}

 We
 derive the elliptic form for the general Painlev\'e 6 equation. We
follow \cite{fuchs}. Let
\be
u=\int_{\infty}^y ~{d\lambda \over \sqrt{\lambda(\lambda-1)(\lambda-x)}}
\label{elliptint1}
\ee
We observe that
$$
  {du \over dx}= {\partial u \over \partial y }~{dy \over dx} + {\partial u
  \over \partial x } ={1\over \sqrt{y(y-1)(y-x)} } ~
{dy \over dx} + {\partial u
  \over \partial x } 
$$
from which we compute
$$
{d^2 u\over dx^2}+{2x-1\over x(x-1)} {du\over dx} +{u\over 4 x(x-1)}=$$
$$
  = {1\over  \sqrt{y(y-1)(y-x)}}\left[{d^2 y\over dx^2}+\left({1\over x}
  +{1\over x-1}+{1\over y-x}\right) ~{dy \over dx}-{1\over 2}\left({1\over y}
  +{1\over y-1}+{1\over y-x}\right)~\left({dy \over dx}\right)^2 \right]
$$
$$
+{\partial^2 u \over \partial x^2}+{2x-1\over x(x-1)}~{\partial u \over
  \partial x}+ {u\over 4x(x-1)} 
$$
By direct calculation we have:
$$
{\partial^2 u \over \partial x^2}+{2x-1\over x(x-1)}~{\partial u \over
  \partial x}+ {u\over 4x(x-1)} 
= -{1\over 2} {\sqrt{y(y-1)(y-x)} \over x(x-1)}
  ~{1\over (y-x)^2}
$$
 Therefore, $y(x)$ satisfies PVI if and only if 
\be
 {d^2 u\over dx^2}+{2x-1\over x(x-1)}~ {du\over dx} +{u\over 4
x(x-1)}={\sqrt{y(y-1)(y-x)} \over 2 x^2 (1-x)^2}~\left[2 \alpha+2\beta{x\over
y^2} +\gamma {x-1\over (y-1)^2} +\left(\delta-{1\over2}\right){x(x-1)\over
(y-x)^2} \right]
\label{ellipainleve1}
\ee

\vskip 0.2 cm
We invert the function $u=u(y)$ by observing that we are dealing with an
elliptic integral. Therefore, we write
$$ 
  y=f(u,x)$$
where $f(u,x)$ is an elliptic function of $u$. This implies that 
$$
  {\partial y\over \partial u}=\sqrt{y(y-1)(y-x)}
$$
 The above equality allows us to  rewrite (\ref{ellipainleve1}) in the
 following way: 
\be
x(1-x)~{d^2 u\over dx^2}+(1-2x)~{du\over dx}-{1\over 4}~u = {1\over 2x(1-x)}
 ~{\partial \over \partial u}\psi(u,x),
\label{ellipainleve2}
\ee
where 
$$
\psi(u,x):= 2 \alpha f(u,x)-2\beta{x\over f(u,x)} +2\gamma{1-x\over f(u,x)-1}
+(1-2\delta){x(x-1)\over f(u,x)-x}
$$

\vskip 0.2 cm
The last step concerns the form of $f(u,x)$. We observe that 
$4\lambda(\lambda-1)(\lambda-x)$ is not in Weierstrass canonical form. We change variable:
$$\lambda= t +{1+x\over 3},$$
and we get the Weierstrass form:
$$
    4\lambda(\lambda-1)(\lambda-x)= 4 t^3-g_2 t-g_3,~~~~ 
 g_2={4\over 3}(1-x+x^2),~~~~g_3:={4\over
    27}(x-2)(2x-1)(1+x)
$$
Thus
$$
  {u\over 2}=\int_{\infty}^{y-{1+x\over 3}}~{dt\over  \sqrt{4 t^3-g_2 t-g_3}} 
$$
which implies 
$$
   f(u,x)= {\cal P}\left({u\over 2}; \omega_1,\omega_2\right)+{1+x\over 3}
$$
 We still need to explain what are the {\it half periods} $\omega_1$,
 $\omega_2$. In order to do that, we first observe that the Weierstrass form is 
$$
   4 t^3-g_2 t-g_3= 4(t-e_1)(t-e_2)(t-e_3)
$$
where
$$ 
    e_1= {2-x\over 3},~~~e_2={2x-1\over 3},~~~e_3=-{1+x\over 3}.
$$
Therefore
$$
   g:=\sqrt{e_1-e_2}=1,~~~{\kappa}^2:={e_2-e_3\over e_1-e_3} =
   x,~~~{{\kappa}^{\prime}}^2 := 1- {\kappa}^2=1-x
$$
 We identify $e_1=\wp(\omega_1)$, $e_2=\wp(\omega_1+\omega_2)$, 
$e_3=\wp(\omega_2)$. Therefore,   the half-periods are 
$$
\omega_1= {1\over g} \int_0^1 ~{d\xi\over \sqrt{(1-\xi^2)(1-\kappa^2 \xi^2)}} 
= \int_0^1 ~{d\xi\over \sqrt{(1-\xi^2)(1-x \xi^2)}} ={\bf K}(x)
$$
$$
 \omega_2= {i\over g} \int_0^1 ~{d\xi\over
\sqrt{(1-\xi^2)(1-{\kappa^{\prime}}^2  
\xi^2)}} 
= i~\int_0^1 ~{d\xi\over \sqrt{(1-\xi^2)(1-(1-x) \xi^2)}} =i{\bf
K}(1-x) 
$$
The elliptic integral ${\bf K}(x)$ coincides, up to a factor, 
 with the hyper-geometric function 
$
F\left({1\over 2},{1\over 2},1;x\right)$, whose convergent series for $|x|<1$ 
is: 
$$
F\left({1\over 2},{1\over 2},1;x\right)=  
\sum_{n=0}^{\infty}{ \left[\left({1 \over 2}\right)_n\right]^2
  \over (n!)^2 } x^n.
$$
Namely: 
$$
{\bf K}(x)= {\pi \over 2} ~F\left({1\over 2},{1\over 2}, 1; x\right)
$$
Moreover, 
${\bf K}(x)$ and ${\bf K}(1-x)$ 
 are two linearly independent solutions of the hyper-geometric equation 
$$
 x(1-x) \omega^{\prime\prime}+(1-2x) \omega^{\prime} -{1\over 4} \omega=0.
$$
Observe that for $|\hbox{arg}(x)|<\pi$: 
$$-\pi F\left({1\over 2},{1\over 2}, 1;1- x\right) =  F\left({1\over
2},{1\over 2},1;x\right)  \ln(x) +
F_1(x)  
$$ 
where 
$$
F_1(x):=  
\sum_{n=0}^{\infty}{ \left[\left({1 \over 2}\right)_n\right]^2
  \over (n!)^2 } 2\left[ \psi(n+{1\over 2}) - \psi(n+1)\right]
x^n,
$$
$$\psi(z) = 
{d \over dz}\ln \Gamma(z),~~~ \psi\left({1\over 2}\right) = -\gamma -2 \ln
2,~~~ \psi(1)=-\gamma,~~~\psi(a+n)=\psi(a)+\sum_{l=0}^{n-1} {1\over a+l}.
$$
 Therefore 
$$
\omega_2(x)= -{i\over 2}[F\left({1\over 2},{1\over 2},
1;x\right)\ln(x)+F_1(x)]
$$ 
In the following we use sometimes the abbreviation  $F(x)$ 
 for $ F\left({1\over 2},{1\over 2},
1;x\right)$. The series of $F(x)$ and $F_1(x)$ 
 converge for $|x|<1$.


\vskip 0.3 cm 

Let 
$${\cal L}(u)
  :=
 x(1-x) {d^2u\over dx^2}+(1-2x) {du\over dx} -{1\over 4} u
$$ 
The PVI equation becomes
$$
{\cal L}(u)= {1\over 2x(1-x)} {\partial \over \partial u} \left\{ 2\alpha\left[
\wp\left({u\over 2};\omega_1,\omega_2\right)-e_3\right]-2\beta {x\over 
\wp\left({u\over 2};\omega_1,\omega_2\right)-e_3}+\right.
$$
\be
 \left. +2\gamma{1-x \over \wp\left({u\over 2};\omega_1,\omega_2\right)-e_1}+(1-2\delta){x(1-x)\over \wp\left({u\over 2};\omega_1,\omega_2\right)-e_2}
\right\}
\label{ellittic1}
\ee
We recall that 
$$
{1\over \wp\left({u\over 2}\right)-e_3}= {1\over (e_1-e_3)(e_2-e_3)}
\left[ \wp\left({u\over 2}+\omega_2\right) -e_3  \right]
$$
$$
{1\over \wp\left({u\over 2}\right)-e_1}= {1\over (e_1-e_2)(e_1-e_3)}
\left[ \wp\left({u\over 2}+\omega_1\right) -e_1  \right]
$$
$$
{1\over \wp\left({u\over 2}\right)-e_2}= {1\over (e_1-e_2)(e_3-e_2)}
\left[ \wp\left({u\over 2}+\omega_1+\omega_2\right) -e_2  \right].
$$
We observe also that $e_1-e_3=1$, $e_2-e_3=x$, $e_1-e_2=1-x$. Therefore, 
(\ref{ellittic1}) becomes: 
$$
{\cal L}(u)= {1\over 2x (1-x)} \left[ 
                      2\alpha {\partial \over \partial u} 
\wp\left({u\over 2} \right) -2\beta  {\partial \over \partial u} 
\wp\left({u\over 2}+\omega_2 \right)+\right.
$$
\be
\left.+2\gamma   {\partial \over \partial u} 
\wp\left({u\over 2}+\omega_1 \right) + (1-2\delta) {\partial \over \partial u} 
\wp\left({u\over 2}+\omega_1 +\omega_2\right)
\right]
\label{fuchseq}
\ee


\section{Proof of  Theorem 1}
 For technical reason it is convenient to introduce the following definitions: 
 \be
{\cal D}_1(r;\nu_1,\nu_2):= 
\left\{ x\in \widetilde{{\bf C}_0}~ \hbox{ such that }
  |x|<r, \left|{e^{-i\pi \nu_1}\over 16^{1-\nu_2}} x^{1-\nu_2} 
\right|<r,
\left| 
{e^{i\pi \nu_1} \over 16^{\nu_2}} x^{\nu_2}\right|<r \right\}
\label{DOMINOO1}
\ee
\be
{\cal D}_2(r;\nu_1,\nu_2):= 
\left\{ x\in \widetilde{{\bf C}_0}~ \hbox{ such that }
  |x|<r, \left|{e^{-i\pi \nu_1}\over 16^{2-\nu_2}} x^{2-\nu_2} 
\right|<r,
\left| 
{e^{i\pi \nu_1} \over 16^{\nu_2-1}} x^{\nu_2-1}\right|<r \right\}
\label{DOMINOO2}
\ee
Note that ${\cal D}_1(r;\nu_1,\nu_2)\equiv {\cal D}(r;\nu_1,\nu_2)$ in 
 (\ref{DOMINOO1talk})
 and 
${\cal D}_2(r;\nu_1,\nu_2)\equiv {\cal D}(r;-\nu_1,2-\nu_2)$. 
Let 
$$
  \tau(x):= {\omega_2(x) \over \omega_1(x)}= {1\over \pi} (\arg x - i \ln{|x|
\over 16}) + O(x)
 $$
be the modular parameter of the elliptic function, obtained expanding 
$\omega_1(x)$, $\omega_2(x)$ at $x=0$.   
We expand the r.h.s of (\ref{fuchseq}) 
in Fourier series. We are dealing with 
Weierstrass $\wp$-functions of the form $
   \wp\left({u\over 2}+\varepsilon_1 \omega_1 + \varepsilon_2 \omega_2\right)
$, where 
$\varepsilon_i\in \{0,1\}$, $i=1,2$. Due to periodicity of $\wp$ we are going to expand
$$ 
 \wp\left({u\over 2}+[\varepsilon_1+2N_1]~ \omega_1 + 
[\varepsilon_2 +2N_2]~\omega_2\right), ~~~~\varepsilon_i\in \{0,1\},~~~
N_i\in{\bf Z},~~~i=1,2 
$$
As it is well known \cite{SANJ}, the Fourier expansion of the $\wp$-function is $$
\wp(z;\omega_1,\omega_2)= \left({\pi \over 2 \omega_1}\right)^2~ 
\left[
-{1\over 3} +{1\over \sin^2\left({\pi z \over 2 \omega_1}\right)} +
8 \sum_{n\geq 1} {n e^{2 i \pi  n \tau} \over 1- e^{2 i \pi n \tau}} ~\left[ 
1-\cos\left({n\pi z\over \omega_1} \right)
\right] 
\right],~~~\hbox{ for } \Im \tau > \left|\Im \left({z\over 2 \omega_1}\right)
 \right|
$$
Therefore, the expansion  is possible provided that 
\be
   \left| \Im \left({u(x)\over 4\omega_1(x)}+{\varepsilon_1+2N_1\over 2} 
+ {\varepsilon_2 +2N_2\over 2}~ \tau(x)\right) \right|<\Im \tau(x)
\label{condizione}
\ee
We observe that for $x\to 0$ 
$$
   \Im \tau(x) = -{1\over \pi} \ln{|x|\over 16}+ O(x) \to +\infty 
$$
The condition (\ref{condizione}) becomes 
\be 
    -\Im \tau < {1\over 2} \Im \nu_1 + {1\over 2} [\varepsilon_2 +2 N_2 +
\Re \nu_2] ~\Im \tau +{1\over 2} \Im \nu_2~ \Re \tau  + \Im \left({v\over 2
\omega_1} \right) < \Im \tau
\label{SAKADOM}
\ee
If $\nu_2$ is real, we divide by $\Im \tau$, we let $x \to 0$ and we obtain 
$$ 
   -2-2N_2-\varepsilon_2<\nu_2<2- 2 N_2 - \epsilon_2
$$
provided that $v(x)$ is bounded.  
In general,   (\ref{SAKADOM}) becomes
$$
   (\Re \nu_2+2 +\varepsilon_2 +2N_2) \ln{|x|\over 16} -\left[
\pi \Im \nu_1 + \Im \left(
{\pi v \over \omega_1}
 \right)
\right] + O(x) 
$$
$$
\leq \Im \nu_2 \arg x \leq 
$$
\be
 (\Re \nu_2 -2+ \varepsilon_2 +2N_2)\ln{|x|\over 16} 
- \left[ \pi \Im \nu_1 + \Im \left(
{\pi v \over \omega_1}
 \right)
\right] + O(x) 
\label{PALLASAN}
\ee
In the domain defined by the above condition, and $|x|<1$, we do the Fourier 
expansion as follows 
$$
{\partial \over \partial u} \wp \left({u\over 2} +[\varepsilon_1+2N_1] \omega_1 
+[\varepsilon_2 +2N_2]\omega_2 \right)=
$$
$$
  = {\pi^3\over 8 \omega_1^3} 
\left\{
8 \sum_{n\geq 1} {n^2 e^{ 2 i \pi n \tau} \over 1 - e^{ 2 i \pi n \tau}} 
\sin
\left[\pi n 
\left(
{u \over 2 \omega_1} + \epsilon_1 + (\epsilon_2 +2N_2)\tau 
 \right) 
\right] - 
{\cos
\left[ {\pi \over 2} 
\left(
{u \over 2 \omega_1} + \epsilon_1 + (\epsilon_2 +2N_2)\tau 
\right)
\right] \over 
\sin^3 
\left[ 
{\pi \over 2} 
\left(
{u \over 2 \omega_1} + \epsilon_1 + (\epsilon_2 +2N_2)\tau 
\right)
\right]}
 \right\}
$$
\be
=
{i\pi^3\over 2 \omega_1^3} \left\{ 
{e^{ i{f(x) \over 2}} +e^{- i{f(x) \over 2}}
\over \left[ e^{ i{f(x) \over 2} } -e^{- i{f(x) \over 2}  }
\right]^3
}
- \sum_{n\geq 1} {n^2 e^{2i \pi n \tau} \over 1 -  e^{2i \pi n \tau}} 
\left[ 
e^{ i n f(x)} -e^{- i nf(x)}
\right]
\right\}
\label{espan}
\ee
where
$$
e^{i f(x)}= \exp\left\{i \pi \left[(\nu_1+\varepsilon_1)+(\nu_2+\varepsilon_2+2N_2) \tau(x) 
+{v(x)\over \omega_1(x)}\right] \right\}
$$
Note that $N_1$ does not appear in the expansion. So, we can take  $N_1=0$
 in the following. 

\vskip 0.2 cm 
Now we make  the {\it further assumption} that 
\be
\left|e^{i f(x)} \right|<1
\label{assume}
\ee
and we  rewrite ({\ref{espan}) as follows:
$$
    {i\pi^3\over 2 \omega_1^3} \left\{
 {e^{2if(x)} + e^{i f(x)} \over 
\left[
e^{i f(x)} -1\right
]^3}  - \sum_{n\geq 1}{ n^2 e^{i\pi n 
\left[ 
-(\nu_1+\varepsilon_1)+(2-\nu_2-\varepsilon_2-2N_2)\tau - {v\over \omega_1}
\right]}\over 1-e^{2\pi n \tau}}
~\left[ 
e^{2inf(x)}-1
\right]
\right\}
$$
Due to (\ref{assume}), the denominator in the first term does not vanish, so that the expansion has no poles.  
The condition  (\ref{assume}) is more restrictive than (\ref{PALLASAN}) and 
yields:  
$$
   (\Re \nu_2 +\varepsilon_2 +2N_2) \ln{|x|\over 16} -\left[
\pi \Im \nu_1 + \Im \left(
{\pi v \over \omega_1}
 \right)
\right] + O(x) 
$$
$$
\leq \Im \nu_2 \arg x \leq 
$$
\be
 (\Re \nu_2 -2+ \varepsilon_2 +2N_2)\ln{|x|\over 16} 
- \left[ \pi \Im \nu_1 + \Im \left(
{\pi v \over \omega_1}
 \right)
\right] + O(x) 
\label{newdomain}
\ee
We observe that, for any complex number $C$,    
$$
e^{i \pi C \tau(x)}=h(x)^C~
 \left[ {x \over 16}\right]^C,
$$
where
$$
h(x):= 
 e^{C\left({F_1(x) \over F(x)}+4\ln2 \right)}= 1 +O(x),~~~~x\to 0
$$
Therefore, 
$$ 
    {\partial\over \partial u} \wp \left({u\over 2}+
 (\varepsilon_1+2N_1)\omega_1 
+(\varepsilon_2+2N_2)\omega_2\right)
={i\pi^3\over 2\omega_1^3}\left\{   {   A^2+  A
\over 
\left(A-1 \right)^3
}
- 
\sum_{n\geq 1} {n^2 B^n
\over 
1-h(x)^{2n} \left({x\over 16}\right)^{2n}}
\left( A^{2n}-1
\right)
\right\}
$$
where
\be
A=A\left(x,e^{i\pi\nu_1}\left({x\over 16} \right)^{
\nu_2+\epsilon_2+2N_2}e^{i\pi{v\over \omega_1}}\right):=
h(x)^{\nu_2+\varepsilon_2+2N_2}  e^{i\pi (\varepsilon_1+\nu_1)}
\left({x\over 16}\right)^{\nu_2+\varepsilon_2+2N_2}e^{i{\pi v \over \omega_1}}
\label{dipendenzaa}
\ee
\be
B=B\left(x,e^{-i\pi\nu_1} \left({x\over 16} \right)^{2-\nu_2-\varepsilon_2-2N_2}
e^{-i\pi{v\over \omega_1}}\right):=h(x)^{2-\nu_2-\varepsilon_2-2N_2} e^{-i\pi(\nu_1+\varepsilon_1)}\left(
{x\over 16}
\right)^{2-\nu_2-\varepsilon_2-2N_2}e^{-i \pi {v\over \omega_1}}
\label{dipendenzab}
\ee

\vskip 0.2 cm 
 In the r.h.s. of 
equation (\ref{fuchseq}), we have four $\wp$-functions with different 
arguments. We recall that we can take $N_1=0$. As for $N_2$, it can be different in each $\wp$-function. There are only two possibilities in order 
 for the four resulting domains 
(\ref{newdomain})  to intersect: the $\wp$-functions must be 
$$
\wp\left({u\over 2} +2N_2\omega_2 \right),~~~ 
~~~\wp\left({u\over 2}+\omega_1 +2N_2\omega_2 \right),~~~
\wp\left({u\over 2}+\omega_2 +2N^{\prime}_2 \omega_2\right),~~~ 
~~~\wp\left({u\over 2}+\omega_1 +\omega_2+2N^{\prime}_2\omega_2 \right),
$$
and 
$$ 
 (a):~~~~~ N_2^{\prime}=N_2
$$
or
$$
(b):~~~~~ N_2^{\prime}=N_2-1
$$
 Thus, the r.h.s. of (\ref{fuchseq}) is expanded in Fourier  
series if $|x|<1$, in the following domains.   

\vskip 0.2 cm
Case (a):
$$
(\Re \nu_2+2N)\ln{|x|\over 16} -\left[
\pi \Im \nu_1 + \Im \left(
{\pi v \over \omega_1}
 \right)
\right] + O(x)<\Im \nu_2 \arg x <
$$
\be < 
(\Re \nu_2-1+2N)\ln{|x| \over 16} -\left[
\pi \Im \nu_1 + \Im \left(
{\pi v \over \omega_1}
 \right)
\right] + O(x)
,~~~~\hbox{ for } N_2^{\prime}=N_2=:N
\label{dominiodim1}
\ee

\vskip 0.2 cm
Case (b):
$$
(\Re \nu_2-1+2N)\ln{|x|\over 16} -\left[
\pi \Im \nu_1 + \Im \left(
{\pi v \over \omega_1}
 \right)
\right] + O(x)<\Im \nu_2 \arg x < 
$$
\be <
(\Re \nu_2-2+2N)\ln{|x|\over 16} -\left[
\pi \Im \nu_1 + \Im \left(
{\pi v \over \omega_1}
 \right)
\right] + O(x)
,~~~~\hbox{ for } N_2^{\prime}+1=N_2=:N
\label{dominiodim2}
\ee
 Note that if $\Im \nu_2=0$ any value of $\arg(x)$ is allowed, but dividing 
the previous expressions  
by $\ln|x|$ and letting $x\to 0$ we have 
$$
   -2N <\nu_2<1-2N \hbox{ for (a)},~~~~~1-2N<\nu_2<2-2N \hbox{ for (b)}
$$
provided that $v(x)$ is bounded fort $x\to 0$. 
Taking (\ref{dipendenzaa})  (\ref{dipendenzab}) into account in the two cases, we conclude  that: 

\vskip 0.2 cm 

In case (a), for $ N_2^{\prime}=N_2=N$, the following powers appear in the 
r.h.s. of (\ref{fuchseq}):
$$
 x^{2-\nu_2-2N},~~~x^{\nu_2+2N},~~~x^{1-\nu_2-2N},~~~x^{1+\nu_2+2N}
$$
The dominant powers are then
$$ 
x^{\nu_2+2N},~~~x^{1-\nu_2-2N}
$$
and we can write  (\ref{fuchseq}) as;
$$ 
   {\cal L}(u)= {1\over 2x(1-x)} {\cal F}\left(x,e^{-i\pi \nu_1} \left(
{x\over 16}
\right)^{1-\nu_2-2N} e^{-i\pi{v\over \omega_1}},e^{i\pi \nu_1} \left(
{x\over 16}\right)^{\nu_2+2N} e^{i\pi{v\over \omega_1}}
\right)
$$
on the domain defined by  $|x|<1$ and (\ref{dominiodim1}), 
assuming that $v(x)$ is 
bounded. 

\vskip 0.2 cm 

In case (b), for  $N_2^{\prime}+1=N_2=N$,  the following powers appear 
in the 
r.h.s. of (\ref{fuchseq}):
$$
 x^{2-\nu_2-2N},~~~x^{\nu_2+2N},~~~x^{3-\nu_2-2N},~~~x^{\nu_2-1+2N}
$$
The dominant powers are then
$$ 
x^{2-\nu_2-2N},~~~x^{\nu_2-1+2N}
$$
and we can write  (\ref{fuchseq}) as;
$$ 
   {\cal L}(u)= {1\over 2x(1-x)} {\cal F}\left(x,e^{-i\pi \nu_1} \left(
{x\over 16}
\right)^{2-\nu_2-2N} e^{-i\pi{v\over \omega_1}},e^{i\pi \nu_1} \left(
{x\over 16}\right)^{\nu_2-1+2N} e^{i\pi{v\over \omega_1}}
\right)
$$
on the domain defined by  $|x|<1$ and (\ref{dominiodim2}), 
assuming that $v(x)$ is 
bounded. 

\vskip 0.2 cm 
Let $\epsilon<1$ be sufficiently small: in both cases, 
${\cal F}(x,y,z)$ is holomorphic for $|x|,|y|,|z|<\epsilon$ 
and ${\cal F}(0,0,0)=0$.  
In (a) we decompose 
$$ 
{\cal F}\left(x,e^{-i\pi \nu_1} \left(
{x\over 16}
\right)^{1-\nu_2-2N} e^{-i\pi{v\over \omega_1}},e^{i\pi \nu_1} \left(
{x\over 16}\right)^{\nu_2+2N} e^{i\pi{v\over \omega_1}}
\right)=
$$
$$
= 
{\cal F}\left(x,e^{-i\pi \nu_1} \left(
{x\over 16}
\right)^{1-\nu_2-2N},e^{i\pi \nu_1} \left(
{x\over 16}\right)^{\nu_2+2N} 
\right)+
{\cal G}\left(x,e^{-i\pi \nu_1} \left(
{x\over 16}
\right)^{1-\nu_2-2N},e^{i\pi \nu_1} \left(
{x\over 16}\right)^{\nu_2+2N},v 
\right)
$$
where 
$$
{\cal G}\left(x,e^{-i\pi \nu_1} \left(
{x\over 16}
\right)^{1-\nu_2-2N},e^{i\pi \nu_1} \left(
{x\over 16}\right)^{\nu_2+2N},v 
\right):=
$$
$${\cal F}\left(x,e^{-i\pi \nu_1} \left(
{x\over 16}
\right)^{1-\nu_2-2N} e^{-i\pi{v\over \omega_1}},e^{i\pi \nu_1} \left(
{x\over 16}\right)^{\nu_2+2N} e^{i\pi{v\over \omega_1}}
\right)  -$$
$$
- {\cal F}\left(x,e^{-i\pi \nu_1} \left(
{x\over 16}
\right)^{1-\nu_2-2N},e^{i\pi \nu_1} \left(
{x\over 16}\right)^{\nu_2+2N} 
\right)
$$ 

We  note that ${\cal L}(u)={\cal L}(u_0+2v)$=${\cal L}(u_0)+
2{\cal L}(v)\equiv 2{\cal L}(v)$. Let us put
$$
   w(x):= x {d\over dx} v(x)
$$
 Equation (\ref{fuchseq}) becomes the system 
$$
\left\{\matrix{x{dv\over dx} = w 
\cr
\cr
x{dw\over dx} = \Phi+\Psi}
\right. 
$$
where 
$$\Phi= \Phi\left(x,e^{-i\pi \nu_1} \left(
{x\over 16}
\right)^{1-\nu_2-2N},e^{i\pi \nu_1} \left(
{x\over 16}\right)^{\nu_2+2N} 
\right):=$$
$$ =
{1 \over 4(1-x)^2}{\cal F}\left(x,e^{-i\pi \nu_1} \left(
{x\over 16}
\right)^{1-\nu_2-2N},e^{i\pi \nu_1} \left(
{x\over 16}\right)^{\nu_2+2N} 
\right)
$$
and
$$\Psi=\Psi\left(x,e^{-i\pi \nu_1} \left(
{x\over 16}
\right)^{1-\nu_2-2N},e^{i\pi \nu_1} \left(
{x\over 16}\right)^{\nu_2+2N},v,w 
\right):=$$
$$=
 {x(w+v/4)\over 1-x} +{1 \over 4(1-x)^2}
{\cal G}\left(x,e^{-i\pi \nu_1} \left(
{x\over 16}
\right)^{1-\nu_2-2N},e^{i\pi \nu_1} \left(
{x\over 16}\right)^{\nu_2+2N},v 
\right)
$$
\vskip 0.2 cm 
 In the same way, in case (b) we obtain a similar system, where 
$$\Phi= \Phi\left(x,e^{-i\pi \nu_1} \left(
{x\over 16}
\right)^{2-\nu_2-2N},e^{i\pi \nu_1} \left(
{x\over 16}\right)^{\nu_2-1+2N} 
                                      \right) 
$$
$$\Psi=\Psi\left(        x,
e^{-i\pi \nu_1}               \left(
{x\over 16}
                              \right)^{2-\nu_2-2N},
e^{i\pi \nu_1}     \left(
{x\over 16}             \right)^{\nu_2-1+2N},v,w 
\right)
$$

\vskip 0.2 cm
By construction, in both cases $\Phi(x,y,z)$ and $\Psi(x,y,z,v,w)$ are holomorphic of their arguments for $|x|,|y|,|z|,|v|,|w|<\epsilon$ for sufficiently 
small $\epsilon<1$ (assuming $v(x)$ bounded in a neighborhood of $x=0$). 
Moreover 
$$
 \Phi(0,0,0)=\Psi(0,0,0,v,w)=\Psi(x,y,z,0,0)=0
$$

\vskip 0.3 cm
We complete the proof of the theorem in case (a). Case (b) is analogous. 
 We reduce the system of differential equations to
 a system of integral equations 
$$
w(x)=\int_{L(x)} {1\over s} \left\{
 \Phi (s, e^{-i\pi \nu_1} \left({s\over 16}\right)^{1-\nu_2-2N}, 
e^{i\pi \nu_1}  \left({s\over 16}\right)^{\nu_2+2N} ) + \right. 
$$
$$ \left. 
 + \Psi(s, e^{-i\pi \nu_1} \left({s\over 16}\right)^{1-\nu_2-2N}, 
e^{i\pi \nu_1}  \left({s\over 16}\right)^{\nu_2+2N},v(s),w(s))\right\}~ds
$$
$$
v(x)=\int_{L(x)}{1\over s}w(s)~ds
$$
The point $x$ and 
the path of integration  are chosen to belong to the domain where $|x|$, 
$|{e^{-i\pi \nu_1}\over 16^{1-\nu_2-2N} } x^{1-\nu_2-2N}|$, 
$|{e^{i\pi \nu_1} \over 16^{\nu_2+2N}} x^{\nu_2+2N}|$, 
$|v(x)|$, $|w(x)|$ are less than $\epsilon$,
 in such a way that $\Phi$ and $\Psi$ are holomorphic. 
That such a domain is not empty will be shown below. In particular, 
we'll show that if we require that  
 $|x|<r$, $|{e^{-i\pi \nu_1}\over 16^{1-\nu_2}} x^{1-\nu_2}|<r$, 
$|{e^{i\pi \nu_1} \over 16^{\nu_2}} x^{\nu_2}|<r$, where $r<\epsilon$ is 
small enough,  
also  the solutions of the integral equations $|v(x)|$ and $|w(x)|$ are  
less than $\epsilon$. Such a domain 
is contained in 
(\ref{dominiodim1}). 

\vskip 0.2 cm 

We choose the path of integration $L(x)$ connecting 0 to $x$, defined by
$$
\arg(s)=\arg(x) +{\Re \nu_2+2N-\nu^* \over \Im \nu_2} \log{|s|\over |x|},
~~~0<|s|\leq |x|,~~~\nu^{*}\in{\bf C}
$$
If $x$ belongs to the domain (\ref{dominiodim1}) in case (a) and to 
(\ref{dominiodim2}) in case (b),  
than  the path  does not leave the domain  when $s\to 0$, provided that 
$  0 \leq \nu^* \leq 1$ in case (a), or $
1 \leq \nu^* \leq 2 
$ in case (b).  
If $\Im \nu_2=0$ we take the radial path $\arg s = \arg x$, $0<|s|\leq |x|$, 
  namely $
\nu^*=\nu_2$. 
The parameterization of the path is 
$$
   s=\rho ~e^{i\left\{ \arg x + {\Re \nu_2+
2N - \nu^* \over \Im \nu_2} \log {\rho 
\over |x|}\right\}},~~~~0<\rho\leq|x|
$$
therefore
$$
|ds|=P(\nu_2+2N,\nu^*)~d\rho,~~~~~P(\nu_2+2N,\nu^*):= \sqrt{1+\left(
 {\Re \nu_2+2N- \nu^* \over \Im \nu_2}\right)^2}
$$
For any complex numbers $A$, $B$ we have
\be
\int_{L(x)} {1\over |s|} \left( |s|+ |A s^{1-\nu_2-2N}|+ |B s^{\nu_2+2N}| 
\right)^n 
~|ds| 
 \leq 
{P(\nu_2+2N,\nu^*) \over n \min(\nu^*,1-\nu^*)} \left(|x|+ 
|A x^{1-\nu_2-2N}|+ |B x^{\nu_2+2N}| \right)^n
\label{CONDIZRIMS1}
\ee
The quantity   $\min(\nu^*,1-\nu^*)$ is substituted by  
$\min(\nu^*-1,2-\nu^*)$ in case (b).  To prove  (\ref{CONDIZRIMS1}) 
 we observe  that on $L(x)$ we have
$$
   |s^{\nu_2+2N+\alpha}|=|x^{\nu_2+2N+\alpha}| {|s|^{\nu^*+\alpha} \over 
|x|^{\nu^*+\alpha}},~~~~\forall\alpha\in{\bf C}.
$$
Therefore
$$
\int_{L(x)} {1\over |s|} |s|^i|As^{1-\nu_2-2N}|^j |Bs^{\nu_2+2N}|^k  ~|ds| 
=
$$
$$ 
  { |A x^{1-\nu_2-2N}|^j~|Bx^{\nu_2+2N}|^k \over
|x|^{(1-\nu^*)j} ~|x|^{\nu^* k}} 
~P(\nu_2+2N,\nu^*)\int_0^{|x|} d\rho~\rho^{i-1+(1-\nu^*)j+\nu^* k} 
$$
$$
= {P(\nu_2+2N,\nu^*)\over i+j(1-\nu^*)+k\nu^*} 
|x|^i ~|A x^{1-\nu_2-2N}|^j~|B x^{\nu_2+2N}|^k
$$
$$\leq 
 {P(\nu_2+2N,\nu^*)\over (i+j+k) \min(\nu^*,1-\nu^*)} 
|x|^i ~|A x^{1-\nu_2-2N}|^j~|B x^{\nu_2+2N}|^k
$$
from which (\ref{CONDIZRIMS1}) follows, provided that $0<\nu^*<1$. For 
$\Im \nu_2=0$, this yields again $0<\nu_2<1$.

\vskip 0.2 cm 
We observe that a solution of the integral equations is also a solution of the 
differential equations, by virtue of the following  lemma:

\vskip 0.2 cm 
\noindent
{\bf Lemma 1:} {\it Let $f(x)$ be a holomorphic function in the domain 
$|x|<\epsilon$,  $|A x^{1-\nu_2-2N}|<\epsilon$, 
$|B x^{\nu_2+2N}|<\epsilon$, such that 
$f(x)=O(|x|+|Ax^{1-\nu_2-2N}|+|Bx^{\nu_2+2N}|)$, $A,B\in {\bf C}$. Let $L(x)$ 
be the path of integration  defined above for $0<\nu^*<1$ and 
$$
 F(x):= \int_{L(x)} {1\over s} f(s) ~ds
$$
Then, $F(x)$ is holomorphic on the domain and ${dF(x)\over dx}= {1\over x} f(x)$

The same holds for a function 
$f(x)=O(|x|+|Ax^{2-\nu_2-2N}|+|Bx^{\nu_2-1+2N}|)$ in 
the domain $|x|<\epsilon$,  $|A x^{2-\nu_2-2N}|<\epsilon$, 
$|B x^{\nu_2-1+2N}|<\epsilon$ ( $L(x)$ being defined for   $1<\nu^*<2$). 
}

\vskip 0.2 cm
\noindent
{\it Proof:} We give the proof in case (a). 
We choose a point $x+\Delta x$ close to $x$ and we prove 
that $\int_{L(x)}- \int_{L(x+\Delta x)} = \int_x^{x+\Delta x} $, where 
the last integral is on a segment from $x$ to $x+\Delta x$. 
 We consider a small disk $U_R$ centered at $x=0$ of small radius $R<|x|$ 
and the points $x_R:= L(x)\cap U_R$, $x^{\prime}_R :=
  L(x+\Delta x)\cap U_R$. Since  the integral of $f/s$ on a finite 
close curve (not containing 0) is zero we have:
\be
 \left( \int_{L(x)}-\int_{L(x+\Delta x)} - \int_x^{x+\Delta x}\right)ds ~{f(s) 
\over s}
=\left( \int_{L(x_R)}-\int_{L(x^{\prime}_R)} - \int_{\gamma(x_R,x^{\prime}_R)}
\right)ds ~{f(s) 
\over s}
\label{provcopia}
\ee
The last integral is on the arc $\gamma(x_R,x^{\prime}_R)$ 
   from $x_R$ to $x^{\prime}_{R}$ on the circle $|s|=R$. 
  We have also kept into account the obvious fact that 
$L(x_R)$ is contained in $L(x)$ and $L(x^{\prime}_R)$ is contained in $L(x+
\Delta x)$. We then take $R\to 0$ and we prove that the r.h.s. of 
(\ref{provcopia}) vanishes.

Taking into account that $f(x)=O(|x|+|Ax^{1-\nu_2-2N}|+|Bx^{\nu_2+2N}|)$ and 
(\ref{CONDIZRIMS1}) we have 
$$ 
 \left| \int_{L(x_R)} {1\over s} f(s) ~ds\right|\leq \int_{L(x_R)} 
{1\over |s|} O(|s| +|As^{1-\nu_2-2N}|+|Bs^{\nu_2+2N}| )~ |ds|
$$
$$
\leq {P(\nu_2+2N,\nu^*)  \over 
\min(\nu^*,1-\nu^*)} O(|x_R| + |Ax_R^{1-\nu_2-2N}|+|Bx_R^{\nu_2+2N}|)= 
$$
$$ = 
{P(\nu_2+2N,\nu^*)
  \over 
\min(\nu^*,2-\nu^*)}~O(R^{\min\{\nu^*,1-\nu^*\}})
$$
The last step follows from  $|x_R^{\nu_2+2N}|= {|x^{\nu_2+2N}| 
\over |x|^{\nu^*}}
 R^{\nu^*}$. So the integral vanishes for $R\to 0$. The same is proved for 
$\int_{L(x+\Delta x)}$. As for the integral on the arc we have 
$$
 |\arg x_R-\arg x^{\prime}_R|= \left| \arg x -\arg(x+\Delta x) +{\Re \nu_2 
+2N -
\nu^* \over \Im \nu_2} \log\left| 1 +{\Delta x\over x}\right|
\right|$$
or $
 |\arg x_R-\arg x^{\prime}_R|= \left| \arg x -\arg(x+\Delta x)\right|$ if $\Im \nu_2=0$. This is independent of $R$, 
therefore the length of the arc is $O(R)$ and 
$$
 \left|\int_{\gamma(x_R,x^{\prime}_R)}{1\over| s|} |f(s)|~|ds| \right|= O(R^{
\min\{\nu^*,1-\nu^*\}}) \to 0 \hbox{ for } R\to 0
$$
\rightline{$\Box$}

\vskip 0.2 cm
Now we prove a fundamental lemma.

\vskip 0.2 cm 
\noindent
{\bf Lemma 2:} 
{\it 
For any $N\in {\bf Z}$ and for  any complex $\nu_1$,  $\nu_2$  such that 
$$\nu_2 \not\in (-\infty,-2N]\cup\{1-2N\}\cup[2-2N,+\infty)$$ 
there exists a sufficiently small
$r_N<1$ such that 
the system of integral equations 
 has a solution $v_1(x;\nu_1,\nu_2+2N)$ 
holomorphic in the domain ${\cal D}_1(r_N;\nu_1,\nu_2+2N)$ defined in 
(\ref{DOMINOO1}). 
 Moreover, there exists a  positive constant $M_1(\nu_2+2N)$ 
depending on $\nu_2+2N$ 
such that 
 $$
|v_1(x;\nu_1,\nu_2+2N)|\leq M_1(\nu_2+2N) 
\left(|x|+\left|{e^{-i\pi \nu_1}\over 16^{1-\nu_2-2N}} 
x^{1-\nu_2-2N} \right|+\left| 
{e^{i\pi \nu_1} \over 16^{\nu_2+2N}} x^{\nu_2+2N}\right| \right)
$$
 in  ${\cal
D}_1(r;\nu_1,\nu_2)$. 

\vskip 0.2 cm 
 The system of integral equations has {\rm another} solution $v_2(x;\nu_1,\nu_2+2N)$ 
 holomorphic 
in ${\cal
D}_2(r;\nu_1,\nu_2)$ defined in (\ref{DOMINOO2}). 
There exists a positive constant  $M_2(\nu_2+2N)$ such 
that 
 $$
|v_2(x;\nu_1,\nu_2+2N)|\leq M_2(\nu_2+2N) 
\left(|x|+\left|{e^{-i\pi \nu_1}\over 16^{2-\nu_2-2N}} 
x^{2-\nu_2-2N} \right|+\left| 
{e^{i\pi \nu_1} \over 16^{\nu_2-1+2N}} x^{\nu_2-1+2N}\right| \right)
$$
 in 
 ${\cal
D}_2(r;\nu_1,\nu_2)$. 
}

\vskip 0.2 cm 

Note that  ${\cal D}_i(r_N;\nu_1,\nu_2+2N) = {\cal D}_i(r_N;\nu_1+2N_1,
\nu_2+2N)$, $i=1,2$, for any $N_1\in {\bf Z}$.  
To prove Lemma 2 we need some sub-lemmas
\vskip 0.2 cm
\noindent
{\bf Sub-Lemma 1:} 
{\it Let $\Phi(x,y,z)$ and $\Psi(x,y,z,v,w)$ be two holomorphic 
functions of their arguments for $|x|,|y|,|z|,|v|,|w|< \epsilon$, satisfying
$$
   \Phi(0,0,0)=0,~~~~\Psi(0,0,0,v,w)=\Psi(x,y,z,0,0)=0
$$
Then, there exists a constant $c>0$ such that:
\be 
     \left| \Phi(x,y,z) \right|\leq c~\left(|x|+|y|+|z|\right)
\label{primacopia1}
\ee
\be
   \left| \Psi(x,y,z,v,w) \right|\leq c~\left(|x|+|y|+|z|\right)
\label{primacopia2}
\ee
\be 
\left| \Psi(x,y,z,v_2,w_2)- \Psi(x,y,z,v_1,w_1)\right|\leq c
\left(|x|+|y|+|z|\right)~\left( |v_2-v_1|+|w_2-w_1|\right)
\label{primacopia3}
\ee
 for $|x|,|y|,|z|,|v|,|w|< \epsilon$.
}

\vskip 0.2 cm
\noindent
{\it Proof:} Let's proof (\ref{primacopia2}). 
$$
 \Psi(x,y,z,v,w)=\int_0^1 ~{d\over d\lambda}
\Psi(\lambda x, \lambda y, \lambda z,v,w)~d\lambda
$$
$$
  x~\int_0^1 {\partial \Psi \over \partial x}(\lambda x, \lambda y, 
\lambda z,v,w) ~d\lambda~+~y\int_0^1 {\partial \Psi \over \partial y}
(\lambda x, \lambda y, 
\lambda z,v,w) ~d\lambda~+~z\int_0^1 {\partial \Psi \over \partial z}
(\lambda x, \lambda y, 
\lambda z,v,w) ~d\lambda
$$
Moreover, for $\delta$ small:
$$ 
{\partial \Psi \over \partial x}(\lambda x, \lambda y, 
\lambda z,v,w)= \int_{|\zeta-\lambda x|=\delta}~
{\Psi(\zeta,\lambda y,\lambda z,
v,w) \over (\zeta-\lambda x)^2}~{d\zeta \over 2 \pi i} 
$$
which implies that ${\partial \Psi \over \partial x}$ is holomorphic and 
bounded when its arguments are less than $\epsilon$. The same holds true for 
 ${\partial \Psi \over \partial y}$ and  ${\partial \Psi \over \partial z}$. This proves (\ref{primacopia2}), $c$ being a constant which bounds 
$\left|{\partial \Psi \over \partial x}\right|$, $\left|
{\partial \Psi \over \partial y}\right|$ $\left|{\partial \Psi 
\over \partial z}\right|$. The inequality (\ref{primacopia1}) is proved in 
the same way.  
We turn to (\ref{primacopia3}). First we prove that for 
$|x|,|y|,|z|,|v_1|,|w_1|,|v_2|,|w_2| < \epsilon$ there exist two 
holomorphic and bounded functions $\psi_1(x,y,z,v_1,w_1,v_2,w_2)$,  
$\psi_2(x,y,z,v_1,w_1,v_2,w_2)$ such that
$$
    \Psi(x,y,z,v_2,w_2)-\Psi(x,y,z,v_1,w_1)
$$
\be
= (v_2-v_1)~ \psi_1(x,y,z,v_1,w_1,v_2,w_2)+ (w_1-w_2)~
\psi_2(x,y,z,v_1,w_1,v_2,w_2)
\label{tmpRIMS1}
\ee 
In order  to prove this, we write 
$$
 \Psi(x,y,z,v_2,w_2)-\Psi(x,y,z,v_1,w_1)=
$$
$$
=
\int_0^1 {d\over d\lambda} \Psi(x,y,z,\lambda v_2 +(1-\lambda)v_1,
\lambda w_2 +(1-\lambda)w_1)~ d\lambda
$$
$$
=(v_2-v_1) \int_0^1 {\partial \Psi \over \partial v} 
 (x,y,z,\lambda v_2 +(1-\lambda)v_1,
\lambda w_2 +(1-\lambda)w_1)~d\lambda~+
$$
$$
+~(w_2-w_1) \int_0^1
 {\partial \Psi \over \partial w} 
 (x,y,z,\lambda v_2 +(1-\lambda)v_1,
\lambda w_2 +(1-\lambda)w_1)~d\lambda
$$
$$
=: (v_2-v_1)~ \psi_1(x,y,z,v_1,w_1,v_2,w_2)+ (w_2-w_1) ~
\psi_2(x,y,z,v_1,w_1,v_2,w_2)
$$
Moreover, for small $\delta$,  
$$
{\partial \Psi \over \partial v}(x,y,z,v,w)= \int_{|\zeta-v|=\delta} 
{\Psi(x,y,z,\zeta,w) \over (\zeta-v)^2} ~{dz\over 2\pi i} 
$$
which implies that $\psi_1$ is holomorphic and bounded for its arguments less 
 than $\epsilon$. We also obtain 
${\partial \Psi \over \partial v}(0,0,0,v,w)=0$, then  
$\psi_1(0,0,0,v_1,w_1,v_2,w_2)=0$. The proof for $\psi_2$ is analogous.
 We use (\ref{tmpRIMS1}) to complete the proof of (\ref{primacopia3}). 
Actually, we observe that 
$$ 
  \psi_i(x,y,z,v_1,w_1,v_2,w_2)= \int_0^1 {d \over d\lambda} 
\psi_i(\lambda x,\lambda y, \lambda z, v_1,w_1,v_2,w_2) ~d\lambda 
$$
$$
 = x\int_0^1 {\partial \psi_i \over \partial x} d\lambda + 
 y\int_0^1 {\partial \psi_i \over \partial y} d\lambda + 
z\int_0^1 {\partial \psi_i \over \partial z} d\lambda 
$$
and we conclude as in the proof of (\ref{primacopia2}).

\rightline{$\Box$}

\vskip 0.2 cm

We solve the system of integral equations 
by successive approximations.  We can choose any path $L(x)$ 
such that $0<\nu^*<1$ in case (a), or $1<\nu^*<2$ in case (b). 
Here we complete the proof only for the case (a).  Therefore, the solution $v(x)$ we are going to find is $v_1(x;\nu_1,\nu_2+2N)$. Case (b) is analogous and we don't need to repeat the proof for it. It yields the function   
$v_2(x;\nu_1,\nu_2+2N)$.

 We choose $\nu^*={1\over 2}$, therefore $\min\{ \nu^*, 1-\nu^*\}={1\over 2}$.
   For convenience, we put  
$$
A:= {e^{-i\pi \nu_1}\over 16^{1-\nu_2-2N}},~~~~
B:= {e^{i\pi \nu_1} \over 16^{\nu_2+2N}}
$$
Therefore, for any $n\geq 1$ the successive approximations are: 
$$
    v_0=w_0=0
$$
\vskip 0.15 cm 
\be
w_n(x)=\int_{L(x)} {1\over t} \left\{
 \Phi (s,A s^{1-\nu_2-2N}, 
B s^{\nu_2+2N} ) + \Psi(s,A s^{1-\nu_2-2N}, B s^{\nu_2+2N},v_{n-1}(s),w_{n-1}(s))
\right\}~ds
\label{RIMSRIMS1}
\ee
\be
v_n(x)=\int_{L(x)}{1\over s}~w_{n}(s)~ ds
\label{RIMSRIMS2}
\ee

\vskip 0.2 cm
\noindent
{\bf Sub-Lemma 2:} {\it There exists a sufficiently small $\epsilon^{\prime}
<\epsilon$ such that for 
any $n\geq 0$  the functions 
$v_n(x)$ and $w_n(x)$ are holomorphic   in the domain 
$$
   {\cal D}_1(\epsilon^{\prime};\nu_1,\nu_2+2N):= 
\left\{ x\in \tilde{\bf C}_0~ 
\hbox{ such that }
  |x|<\epsilon^{\prime}, \left|A x^{1-\nu_2-2N} \right|<\epsilon^{\prime},
\left| B  x^{\nu_2+2N}\right|<\epsilon^{\prime} \right\}
$$
They are also correctly bounded, namely $|v_n(x)|<\epsilon$, 
$|w_n(x)|<\epsilon$ for any
 $n$. They  satisfy 
\be 
 |v_n-v_{n-1}| \leq { (2c)^n (2P(\nu_2+2N))^{2n} \over n!} \left(
|x|+|Ax^{1-\nu_2-2N}|+|Bx^{\nu_2+2N}|
\right)^n
\label{InEqAlItY1}
\ee
\be
 |w_n-w_{n-1}| \leq { (2c)^n (2P(\nu_2+2N))^{2n} \over n!} \left(
|x|+|Ax^{1-\nu_2-2N}|+|Bx^{\nu_2+2N}|
\right)^n
\label{InEqAlItY}
\ee
where $P(\nu_2+2N) :=P(\nu_2+2N,\nu^*=1/2)$ and $c$ is the constant appearing in 
Sub-Lemma 1. Moreover 
$$
    x {dv_n \over dx}=w_n
$$
}

\vskip 0.2 cm 
\noindent
{\bf Proof:} We proceed by induction. 
$$ 
  w_1 = \int_{L(x)} {1\over s} 
 \Phi (s,A s^{1-\nu_2-2N}, 
B s^{\nu_2+2N} ) ~ds,
~~~~v_1=\int_{L(x)}{1\over s} w_1(s) ds 
$$
 It follows from Lemma 1 and (\ref{primacopia1}) that $w_1(x)$ is 
holomorphic for $|x|,|Ax^{1-\nu_2-2N}|, |Bx^{\nu_2+2N}|<\epsilon$. 
From (\ref{CONDIZRIMS1})  and (\ref{primacopia1}) we have
$$
   |w_1(x)|\leq \int {1\over |s| } |\Phi(s,As^{1-\nu_2-2N},Bs^{\nu_2+2N})|
~|ds| 
$$
$$
\leq 2c P(\nu_2+2N) 
(|x|+|Ax^{1-\nu_2-2N}|+|Bx^{\nu_2+2N}|) 
\leq 6cP(\nu_2+2N) \epsilon^{\prime}<\epsilon
$$
 on 
 $ {\cal D}_1(\epsilon^{\prime};\nu_1,\nu_2+2N)$, provided 
that $\epsilon^{\prime}$ is small enough.   
By Lemma 1, also $v_1(x)$ is holomorphic   
for $|x|,|Ax^{1-\nu_2-2N}|, |Bx^{\nu_2+2N}|<\epsilon$ and
$$
 x {d v_1 \over dx} =w_1
$$
By (\ref{CONDIZRIMS1}) we also have 
$$
   |v_1(x) | \leq c(2P(\nu_2+2N))^2  (|x|+|Ax^{1-\nu_2-2N}|+|Bx^{\nu_2+2N}|)
\leq 12cP(\nu_2)^2 \epsilon^{\prime}<\epsilon
$$ 
on $ {\cal D}_1(\epsilon^{\prime};\nu_1,\nu_2+2N)$, if $\epsilon^{\prime}$ is 
small enough.  
Note that $P(\nu_2+2N)\geq 1$, so (\ref{InEqAlItY}) (\ref{InEqAlItY1}) are
 true for $n=1$. Now we 
 suppose that the statement of the sub-lemma 
is true for $n$ and we prove it for $n+1$. Consider: 
$$
|w_{n+1}(x) - w_n(x)| = \left| 
 \int_{L(x)}{1\over s}~ \left[
            \Psi(s,As^{1-\nu_2-2N},Bs^{\nu_2+2N},v_n,w_n)- \right. 
\right.
$$
$$ \left. \left. -
\Psi(s,As^{1-\nu_2-2N},Bs^{\nu_2+2N},v_{n-1},w_{n-1}) 
\right] 
~ds
\right|
$$
By (\ref{primacopia3}) the above is  
$$
\leq 
     c~\int_{L(x)} {1\over |s|} (|s|+|As^{1-\nu_2-2N}|+|Bs^{\nu_2+2N}|)~
(|v_n-v_{n-1}| + |w_n-w_{n-1}|) ~|ds|
$$
By induction this is 
$$
 \leq 2c~ {(2c)^n (2P(\nu_2+2N))^{2n} \over n!} ~ \int_{L(x)} {1\over |s|} 
(|s|+|As^{1-\nu_2-2N}|+|Bs^{\nu_2+2N}|)^{n+1} ~|ds|
$$
$$
\leq 
      2c~  {(2c)^n (2P(\nu_2+2N))^{2n} \over n!} ~ {2P(\nu_2)\over n+1} 
~(|x|+|Ax^{1-\nu_2-2N}|+|Bx^{\nu_2+2N}|)^{n+1}
$$
$$
\leq 
{(2c)^{n+1} (2P(\nu_2+2N))^{2(n+1)} \over (n+1)!}
(|x|+|Ax^{1-\nu_2-2N}|+|Bx^{\nu_2+2N}|)^{n+1}
$$
This proves (\ref{InEqAlItY}). Now we estimate 
$$
   |v_{n+1}(x)-v_n(x)|\leq \int_{L(x)} |w_{n+1}(s)-w_n(s)|~|ds|
$$
$$
\leq {(2c)^{n+1} (2P(\nu_2+2N))^{2n+1} \over (n+1)!}~\int_{L(x)} 
{1\over |s| } 
(|s|+|As^{1-\nu_2-2N}|+|Bs^{\nu_2+2N}|)^{n+1}~|ds|
$$
$$
\leq 
{(2c)^{n+1} (2P(\nu_2+2N))^{2(n+1)} \over (n+1)~(n+1)!}~
(|x|+|Ax^{1-\nu_2-2N}|+|Bx^{\nu_2+2N}|)^{n+1}
$$
$$
\leq 
{(2c)^{n+1} (2P(\nu_2+2N))^{2(n+1)} \over ~(n+1)!}~
(|x|+|Ax^{1-\nu_2-2N}|+|Bx^{\nu_2+2N}|)^{n+1}
$$
This proves (\ref{InEqAlItY1}). From Lemma 1 we also conclude that 
$w_n$ and $v_n$ are holomorphic in ${\cal D}(\epsilon^{\prime},\nu_1,\nu_2)$ 
and 
$$
  x{dv_n \over dx}=w_n
$$
Finally we see that  
$$ 
|v_n(x)|\leq \sum_{k=1}^n |v_k(x)-v_{k-1}(x)|\leq \exp\{ 2c(2P(\nu_2+2N))^2
 (|x|+|Ax^{1-\nu_2-2N}|+|Bx^{\nu_2+2N}|)\}-1 \leq
$$
$$\leq 
 \exp\{ 24cP^2(\nu_2+2N)
                                  \epsilon^{\prime}\}-1
$$
and the same for $|w_n(x)|$. Therefore, if $\epsilon^{\prime}$ is small 
enough  we have $|v_n(x)|<\epsilon$, $|w_n(x)|<\epsilon$ on 
 ${\cal D}_1(\epsilon^{\prime},\nu_1,\nu_2)$. 

Note that $P(\nu_2+2N)$ grows as $N^2$, so $\epsilon^{\prime}$ 
decreases as $N^{-2}$. 

\rightline{$\Box$}

\vskip 0.2 cm

Let's define 
$$
 v(x):=\lim_{n\to \infty} v_n(x),~~~~~w(x):=\lim_{n\to \infty} w_n(x) 
$$
if they exist. 
We can also rewrite
$$
 v(x) = \lim_{n\to \infty} v_n(x)= \sum_{n=1}^{\infty} (v_n(x)-v_{n-1}(x)). 
$$
 We see that the series converges uniformly in 
 ${\cal D}(\epsilon^{\prime},\nu_1,\nu_2+2N)$ because 
$$
|\sum_{n=1}^{\infty} (v_n(x)-v_{n-1}(x)) | 
$$
$$
\leq 
\sum_{n=1}^{\infty} 
{(2c)^{n} (2P(\nu_2+2N))^{2n} \over ~n!}~
(|x|+|Ax^{1-\nu_2-2N}|+|Bx^{\nu_2+2N}|)^{n}
$$
$$
 = \exp\{ 8 cP^2(\nu_2)
 (|x|+|Ax^{1-\nu_2-2N}|+|Bx^{\nu_2+2N}|) \} -1
$$
The same holds 
for $w_n(x)$. Therefore, $v(x)$ and $w(x)$ define holomorphic functions in  
 ${\cal D}_1(\epsilon^{\prime},\nu_1,\nu_2+2N)$. From Sub-Lemma 2 
we also have 
$$
    x {d v(x) \over dx } = w(x)
$$
in ${\cal D}_1(\epsilon^{\prime},\nu_1,\nu_2+2N)$. 

\vskip 0.2 cm
We show that $v(x),w(x)$ solve the initial integral equations. The l.h.s. of 
(\ref{RIMSRIMS1}) converges to $w(x)$ for $n\to \infty$. Let's prove that the 
r.h.s.  also converges to 
$$\int_{L(x)} {1\over s} \left\{
 \Phi (s,A s^{1-\nu_2-2N}, 
B s^{\nu_2+2N} ) + \Psi(s,A s^{1-\nu_2-2N}, B s^{\nu_2+2N},v(s),w(s))
\right\}~ds. $$ 
 We have to evaluate the following difference: 
$$
   \left|
\int_{L(x)} {1\over s} \Psi(s,A s^{2-\nu_2-2N}, B s^{\nu_2+2N},v(s),w(s))
~ds ~-~ 
\int_{L(x)} {1\over s} \Psi(s,A s^{2-\nu_2-2N}, B s^{\nu_2+2N},v_n(s),w_n(s))
~ds
\right|
$$
By (\ref{primacopia3}) the above is 
\be
\leq ~c~ \int_{L(x)} {1\over |s|} (|s|+|As^{1-\nu_2-2N}|+|Bs^{\nu_2+2N}|)~(
|v-v_n|+|w-w_n|) ~|ds|
\label{palla}
\ee
Now we observe that
$$
|v(x)-v_n(x)| \leq \sum_{k=n+1}^{\infty} |v_k-v_{k-1}| 
$$
$$
  = \sum_{k=n+1}^{\infty} {(2c)^k (2P(\nu_2+2N))^{2k} \over k!} 
(|x|+|Ax^{1-\nu_2-2N} |
+ |Bx^{\nu_2+2N}| )^k
$$
$$
\leq 
     (|x|+|Ax^{2-\nu_2-2N} |
+ |Bx^{\nu_2+2N}| )^{n+1} 
  ~\sum_{k=0}^{\infty} {(2c)^{k+n+1} P(\nu_2+2N)^{2(k+n+1)} \over (k+n+1)!} 
(|x|+|Ax^{2-\nu_2-2N} |
+ |Bx^{\nu_2+2N}| )^k
$$
The series converges. Its sum is less than some constant $S(\nu_2+2N)$ 
 independent of $n$. We obtain 
$$
  |v(x)-v_n(x)| \leq S(\nu_2+2N)   (|x|+|Ax^{1-\nu_2-2N} |
+ |Bx^{\nu_2+2N}| )^{n+1}. $$
The same holds for $|w-w_n|$. Thus, (\ref{palla}) is 
$$
\leq 
  2c~S(\nu_2+2N) 
~\int_{L(x)} {1\over |s| } (|s|+|As^{1-\nu_2-2N}|+|Bs^{\nu_2+2N}|)^{n+2}
~|ds| $$
$$
\leq 
  {2cS(\nu_2+2N) ~2P(\nu_2+2N)\over n+2}~  (|x|+|Ax^{1-\nu_2-2N}|+
|Bx^{\nu_2+2N}|)^{n+2}
$$
Namely: 
$$
   \left|
\int_{L(x)} {1\over s} \Psi(s,A s^{1-\nu_2-2N}, B s^{\nu_2+2N},v(s),w(s))
~ds ~-~ 
\int_{L(x)} {1\over s} \Psi(s,A s^{1-\nu_2-2N}, B s^{\nu_2+2N},v_n(s),w_n(s))
~ds
\right|
$$
$$
\leq 
{2cS(\nu_2) ~2P(\nu_2+2N)\over n+2}~(3\epsilon^{\prime})^{n+2}
$$
In a similar way, the r.h.s. of (\ref{RIMSRIMS2}) is 
$$
    \left|\int {1\over s} (w(s)-w_n(s))~ds\right|\leq {S(\nu_2+2N)~
2 P(\nu_2) \over n+1} 
(3\epsilon^{\prime})^{n+1}
$$
Therefore, the r.h. sides  of  (\ref{RIMSRIMS1}) (\ref{RIMSRIMS2}) 
converge on the domain ${\cal D}(r,\nu_1,\nu_2)$ for 
$r< \min\{\epsilon^{\prime}, 1/3\}$. We note that $r=r_N$, namely it 
depends on $N$ in the same way as  $\epsilon^{\prime}$ does. Thus,   it  
decreases as $N^{-2}$ as $N$ increases. 
 
We finally observe that $|v(x)|$ and $|w(x)|$ are bounded on ${\cal D}(r)$. 
For example  
$$
|v(x)| \leq  (|x|+|Ax^{1-\nu_2-2N} |
+ |Bx^{\nu_2+2N}| ) 
  ~\sum_{k=0}^{\infty} {(2c)^{k+1} (2P(\nu_2))^{2(k+1)} \over (k+1)!} 
(|x|+|Ax^{1-\nu_2-2N} |
+ |Bx^{\nu_2+2N}| )^k
$$
$$
=: M(\nu_2+2N)  (|x|+|Ax^{1-\nu_2-2N} |
+ |Bx^{\nu_2+2N}| )
$$
where  the sum of the series is less than a constant $M(\nu_2+2N)$. 
We have proved Lemma 2 for the case (a). Case (b) is analogous. 

\rightline{$\Box$}

\vskip 0.3 cm 
\noindent 
{\bf Remark:} 
 The domain  ${\cal D}_1(r_N,\nu_1,\nu_2)$ is   
$$   (\Re \nu_2+2N) \ln{|x|\over 16} -\pi \Im \nu_1
-\ln r_N < \Im \nu_2
   \arg x <
$$
$$ < (\Re \nu_2-1 +2N)\ln{|x|\over 16}-\pi \Im \nu_1
+\ln r_N,~~~~|x|<r_N
$$
 while 
the domain  ${\cal D}_2(r_N,\nu_1,\nu_2)$ is 
$$   (\Re \nu_2-1+2N) \ln{|x|\over 16}-\pi \Im \nu_1
-\ln r_N < \Im \nu_2
   \arg x < $$
$$ 
< (\Re \nu_2 -2+2N)\ln{|x|\over 16}-\pi \Im \nu_1
+\ln r_N, ~~~~|x|<r_N
$$
 We invite  the reader to observe the $\ln r_N$ terms.  

Also, note that for $\nu_2\in {\bf R}$ the domain is specified by 
$\max\left\{ |x|,   \left|{e^{-i\pi \nu_1}\over 16^{1-\nu_2}} x^{1-\nu_2} 
\right|,
\left| 
{e^{i\pi \nu_1} \over 16^{\nu_2}} x^{\nu_2}\right|\right\}<r$. We can make 
the 
notation easer by simply re-writing $|x|<r$ 
for sufficiently small $r$. 
So, we obtain the domain ${\cal D}_0(r)$.

\vskip 0.3 cm 

We have proved that 
the structure of the integral equations implies that $v(x)$ is bounded 
(namely $|v(x)|=O(r_N)$). 
The proof of Lemma 2 only makes use of the properties of 
$\Phi$ and $\Psi$, regardless of how these functions have been constructed. 

 In our case, 
 $\Phi$ and $\Psi$ have been 
constructed from the Fourier expansion of elliptic functions. 
We see that the domain 
 (\ref{dominiodim1})
 contains ${\cal D}_1(r_N;\nu_1,\nu_2+2N)$ and (\ref{dominiodim2}) contains  ${\cal D}_2(r_N;\nu_1,\nu_2+2N)$ because 
  the term $\Im {\pi v \over 2 \omega_1}$  in (\ref{dominiodim1}) and 
(\ref{dominiodim2}) is 
$O(r_N)$, while in ${\cal D}_i(r,\nu_1,\nu_2+2N)$, $i=1,2$,  
the term $\ln r_N$ appear (see the Remark).   

\vskip 0.2 cm 
 To conclude the proof of Theorem 1, we have to work out the series 
 of $v(x)$. Let's do that in case (a), namely for $v_1(x;\nu_1,\nu_2+2N)$. 
Case (b) is analogous.   We observe that $w_1$ and $v_1$ 
are series of the type 
\be 
  \sum_{p,q,r\geq 0} c_{pqr}(\nu_2+2N) ~x^p~(Ax^{1-\nu_2-2N})^q~
(Bx^{\nu_2+2N})^r
\label{formMrims}
\ee
where $c_{pqr}(\nu_2+2N)$ is rational in $\nu_2$. This follows from 
$$
   w_1(x)= \int_{L(x)} \Phi(s,As^{1-\nu_2-2N},Bs^{\nu_2+2N})~ds
$$
and from the fact that $\Phi(x,Ax^{1-\nu_2-2N},Bx^{\nu_2+2N})$ itself 
 is a series 
 (\ref{formMrims}) by construction, with  coefficients 
which are  rational functions of  $\nu_2+2N$.  
The same holds true for  $\Psi$.  We conclude that 
  $w_n(x)$ and $v_n(x)$ have the form (\ref{formMrims}) for any $n$. This 
implies that the limit $v(x)$ is also a series of type (\ref{formMrims}). We 
can reorder such a series. Consider the term 
$$
   c_{pqr}(\nu_2) ~  x^p~(Ax^{1-\nu_2-2N})^q~(Bx^{\nu_2+2N})^r,$$
and recall that by definition $B={1\over  16^2 ~A}$. We absorb $16^{-2r}$ into 
$ c_{pqr}(\nu_2)$ and we study the factor 
$$
    A^{q-r} x^{p+(1-\nu_2-2N)q+(\nu_2+2N) r} =  A^{q-r} x^{p+q+(r-q)(\nu_2+2N)}
$$
We have three cases: 

1) $r=q$, then we have $x^{p+q}=: x^n$, $n=p+q$. 

2) $r>q$, then we have $x^{p+q}~\left[{1\over A} x^{\nu_2+2N}\right]^{r-q}=:
x^n~\left[{1\over A} x^{\nu_2+2N}\right]^m$, $n=p+q$, $m=q-r$.

3) $r<q$,  then we have $A^{q-r} ~x^{p+r} ~\left[ A x^{1-\nu_2-2N}\right]^{q-r}
 =: x^n~\left[ A x^{1-\nu_2-2N}\right]^m$, $n= {p+r}$, $m=q-r$. 

\noindent
Therefore, the series of the type  (\ref{formMrims}) can be re-written in the 
 form of $v(x)$ of 
 Theorem 1. The series of $v(x)$ is uniquely determined by the 
Painlev\'e equation, because it is constructed from $\Phi$ and $\Psi$ by 
successive approximations of the solutions of the integral equations.

\vskip 0.2 cm 
\noindent

This concludes the proof. We only remark that if $\nu_2$ is real, we make a  
  translation $\nu_2\mapsto \nu_2+2N$ which yields $0<\nu_2<1$ or $1<\nu_2<2$ 
and all the  formulae can be read for $N=0$. 
\vskip 0.2 cm 

 The statement of Theorem 1 as it follows from our proof is 

\vskip 0.3 cm 
\noindent
 { \it  Let $\nu_1$, $\nu_2$ be two complex numbers. 

\vskip 0.2 cm 
{\bf I)} For any $N\in {\bf Z}$ and for any 
complex  $\nu_1$, $\nu_2$ such that 
$$ 
\Im \nu_2\neq 0
$$
there exist a  positive number $r_N<1$  and a  transcendent 
$$
 y(x) = \wp \Bigl(\nu_1\omega_1(x)+\nu_2\omega_2(x) +v_1(x;\nu_1,\nu_2+2N);~ \omega_1(x),\omega_2(x)
 \Bigr)+{1+x\over 3}
$$
such that $v_1(x;\nu_1,\nu_2+2N)$ is holomorphic in the domain $ 
{\cal D}_1(r_N;\nu_1,\nu_2+2N)  
$,  
where it has  convergent expansion 
$$ 
v_1(x;\nu_1,\nu_2+2N)= \sum_{n\geq 1} a^{(1)}_n x^n + \sum_{n\geq 0, m\geq 1} b^{(1)}_{nm} x^n 
\left[ e^{-i\pi \nu_1} \left({x\over 16}\right)^{1-\nu_2-2N}\right]^m + 
$$
$$
+\sum_{n\geq 0,m\geq 1} c^{(1)}_{nm} x^n 
\left[ e^{i\pi \nu_1} \left({x\over 16}\right)^{\nu_2+2N}\right]^m 
$$ 

  There also exists a    transcendent 
$$
 y(x) = \wp \Bigl(\nu_1\omega_1(x)+\nu_2\omega_2(x) +
v_2(x;\nu_1,\nu_2+2N);~ \omega_1(x),\omega_2(x)
 \Bigr)+{1+x\over 3}
$$
such that $v_2(x;\nu_1,\nu_2+2N)$ is holomorphic in the domain 
 $ 
{\cal D}_2(r_N;\nu_1,\nu_2+2N)  
$,   
where it has  convergent expansion 
 $$ 
v_2(x;\nu_1,\nu_2+2N)= 
\sum_{n\geq 1} a^{(2)}_n x^n + \sum_{n\geq 0, m\geq 1} b^{(2)}_{nm} x^n 
\left[ e^{-i\pi \nu_1} \left({x\over 16}\right)^{2-\nu_2-2N}\right]^m 
+
$$
$$
+\sum_{n\geq 0,m\geq 1} c^{(2)}_{nm} x^n 
\left[ e^{i\pi \nu_1} \left({x\over 16}\right)^{\nu_2-1+2N}\right]^m 
$$
\vskip 0.2 cm 
For both transcendents, 
 the coefficients $a^{(i)}_n$, 
$b^{(i)}_{nm}$, $c^{(i)}_{nm}$, $i=1,2$, are rational functions of $\nu_2+2N$. 
On  ${\cal D}_i(r_N;\nu_1,\nu_2+2N)$ there exists a positive constant 
$M_i(\nu_2+2N)$ such that 
\be
  |v_1(x;\nu_1,\nu_2+2N) | \leq M_1(\nu_2+2N) \left(|x|+\left| e^{-i\pi \nu_1} \left({x\over 16}
\right)^{1-\nu_2-2N}\right| + \left|  e^{i\pi \nu_1} 
\left({x\over 16}\right)^{\nu_2+2N} \right| 
 \right)
\label{bundo1}
\ee
\be
  |v_2(x;\nu_1,\nu_2+2N) | \leq M_2(\nu_2+2N) \left(|x|+\left| e^{-i\pi \nu_1} \left({x\over 16}\right)^{2-\nu_2-2N}\right| + \left|  e^{i\pi \nu_1} \left({x\over 16}\right)^{\nu_2-1+2N} \right| 
 \right)
\label{bundo2}
\ee

\vskip 0.2 cm  

{\bf II)} For any complex $\nu_1$ and for any real $\nu_2$ 
such that 
$$
0<\nu_2<1~~~\hbox{ or }~~~ 1<\nu_2<2
$$
there exists a positive  $r<1$ and a  transcendent 
$$
 y(x) = \wp \Bigl(\nu_1\omega_1(x)+\nu_2\omega_2(x) +v(x;\nu_1,\nu_2);~ \omega_1(x),\omega_2(x)
 \Bigr)+{1+x\over 3}
$$
such that $v(x;\nu_1,\nu_2)$ is holomorphic in ${\cal D}_0(r;\nu_1,\nu_2)$,  
where it has the  convergent expansion  
$$ 
v(x;\nu_1,\nu_2)= \sum_{n\geq 1} a_n x^n + \sum_{n\geq 0, m\geq 1} b_{nm} x^n 
\left[ e^{-i\pi \nu_1} \left({x\over 16}\right)^{1-\nu_2}\right]^m 
+
$$
$$ 
+\sum_{n\geq0,m\geq 1} c_{nm} x^n 
\left[ e^{i\pi \nu_1} \left({x\over 16}\right)^{\nu_2}\right]^m 
$$
if  $0<\nu_2<1$; and 
 $$ 
v(x)= \sum_{n\geq 1} a_n x^n + \sum_{n\geq 0, m\geq 1} b_{nm} x^n 
\left[ e^{-i\pi \nu_1} \left({x\over 16}\right)^{2-\nu_2}\right]^m 
+
$$
$$+
\sum_{n\geq 0,m\geq 1} c_{nm} x^n 
\left[ e^{i\pi \nu_1} \left({x\over 16}\right)^{\nu_2-1}\right]^m 
$$
if  $1<\nu_2<2$. The coefficients are rational functions of $\nu_2$ and $|v(x)|$ is bounded as in (\ref{bundo1}), with $N=0$, for $0<\nu_2<1$ and as in   
(\ref{bundo2}), with $N=0$, for $1<\nu_2<2$. 
}

\vskip 0.3 cm 
 The above statement 
 is completely equivalent to the statement of Theorem 1 in the Introduction,  
due to  
Observation 1 and 2  there. 
In particular, $v_1(x;\nu_1,\nu_2)= v(x;\nu_1,\nu_2)$ of Theorem 1 and  $v_2(x;\nu_1,\nu_2)= - v(x;-\nu_1,2-\nu_2)$

\vskip 0.3 cm  
A technique similar to that used in the proof of Theorem 1 was first introduced by S. Shimomura in \cite{Sh} for a class of functional equations.


\section{ Proof of Theorem 2}\label{theorem2}

Theorem 2 is stated in the Introduction. Here we give its proof. 
 Let $x_0\in {\cal D}(r;\nu_1,\nu_2)$.   
We expand $y(x)$ in Fourier series. This is possible on 
$ {\cal D}(r;\nu_1,\nu_2)$   as it is 
clear from the proof of Theorem 1. Let $f:= 
  \nu_1 +\nu_2\tau +{v\over \omega_1}$. 
$$
  y(x) = \wp(\nu_1\omega_1+\nu_2\omega_2+v;\omega_1,\omega_2)+{1+x\over 3}
=
$$
$$
\left({\pi \over 2 \omega_1}\right)^2
\left\{ 
-{1\over 3} +\sin^{-2}\left( 
{\pi \over 2} f
\right)
+8  \sum_{n\geq 1} {n e^{2 i \pi  n \tau} \over 1 - e^{2 i \pi  n \tau} }
\left[ 
1-\cos
\left(
\pi n 
f 
\right)
\right]
\right\}
              +{1+x\over 3}=
$$
$$
=\left({\pi \over 2 \omega_1}\right)^2
\left\{ 
-{1\over 3} - 4 
\left[
e^{i{\pi \over 2}f}- 
e^{-i{\pi \over 2}f}
\right]^{-2}
 +
4 \sum_{n\geq 1} {n e^{2 i \pi  n \tau} \over 1 - e^{2 i \pi  n \tau} }
\left[ 
2-e^{i\pi n f }- e^{-i\pi n f}
\right]
\right\}
            +{1+x\over 3}
$$
Observe that for $x\to 0$ ${\pi^2 \over 4 \omega_1^2}= {1 \over F(x)^2} 
= 1-{1\over 2} x +O(x^2)$. We recall that $e^{i\pi C \tau(x)} =h(x)^C (x/16)^C$, where $h(x) = 1+O(x)$.   We define
$$
{\cal X}(x) : = h(x)^{2-\nu_2} e^{-i\pi \nu_1}  \left( {x\over 16}
\right)^{2-\nu_2}
 e^{-i\pi {v(x) \over \omega_1(x)}}
$$
$$ 
   {\cal Y}(x):= h(x)^{\nu_2} e^{i \pi \nu_1} \left({ 
x\over 16}\right)^{\nu_2} e^{i \pi {v(x) \over \omega_1(x)}}
$$
Thus
$$
y(x) = {x\over 2} +O(x^2) +4\left(1-{x\over 2}+O(x^2)\right) 
\left\{
\sum_{n\geq 1}{ n {\cal X}^n(x) \over 1 - h(x)^{2n}(x/16)^{2n}} 
\left[
2 {\cal Y}^n(x)- {\cal Y}^{2n}(x) -1 
\right] 
-{ {\cal Y}(x) 
\over \left[ 1- {\cal Y}(x) \right]^2} 
\right\}
$$
Note that the denominator $[ 1- {\cal Y}(x) ]^2$ in the last term 
does not vanish on ${\cal D}$.  
 
The above formula immediately yields (\ref{critA})  and 
(\ref{critE}), (\ref{critEE}). It is enough 
 to  recall that in ${\cal D}(r;\nu_1,\nu_2)$ the dominant powers are 
$x^{\nu_2}$, $x^{1-\nu_2}$. We also need to observe that 
$$
   |x^{\nu_2+\alpha}|=|x|^{{\cal V}+\alpha}~{|x_0^{\nu_2+\alpha}|  \over 
|x_0|^{{\cal V}+\alpha}},~~~~\forall\alpha\in{\bf C}.
$$
and that $v(x)\to 0$ for $x\to 0$, according to (\ref{bundooo1}).  
Therefore,  (\ref{critA})  and 
(\ref{critE}), (\ref{critEE}) follow simply taking the leading terms in the 
expansion. 

As for (\ref{critB}), we observe that if ${\cal V}=1$,   $v(x)$ does not 
vanish, namely:  
 $$\phi(x):= 
\sum_{m\geq 1} 
b_{0m} \left[ e^{-i\pi \nu_1} \left({x\over 16}\right)^{1-\nu_2}\right]^m
\not \to 0 \hbox{ as } x\to 0
$$
and $e^{i\pi {v\over \omega_1}}= e^{2i\phi_1}(1+O(x))$. 
This implies that the dominant terms in the Fourier expansion are:
$$ 
y(x) = \left[{x\over 2} - 4 e^{i\pi \nu_1} \left({x\over 16}\right)^{\nu_2} 
e^{2i\phi(x)}-4 e^{-i\pi \nu_1} \left({x\over 16}\right)^{2-\nu_2} 
e^{-2i\phi(x)}
 \right] ~(1+O(x))
$$
from which (\ref{critB}) follows. Note that $x$, $x^{\nu_2}$ and 
$x^{2-\nu_2}$ are of the same order.

As for (\ref{critC}), $v(x)$ does not vanish when ${\cal V}=0$, because 
$x^{\nu_2}\not\to 0$ . Namely 
$$
\psi(x):= \sum_{m\geq 1} c_{0m}\left[ e^{i\pi \nu_1} \left(
{x\over 16}\right)^{\nu_2}\right]^m  \not \to 0  \hbox{ as } x\to 0
$$
In this case, it is convenient to keep the term $\sin^{-2}(\pi f/2)$. It is 
$$
  \sin^{-2}\left[-i{\nu_2\over 2} \ln{x\over 16} +{\pi \nu_1\over 2} 
+\psi(x)+O(x) \right]=
$$
$$
= \sin^{-2}\left[-i{\nu_2\over 2} \ln{x\over 16} +{\pi \nu_1\over 2} 
+\psi(x) \right](1+O(x))
$$
The last step is possible because 
$ -i{\nu_2\over 2} \ln{x\over 16} +{\pi \nu_1\over 2} 
+\psi(x) \neq 0$ on ${\cal D} (r;\nu_1,\nu_2)$. So, (\ref{critC}) follows.  Note   that 
$\sin\left[-i{\nu_2\over 2} \ln{x\over 16} +{\pi \nu_1\over 2} 
+\psi(x) \right]\neq 0$ on the domain, namely $y(x)$ in (\ref{critC}) does 
not have poles in ${\cal D} (r;\nu_1,\nu_2)$. 

\rightline{$\Box$} 

\vskip 0.3 cm 
\noindent
{\it Remark:} 
 In (\ref{critC}) the $\sin^{-2}(...)$ is 
never infinite on the domain. Actually, it is clear from the proof of 
Theorem 1 that the domains ${\cal D}_i(r_N;\nu_1,\nu_2+2N)$, $i=1,2$,  where exactly 
chosen in such a way that $\sin(...) \neq 0$. Therefore, the {\it movable 
poles} of the transcendent  must lie outside 
${\cal D}(r;\nu_1,\nu_2)$.

\vskip 0.3 cm


\section{Special Cases}\label{special}

 There are some special cases when the domain ${\cal D}(r;\nu_1,\nu_2)$ 
can be 
enlarged. 

\vskip 0.3 cm 
\noindent
{\bf FIRST CASE:  $\beta= 1-2\delta=0$.}
\vskip 0.2 cm 
 Theorem 1 can be stated as follows: 
\vskip 0.3 cm 
\noindent 
{\bf Theorem 1 -- special case $\beta= 1-2\delta=0$ :} { \it 
 For any 
complex  $\nu_1$, $\nu_2$ with the constraint 
$$
 0 <\nu_2<2~~~~ \hbox{ if $\nu_2$ is real}
$$ 
there exist a  positive number $r<1$  and a  transcendent 
$$
 y(x) = \wp \Bigl(\nu_1\omega_1(x)+\nu_2\omega_2(x) +v(x;\nu_1,\nu_2);~ \omega_1(x),\omega_2(x)
 \Bigr)+{1+x\over 3}
$$
such that $v(x;\nu_1,\nu_2)$ is holomorphic in the domain 
$$ 
{\cal D}(r;\nu_1,\nu_2) = \left\{ 
x\in \tilde{{\bf C}_0} ~|~ |x|<r,~\left| e^{-i\pi \nu_1} \left({x\over 16} 
\right)^{2-\nu_2}\right|<r,~\left| e^{i\pi \nu_1} \left({x\over 16} 
\right)^{\nu_2}\right|<r
\right\} ~\hbox{ if } \Im \nu_2 \neq 0
$$ 
or in the domain ${\cal D}_0(r)$ if $0<\nu_2<2$. It has convergent expansion:   
$$
 v(x;\nu_1,\nu_2) = 
\sum_{n\geq 1} a_n x^n + \sum_{n\geq 0,m\geq 1} b_{nm}
x^n\left[e^{-i\pi \nu_1}\left({x\over 16} \right)^{2-\nu_2} \right]^m
+
 \sum_{n\geq 0,m\geq 1} c_{nm}
x^n\left[e^{i\pi \nu_1}\left({x\over 16} \right)^{\nu_2} \right]^m
$$      
}

\vskip 0.3 cm 

Note that the domain is larger than the general case and that $\nu_2=1$ is 
now allowed.  

\vskip 0.3 cm 
\noindent
{\it Proof:} 
 If we go back to the proof of Theorem 1 we see that in
 the equation (\ref{fuchseq}) we only have the functions $ \wp(u/2+
\varepsilon_1
 \omega_1+2N_2 \omega_2)$. This means that we do not have $N_2^{\prime}$. 
 Let us put $N_2=N$. Then, we do the proof of Theorem 1 in the domain 
 (\ref{newdomain}) with $N_2=N$ and $\varepsilon_2=0$. A remarkable fact happens. While in proof for the general case we have to distinguish  (for given $N$, $\nu_1$, $\nu_2$) between  ${\cal D}_1(r_N; \nu_1,\nu_2+2N)$ and 
${\cal D}_2(r_N; \nu_1,\nu_2+2N)$,  
in the present special case we prove the existence of $v(x)$ on the larger domain 
\be
{\cal D}(r_N;\nu_1,\nu_2+2N):= \left\{ 
x\in \tilde{{\bf C}_0} ~|~ |x|<r_N,~\left| e^{-i\pi \nu_1} \left({x\over 16} 
\right)^{2-\nu_2-2N}\right|<r,~\left| e^{i\pi \nu_1} \left({x\over 16} 
\right)^{\nu_2+2N}\right|<r
\right\}
\label{31anni!}
\ee 
where 
  $v(x;\nu_1,\nu_2+2N)$ is represented as:   
$$ 
\sum_{n\geq 1} a_n x^n + \sum_{n\geq 0,m\geq 1} b_{nm}
x^n\left[e^{-i\pi \nu_1}\left({x\over 16} \right)^{2-\nu_2-2N} \right]^m
+
 \sum_{n\geq 0,m\geq 1} c_{nm}
x^n\left[e^{i\pi \nu_1}\left({x\over 16} \right)^{\nu_2+2N} \right]^m
$$
 We note that, if $\Im \nu_2=0$, the  constraint $
 0 <\nu_2<2$ holds. 
Namely, now $\nu_2=1$ is allowed. When $\nu_2$ is real, 
 $v(x)$ is equal to (\ref{31anni!}) with $N=0$ and 
the domain is  simply given by the condition $\max\{ |x|, \left| 
e^{-i\pi \nu_1} \left({x\over 16} 
\right)^{2-\nu_2}\right|,\left| e^{i\pi \nu_1} \left({x\over 16} 
\right)^{\nu_2}\right|\} <r$, namely $|x|<r$ for $r$ small enough. 
This is the domain ${\cal D}_0(r)$.

\rightline{$\Box$} 
\vskip 0.3 cm

The critical  behavior for $x\to 0$  can be studied along 
$$
   \arg x = \arg x_0 +{\Re \nu_2  - {\cal V} \over \Im \nu_2} \ln {|x| 
\over |x_0|} 
$$
contained in ${\cal D}(r;\nu_1,\nu_2)$ for 
$x_0\in {\cal D}(r;\nu_1,\nu_2)$ and $
0\leq {\cal V}\leq 2
$. 
If $\Im \nu_2=0$ any path going to $x=0$ is contained in 
${\cal D}_0(r)$.  
The critical behavior is obtained through the Fourier expansion. We give the result.

If $\Im \nu_2\neq 0$, the transcendent $y(x)=\wp\bigl(
\nu_1\omega_1(x)+\nu_2\omega_2(x)+v(x;\nu_1,\nu_2)\bigr)+(1+x)/3 $ 
has the following behaviors on ${\cal D}(r;\nu_1,\nu_2)$:

\vskip 0.2 cm 
\noindent 
For $0<{\cal V}<1$:
\be
y(x)= -{1\over 4} \left[{e^{i\pi \nu_1} \over 16^{\nu_2-1}}\right] ~
          x^{\nu_2} ~\left(1+ O(|x^{\nu_2}|)\right).
\label{critAs1}
\ee

\vskip 0.2 cm 
\noindent 
For $1<{\cal V}<2$:
\be
y(x)= -{1\over 4} \left[{e^{i\pi \nu_1} \over 16^{\nu_2-1}}\right]^{-1} ~
          x^{2-\nu_2} ~\left(1+ O(|x^{2-\nu_2}|)\right). 
\label{critAAs1}
\ee

\vskip 0.2 cm 
\noindent
For ${\cal V}=1$:
\be
y(x) =
x~\sin^2
                                                             \left(
i{1-\nu_2\over 2}\ln{|x|\over 16}+{\pi \nu_1 \over 2}
                                                              \right)(1+O(x)).
\label{critBs1}
\ee
\vskip 0.2 cm 
\noindent
For ${\cal V}=0$:
\be
y(x)=  \left[
       {x\over 2} +\sin^{-2}
       \left( 
  -i{\nu_2\over 2} \ln {x\over 16} +{\pi \nu_1 \over 2} + 
  \sum_{m\geq 1} c_{0m} 
\left[
e^{i\pi \nu_1} 
\left(
                     {x\over 16}
\right)^{\nu_2}
\right]^m
\right)           
       \right] \bigl(1+O(x)\bigr).
\label{critCs1}
\ee
\vskip 0.2 cm 
\noindent
For ${\cal V}=2$:
\be
y(x)=  \left[
       {x\over 2} +\sin^{-2}
       \left( 
  i{2-\nu_2\over 2} \ln {x\over 16} +{\pi \nu_1 \over 2} + 
  \sum_{m\geq 1} b_{0m} 
\left[
e^{-i\pi \nu_1} 
\left(
                     {x\over 16}
\right)^{2-\nu_2}
\right]^m
\right)           
       \right] \bigl(1+O(x)\bigr).
\label{critDs1}
\ee

\vskip 0.2 cm 
\noindent
For $\nu_2$  real, satisfying the constraints $0<\nu_2<2$, the transcendent $y(x)=\wp\bigl(
\nu_1\omega_1(x)+\nu_2\omega_2(x)+v(x;\nu_1,\nu_2)\bigr)+(1+x)/3 $ 
has the following behaviors on ${\cal D}_0(r)$: 
\be
 y(x)= -{1\over 4} \left[{e^{i\pi \nu_1} \over 16^{\nu_2-1}}\right] ~
          x^{\nu_2} ~\left(1+ O(|x^{\nu_2}|)\right),
~~~~0<\nu_2<1
\label{critEs1}
\ee
\be
y(x)= -{1\over 4} \left[{e^{i\pi \nu_1} \over 16^{\nu_2-1}}\right]^{-1} ~
          x^{2-\nu_2} ~\left(1+ O(|x^{2-\nu_2}|)\right),
~~~~1<\nu_2<2
\label{critEEs1}
\ee
\be
y(x) =
x~\sin^2
                                                             \left(
{\pi \nu_1 \over 2}
                                                              \right) \bigl(1+O(x)\bigr),
~~~~\nu_2=1, ~~~~(\nu_1\neq 0)
\label{critEEEs1}
\ee

\vskip 0.5 cm 
\noindent
{\bf SECOND CASE:  $\alpha= \gamma=0$.}
\vskip 0.2 cm 
Theorem 1 is restated as follows:

\vskip 0.3 cm 
\noindent 
{\bf Theorem 1 -- special case $\alpha= \gamma=0$:} { \it 
 For any 
complex  $\nu_1$, $\nu_2$ with the constraint 
$$
 -1<\nu_2<1 ~~~~\hbox{ if $\nu_2$ is real}
$$ 
there exist a  positive number $r<1$  and a  transcendent 
$$
 y(x) = \wp \Bigl(\nu_1\omega_1(x)+\nu_2\omega_2(x) +v(x;\nu_1,\nu_2);~ \omega_1(x),\omega_2(x)
 \Bigr)+{1+x\over 3}
$$
such that $v(x;\nu_1,\nu_2)$ is holomorphic in the domain 
$$
{\cal D}(r;\nu_1,\nu_2)= 
\left\{
x\in \tilde{{\bf C}_0} ~|~ |x|<r,~\left| e^{-i\pi \nu_1} \left({x\over 16} 
\right)^{1-\nu_2}\right|<r,~\left| e^{i\pi \nu_1} \left({x\over 16} 
\right)^{\nu_2+1}\right|<r
\right\}
$$ 
or in the domain ${\cal D}_0(r)$ if $-1<\nu_2<1$, where  
it has convergent expansion : 
  $$
 v(x;\nu_1,\nu_2) 
= \sum_{n\geq 1} a_n x^n + \sum_{n\geq 0,m\geq 1} b_{nm}
x^n\left[e^{-i\pi \nu_1}\left({x\over 16} \right)^{1-\nu_2} \right]^m
+
 \sum_{n\geq 0,m\geq 1} c_{nm}
x^n\left[e^{i\pi \nu_1}\left({x\over 16} \right)^{\nu_2+1} \right]^m
$$     
}

\vskip 0.3 cm 
Note that now $\nu_2=0$ is allowed. 

\vskip 0.3 cm
\noindent
{\it Proof: }
If we go back to the proof of Theorem 1, we see that in
 the equation (\ref{fuchseq}) we only have the functions $ \wp(u/2+
\varepsilon_1
 \omega_1+(2N^{\prime}_2+1) \omega_2)$. This means that we do not have 
$N_2$. This time, let us put $N_2^{\prime}=N$. 
 Therefore we do the proof of Theorem 1 for the domain 
 (\ref{newdomain}) with $N^{\prime}_2=N$ and $\varepsilon_2=1$. Therefore, 
 on the new (bigger) domain
\be
{\cal D}(r_N;\nu_1,\nu_2+2N)= 
\left\{
x\in \tilde{{\bf C}_0} ~|~ |x|<r_N,~\left| e^{-i\pi \nu_1} \left({x\over 16} 
\right)^{1-\nu_2-2N}\right|<r,~\left| e^{i\pi \nu_1} \left({x\over 16} 
\right)^{\nu_2+1+2N}\right|<r
\right\}
\label{31anni!1}
\ee
there exist   one 
 $v(x)$,  holomorphic with convergent expansion:  
$$
 v(x;\nu_1,\nu_2+2N) 
= \sum_{n\geq 1} a_n x^n + \sum_{n\geq 0,m\geq 1} b_{nm}
x^n\left[e^{-i\pi \nu_1}\left({x\over 16} \right)^{1-\nu_2-2N} \right]^m
+
$$
$$
+
 \sum_{n\geq 0,m\geq 1} c_{nm}
x^n\left[e^{i\pi \nu_1}\left({x\over 16} \right)^{\nu_2+1+2N} \right]^m
$$
In a sense, in  the case $\beta=1-2\delta=0$ the disjoint union 
 ${\cal D}_1(r_N;\nu_1,\nu_2+2N) \cup {\cal D}_2(r_{N};\nu_1,\nu_2+2N)$ 
was replaced
  by the bigger domain (\ref{31anni!}). In the present case  the 
disjoint union  
${\cal D}_1(r_N;\nu_1,\nu_2+2N) \cup {\cal D}_2(r_{N+1};\nu_1,\nu_2+2[N+1])$ 
is replaced by the bigger domain (\ref{31anni!1}).

 If $\nu_2$ is real, this time $\nu_2=0$ is allowed, but not $\nu_2=1$. 
 Then, it is convenient to choose the 
constraint 
$$
 -1<\nu_2<1
$$
by doing  $\nu_2\mapsto \nu_2-2$ if $1<\nu_2<2$. The domain is simply specified by $\max \left\{ |x|,~ \left| 
e^{-i\pi \nu_1} \left({x\over 16} 
\right)^{1-\nu_2}\right|,\right. $ $\left.\left| e^{i\pi \nu_1} \left({x\over 16} 
\right)^{\nu_2+1}\right| \right\} <r$, namely $|x|<r$ for $r$ small enough. 
This is ${\cal D}_0(r)$.    

\rightline{$\Box$}

\vskip 0.3 cm 
 We  compute  the critical behaviors for $x\to 0$  along 
$$ 
\arg x = \arg x_0 + {\Re \nu_2 -{\cal V} \over \Im \nu_2} 
\ln {|x|\over |x_0|}, 
$$
which is contained in 
${\cal D}(r;\nu_1,\nu_2)$ if and only if 
$x_0\in {\cal D}(r;\nu_1,\nu_2)$ and $
-1 \leq {\cal V}\leq 1
$. 
If $\Im \nu_2=0$ any path converging  to $x=0$ is contained in ${\cal D}_0 
(r)$.   
We use Fourier expansion;  the result is the following. 
\vskip 0.2 cm 
\noindent 
For $0<{\cal V}<1$:
\be
y(x)= -{1\over 4} \left[{e^{i\pi \nu_1} \over 16^{\nu_2-1}}\right] ~
          x^{\nu_2} ~\left(1+ O(|x^{1-\nu_2}|)\right).
\label{critAs2}
\ee

\vskip 0.2 cm 
\noindent 
For $-1<{\cal V}<0$:
\be
y(x)= -{1\over 4} \left[{e^{i\pi \nu_1} \over 16^{\nu_2+1}}\right]^{-1} ~
          x^{-\nu_2} ~\left(1+ O(|x^{\nu_2+1}|)\right).
\label{critAAs2}
\ee

\vskip 0.2 cm 
\noindent
For ${\cal V}=1$:
\be
y(x) =
x~\sin^2
                                                             \left(
i{1-\nu_2\over 2}\ln{|x|\over 16}+{\pi \nu_1 \over 2}+ 
\sum_{m\geq 1} b_{0m}
                                 \left[
e^{-i\pi \nu_1}
                                 \left(
{x\over 16}
                                        \right)^{1-\nu_2}
                                       \right]^m
                                                              \right) \bigl(1+O(x)\bigr). 
\label{critBs2}
\ee
\vskip 0.2 cm 
\noindent
For ${\cal V}=-1$:
\be
y(x)=x~\sin^2\left(
-i{\nu_2+1\over 2}\ln{|x|\over 16} +{\pi \nu_1 \over 2}+ 
\sum_{m\geq 1} c_{0m}\left[e^{i\pi \nu_1}\left({x\over 16}
\right)^{\nu_2+1}
\right]^m
\right) \bigl(1+O(x)\bigr).
\label{critBBs2}
\ee
\vskip 0.2 cm 
\noindent
For ${\cal V}=0$:
\be
y(x)=  \left[
       {x\over 2} +\sin^{-2}
       \left( 
  -i{\nu_2\over 2} \ln {x\over 16} +{\pi \nu_1 \over 2}
\right)           
       \right] \bigl(1+O(x)\bigr).
\label{critCs2}
\ee

\vskip 0.2 cm 
\noindent
For $\nu_2$  real, the  behaviors on ${\cal D}_0(r)$ are:
\be
 y(x)= -{1\over 4} \left[{e^{i\pi \nu_1} \over 16^{\nu_2-1}}\right] ~
          x^{\nu_2} ~\left(1+ O(|x^{1-\nu_2}|)\right),
~~~ \hbox{ if }0<\nu_2<1.
\label{critEs2}
\ee
\be
y(x)= -{1\over 4} \left[{e^{i\pi \nu_1} \over 16^{\nu_2+1}}\right]^{-1} ~
          x^{-\nu_2} ~\left(1+ O(|x^{\nu_2+1}|)\right),
~~~\hbox{ if } -1<\nu_2<0.
\label{critEEs2}
\ee
\be 
y(x) = \left[{1\over \sin^2\left({\pi \nu_1\over 2}\right) } + {x\over 2}
\right] \bigl(1+O(x)\bigr),~~~\hbox{ if }\nu_2=0,~~(\nu_1\neq 0).
\label{critEEEs2}
\ee


\section{Points $x=1,\infty$ -- Comments to  Theorem 3 and Critical behavior}\label{theorem3}

Two independent solutions of the hyper-geometric equation ${\cal L}(u)=0$ are 
the hyper-geometric function 
$F(x;1/2,1/2,1)$ and 
$$
g\left({1\over 2},{1\over 2},1;x\right):= F_1(x) + 
F\left({1\over 2},{1\over 2},1;x\right)\ln x
$$ 
The 
following connection formulae hold \cite{Norlund}:

\vskip 0.2 cm 
\noindent
i) Connection 0 -- 1  
$$
    F\left({1\over 2}, {1\over 2},1;x\right)= -{1\over \pi} 
g\left({1\over 2}, {1\over 2},1;1-x\right),~~~~|\arg(1-x)|<\pi
$$
$$
    g\left({1\over 2}, {1\over 2},1;x\right)= -{ \pi} 
F\left({1\over 2}, {1\over 2},1;1-x\right),~~~~|\arg x|<\pi
$$

\vskip 0.2 cm 
\noindent
ii) Connection 0 -- $\infty$ 
$$
     F\left({1\over 2}, {1\over 2},1;x\right)=
{x^{-{1\over 2}}\over \pi}
\left[ 
i g\left({1\over 2}, {1\over 2},1;{1\over x}\right)+\pi F\left({1\over 2}, 
{1\over 2},1;{1 \over x}\right) 
\right],~~~~-2\pi <\arg x <0
$$
$$
  g\left({1\over 2}, {1\over 2},1;x\right)= x^{-{1\over 2}}  
g\left({1\over 2}, {1\over 2},1;{1\over x}\right),~~~~|\arg x|<\pi
$$

 We recall that $\omega_1(x)$ and $\omega_2(x)$ are $\pi/2~F(1/2,1/2,1;x)$ 
and $-i/2~g(1/2,1/2,1;x)$ respectively. 
This  representation is convenient in a neighborhood 
of $x=0$. In this section we will use the notation 
$\omega^{(0)}_1$ and $\omega^{(0)}_2$ instead of $\omega_1$ and $\omega_2$.  
Namely 
$$ 
\omega^{(0)}_1(x)= {\pi \over 2} F\left({1\over 2},{1\over 2},1;x\right),
~~~~\omega^{(0)}_2(x)= -{i\over 2} g\left({1\over 2},{1\over 2},1;x\right)
$$
Taking into account the above connection formulae, we define in a neighborhood
 of $x=1$:
$$
\omega^{(1)}_1(x):=\omega_2^{(0)}(x)=  
{i\pi \over 2} F\left({1\over 2},{1\over 2},1;1-x\right), ~~~|\arg x |<\pi 
$$
$$
\omega^{(1)}_2(x):=\omega^{(0)}_1(x) = 
-{1\over 2} g\left({1\over 2},{1\over 2},1
;1-x\right),~~~|\arg(1-x)|<\pi 
$$  We also define, in a neighborhood of $x=\infty$: 
$$
\omega^{(\infty)}_1(x) := \omega^{(0)}_1(x) +  \omega^{(0)}_2(x) =
{\pi \over 2} x^{-{1\over 2}} F\left({1\over 2},{1\over 2},1;{1\over x}\right), ~~~-\pi<\arg x<0
$$
$$
\omega^{(\infty)}_2(x):=  \omega^{(0)}_2(x)= -{i\over 2} x^{-{1\over 2}}
 g\left({1\over 2},{1\over 2},1;{1\over x}\right),~~~ |\arg x  |<\pi
$$
 The above functions have  branch cuts specified by 
the constrains on $\arg x$. Once they are so defined, they 
 are continued on the 
universal covering of ${\bf P}^1\backslash \{ 0,1,\infty\}$.

In the following, we use superscripts 0, 1, $\infty$ for solutions $u(x)$ of 
(\ref{fuchseq}) represented around $x=0,1,\infty$ respectively. 
At $x=0$ we looked for solutions 
$${u^{(0)}\over 2}= \nu^{(0)}_1 \omega^{(0)}_1 +\nu^{(0)}_2\omega^{(0)}_2 +
v^{(0)}$$ 
and we found the representations of $v^{(0)}(x)$ in Theorem 1. 
 Now we look for a solution in a neighborhood of $x=1$ of the type
$$
{u^{(1)}\over 2}= \nu^{(1)}_1 \omega^{(1)}_1 +\nu^{(1)}_2\omega^{(1)}_2 +
v^{(1)}
$$ 
 and for a solution in a neighborhood of $x=\infty$ of the type
$$
{u^{(\infty)}\over 2}= \nu^{(\infty)}_1 \omega^{(\infty)}_1 +\nu^{(\infty)}_2
\omega^{(\infty)}_2 +
v^{(\infty)}
$$ 
  They  yield Painlev\'e transcendents 
$y(x)= \wp (u^{(1)}/2; \omega_1^{(0)},
\omega^{(0)}_2)+(1+x)/3\equiv  \wp (u^{(1)}/2; \omega_1^{(1)},
\omega^{(1)}_2)+(1+x)/3$ and   $y(x)= \wp (u^{(\infty)}/2; \omega_1^{(0)},
\omega^{(0)}_2)+(1+x)/3\equiv  \wp (u^{(\infty)}/2; \omega_1^{(\infty)},
\omega^{(\infty)}_2)+(1+x)/3$.

The definitions of $\omega^{(i)}_j$ are convenient because 
 the proof of Theorem 1 can be repeated with no changes 
in a neighborhood of 
$x=1$ and $\infty$.  This yields the statement of Theorem 3 in the Introduction.


\subsection{Critical Behaviors at $x=1,\infty$}\label{HIEI-ZAN}

 From Theorem 3 and the Fourier 
expansion of the Weierstrass function we obtain the critical behaviors of the 
transcendents of Theorem 3.  

\vskip 0.3 cm 
\noindent 
{\bf In a neighborhood of $x=1$:}  
\vskip 0.2 cm 
 If $\Im \nu^{(1)}_2=0$, we let $x\to 1$ along any regular path. Otherwise, we consider the 
paths  
$$ 
 \arg(1-x)=\arg(1-x_0) + {\Re \nu_2^{(1)} -{\cal V}
\over 
 \Im  \nu_2^{(1)}} \ln{ |1-x|\over |1-x_0|},~~~0\leq {\cal V}\leq 1
$$
joining $x_0$ to $x=1$. They lie  in 
${\cal D}(r;\nu_1^{(1)},\nu_2^{(1)})$ if 
$x_0 \in {\cal D}(r;\nu_1^{(1)},\nu_2^{(1)})$.

\vskip 0.2 cm
\noindent
Let $\nu_1^{(1)},\nu_2^{(1)}$ be given. If $\Im \nu_2^{(1)}\neq 0$, the transcendent defined in 
Theorem 3 in 
${\cal D}(r;\nu_1^{(1)},\nu_2^{(1)})$ has the following behaviors for $x\to 1$:

\vskip 0.2 cm
\noindent
For $0<{\cal V}<1$: 
\be 
 y(x) = 1+{1\over 4} \left[{e^{-i\pi \nu_1^{(1)}}\over 16^{\nu_2^{(1)}-1}}
\right](1-x)^{\nu^{(1)}_2}  ~\left( 1+O(|1-x|^{\nu^{(1)}_2} + 
|1-x|^{1-\nu_2^{(1)} })\right).
\label{fat1}
\ee

\vskip 0.2 cm
\noindent
For ${\cal V}=0$:
$$
y(x)
=
1-\left\{
{1-x\over 2}+\right.
$$
$$\left.
 +\sin^{-2}\left( 
i{\nu_2^{(1)}\over 2} \ln {1-x\over 16} +{\pi \nu_1^{(1)} \over 2} 
-i 
\sum_{m\geq 1} c_{0m}\left[ 
e^{-i\pi \nu_1^{(1)}} \left({1-x\over 16}\right)^{\nu_2^{(1)}}
\right]^m
\right)
 \right\}~\bigl(
1+O(1-x)\bigr).
$$

\vskip 0.2 cm
\noindent
For ${\cal V}=1$: 
$$
y(x) 
=
1-(1-x)~\sin^2\left( 
i{\nu_2^{(1)}-1\over 2}\ln{1-x\over 16}+{\pi\nu_1^{(1)}\over 2} 
\right.-
$$
$$ \left.
-i
\sum_{m\geq 1} 
b_{0m}\left[
e^{i\pi \nu_1^{(1)}} \left({1-x\over 16}\right)^{1-\nu_2^{(1)}}
\right]^m
\right)~\bigl(
1+O(1-x)
\bigr).
$$

\vskip 0.2 cm
\noindent
If $\Im \nu_2=0$, the transcendent defined on ${\cal D}_0(r)$ 
has behavior  (\ref{fat1})   
 when $0<\nu_2^{(1)}<1$, or behavior 
\be 
 y(x) = 1+{1\over 4} \left[{e^{-i\pi \nu_1^{(1)}}\over 16^{\nu_2^{(1)}-1}}
\right]^{-1}(1-x)^{2-\nu^{(1)}_2}  ~\left( 1+O(|1-x|^{2-\nu^{(1)}_2} + 
|1-x|^{\nu_2^{(1)}-1 }|)\right)
\label{fat2}
\ee
when  $1<\nu_2^{(1)}<2$.

\vskip 0.3 cm 
\noindent 
{\bf In a neighborhood of $x=\infty$:}  
\vskip 0.2 cm 
 If $\Im \nu^{(\infty)}_2=0$, we let $x\to \infty$ along any path. 
Otherwise, we consider 
the paths  
$$
\arg x = \arg x_0 + {\Re\nu_2^{(\infty)} -{\cal V}\over \Im \nu_2^{(\infty)}
} \ln {|x|\over |x_0|},~~~0\leq {\cal V}\leq 1
$$
joining $x_0\in  D(r; \nu_1^{(\infty)},\nu_2^{(\infty)}) $ to infinity. 
 \vskip 0.2 cm 
\noindent 
Let $\nu_1^{(\infty)},\nu_2^{(\infty)}$ be given. If $\Im \nu_2^{(\infty)}\neq 0$ the transcendent of Theorem 3 defined in 
$ D(r; \nu_1^{(\infty)},\nu_2^{(\infty)}) $ has the following behaviors for $x\to \infty$: 

\vskip 0.2 cm
\noindent
  For $0<{\cal V}<1$: 
\be
y(x) 
=
-{1\over 4} \left[{
e^{i\pi \nu_1^{(\infty)}}\over 16^{\nu_2^{(\infty)}-1}
}\right]
x^{1-\nu_2^{(\infty)}} \left(
1+O(|x|^{-\nu_2^{(\infty)} }+|x|^{\nu_2^{(\infty)} -1})
\right).
\label{fatt1}
\ee

\vskip 0.2 cm
\noindent
For ${\cal V}=0$:
$$
y(x) 
=\left(1+O(x^{-1})\right)
 \left\{{1\over 2} +x\sin^{-2} \left(
-i{\nu_2^{(\infty)} \over 2} \ln {16\over x} +{\pi \nu_1^{(\infty)}\over 2} 
+\right. \right.$$
$$ \left. \left.
\sum_{m\geq 1} c_{0m} \left[
e^{i\pi \nu_1^{(\infty)}} \left({16\over x}\right)^{\nu_2^{(\infty)}}
\right]^m
\right)
\right\}.
$$

\vskip 0.2 cm
\noindent
For ${\cal V}=1$:
$$
y(x) =
\left( 1+O(x^{-1})\right)
\sin^2 \left(
i{1-\nu_2^{(\infty)} \over 2} \ln {16\over x} +{\pi \nu_1^{(\infty)}\over 2} 
+
\sum_{m\geq 1} b_{0m} \left[
e^{-i\pi \nu_1^{(\infty)}} \left({16\over x}\right)^{1-\nu_2^{(\infty)}}
\right]^m
\right).
$$

\vskip 0.2 cm
\noindent
If $\Im \nu_2=0$, the transcendent defined on ${\cal D}_0(r)$ has behavior 
 (\ref{fatt1}) when $0<\nu_2^{(\infty)}<1$, or 
it has behavior
\be
y(x) 
=
-{1\over 4} \left[{
e^{i\pi \nu_1^{(\infty)}}\over 16^{\nu_2^{(\infty)}-1}
}\right]^{-1}
x^{\nu_2^{(\infty)}-1} \left(
1+O(|x|^{\nu_2^{(\infty)} -2 }+|x|^{1-\nu_2^{(\infty)} })
\right)
\label{fatt2}
\ee
 when $1<\nu_2^{(\infty)}<2$.


\section{Some Considerations on Analytic Continuation} 

 We can easily study the effect of a small loop around a critical point on 
 the transcendents of Theorem 1 and 3. 

Consider a transcendent of Theorem 1
$$
  y(x)= \wp\Bigl( \nu_1 \omega^{(0)}_1 +\nu_2 \omega^{(0)}_2 +
v(x;\nu_1,\nu_2);~\omega^{(0)}_1,\omega^{(0)}_2 \Bigr) 
+{1+x\over 3}
$$
defined on ${\cal D}(r;\nu_1,\nu_2)$. 
If $\nu_2 $ is real, it is defined on ${\cal D}_0(r)$.

We do the loop $x\mapsto x e^{2\pi i }$,  where $x\in 
{\cal D}(r;\nu_1,\nu_2)$ (or ${\cal D}_0(r)$),  therefore $|x|<1$.  From the 
monodromy  properties of 
$F(1/2,1/2,1;x)$ and $g(1/2,1/2,1;x)$ we have
$$
   \omega_1^{(0)}(x) \mapsto  \omega_1^{(0)}(xe^{2\pi i }) = \omega_1^{(0)}(x),
$$
$$  \omega_2^{(0)}(x) 
\mapsto   \omega_2^{(0)}(xe^{2\pi i }) =  \omega_2^{(0)}(x) +2\omega_1^{(0)}(x)
$$
We also have 
$$
  v(x;\nu_1,\nu_2) \mapsto v(xe^{2\pi i};\nu_1,\nu_2)
\equiv v(x; \nu_1+2\nu_2,
\nu_2)
$$
The last step of the above equalities follows from the explicit 
expansion of $v(xe^{2\pi i};\nu_1,\nu_2)$  given in Theorem 1, {\it in the 
hypothesis} that also $x e^{2\pi i} \in {\cal D}(r;\nu_1,\nu_2)$ -- or, 
which is the same thing, $x \in {\cal D}(r;\nu_1+2\nu_2,\nu_2)$ 
[this always happens if $\nu_2$ is real].  
 Therefore 
$$
y(x) \mapsto y(xe^{2\pi i})=:y^{\prime}(x)=
$$
$$=
\wp\Bigl( \nu_1\omega_1^{(0)}(x)+\nu_2 [\omega_2^{(0)}(x)+2\omega^{(0)}_1(x)]+
 v(x;\nu_1+2\nu_2,
\nu_2);~\omega^{(0)}_1(x), \omega_2^{(0)}(x) +2\omega_1^{(0)}(x)
\Bigr)
+{1+x\over 3}
$$
$$
\equiv \wp\Bigl( (\nu_1+2\nu_2)\omega_1^{(0)}(x)+\nu_2 \omega_2^{(0)}(x)+ 
v(x;\nu_1+2\nu_2,
\nu_2);~\omega^{(0)}_1(x),\omega^{(0)}_2(x) \Bigr)
+{1+x\over 3}
$$
So, the effect of the loop is simply the transformation 
$$
  (\nu_1,\nu_2)\mapsto (\nu_1+2\nu_2,\nu_2)
$$ 
From the above, we also obtain the critical  behavior according to Theorem 2.

We remark that the above considerations hold only if both $x$ and $x e^{2\pi i}$ belongs to ${\cal D}(\nu_1,\nu_2)$ -- namely $x \in {\cal D}(\nu_1+2\nu_2,\nu_2)$.  
Otherwise, we can not represent $v(x e^{2\pi i};\nu_1,\nu_2)$ 
as $v(x; \nu_1+2\nu_2,\nu_2)$.  
 This last case may actually occur when $\Im \nu_2\neq 0$. If  the 
 point $x e^{2\pi i}$ lies outside the domain  it  may 
 be a  pole of $y(x)$.

\vskip 0.3 cm 
 The same procedure is applied to $$
  y(x) =\wp\Bigl( \nu^{(1)}_1 \omega^{(1)}_1 +\nu^{(1)}_2 \omega^{(1)}_2 +
v(x; \nu^{(1)}_1 ,\nu^{(1)}_2);~\omega^{(1)}_1,\omega^{(1)}_2
 \Bigr) +{1+x\over 3}
$$ 
 Now we denote by 
 ${\cal D}
(r;\nu^{(1)}_1,\nu^{(1)}_2)$ 
 the domain around $x=1$.   
The effect of the loop $(1-x) \mapsto (1-x) e^{2\pi i}$, $x\in {\cal D}
(r;\nu^{(1)}_1,\nu^{(1)}_2)$ is 
$$ 
(\nu^{(1)}_1,\nu^{(1)}_2) \mapsto (\nu^{(1)}_1-2\nu^{(1)}_2,\nu^{(1)}_2)
$$
if also $x \in {\cal D}(r;\nu^{(1)}_1-2\nu^{(2)}_1,\nu^{(1)}_2)$.

\vskip 0.3 cm 
 The  procedure is applied to $$
  y(x)= \wp\Bigl( \nu^{(\infty)}_1 \omega^{(\infty)}_1 +\nu^{(\infty)}_2
 \omega^{(\infty)}_2 +
v(x; \nu^{(\infty)}_1,\nu^{(\infty)}_2);~\omega^{(\infty)}_1,\omega^{(\infty)}_2\Bigr) +{1+x\over 3}
$$
 Now  ${\cal D}
(r,\nu^{(\infty)}_1,\nu^{(\infty)}_2)$ is the domain around $x=\infty$. 
The effect of the loop $ x \mapsto x e^{-2\pi i}$, $x\in {\cal D}
(r;\nu^{(\infty)}_1,\nu^{(\infty)}_2)$ is 
$$ 
(\nu^{(\infty)}_1,\nu^{(\infty)}_2) \mapsto (\nu^{(\infty)}_1+2
\nu^{(\infty)}_2,\nu^{(\infty)}_2)
$$
provided that  also $x \in {\cal D}(r;\nu^{(\infty)}_1+2\nu^{(\infty)}_1,
\nu^{(\infty)}_2)$.


\section{Connection Problem}\label{twoquestions}

\vskip 0.3 cm 
The analysis so far developed is only local.  Two natural questions arise. 

\vskip 0.2 cm
\noindent 
{\bf Question 1).} This question makes sense in $\nu_2$ is not real. 
 Suppose that $y(x)$ has the representation of Theorem 1 close to $x=0$: 
 $y(x)=\wp\bigl(\nu_1\omega_1+\nu_2\omega_2+v(x;\nu_1,\nu_2)\bigr)
+(1+x)/3$   
in  ${\cal D}(r;\nu_1,\nu_2)$.  
We ask if $y(x)$ may  have a 
representation  $y(x)=\wp\bigl(\nu_1^{\prime}\omega_1+\nu_2\omega_2+v(x;\nu_1,\nu_2+2N)\bigr)
+(1+x)/3$   in   ${\cal D}(r;\nu^{\prime}_1,
\nu_2+2N)$ for some integer  $N$. If this happens, we should be able to express $\nu^{\prime}_1$ as a function of $\nu_1,\nu_2,N$.

\vskip 0.2 cm 
\noindent
{\bf Question 2).} [Connection Problem].  
 May a given  transcendent have  three 
representations of Theorem 1 and Theorem 3 at $x=0,1,\infty$ at the same time?  If this happens, which is the  relation between the three sets 
$(\nu_1^{(0)},\nu_2^{(0)})$, $(\nu_1^{(1)},\nu_2^{(1)})$, $(\nu_1^{(\infty)},
\nu_2^{(\infty)})$?

\vskip 0.2 cm

The answer to this questions is positive, but it can not be achieved with  the 
method so far used, which only allows to obtain local results. 
For example, suppose we want to solve question 2) for 
 $y(x)$  represented as in Theorem 1. We can write
$$
y(x)= \wp(\nu^{(0)}_1\omega^{(0)}_1(x)+\nu_2^{(0)}\omega_2^{(0)}(x)+v^{(0)}(x))+
{1+x\over 3}
$$
\be
\equiv 
\wp(\nu^{(0)}_2\omega^{(1)}_1(x)+\nu_1^{(0)}\omega_2^{(1)}(x)+v^{(0)}(x))+
{1+x\over 3}
\label{impossibile}
\ee
  There also exists a  transcendent
\be
y(x)= 
\wp(\nu_1^{(1)}\omega^{(1)}_1(x)+\nu_2^{(1)}\omega_2^{(1)}(x)+v^{(1)}(x))+
{1+x\over 3}
\label{impossibileno}
\ee
where $v^{(1)}(x)$ is represented by Theorem 3. In particular, $v^{(1)}(x)$ 
  is bounded as $x\to 1$.   
On the other hand, the function $v^{(0)}$ is only 
known for small $x$, and it may diverge as $x\to 1$. In the best hypothesis, 
we may suppose that after the  rescaling 
$$ 
v^{(0)}(x)\mapsto v^{(0)}(x)-\delta\nu_1~ \omega_1^{(1)}(x) 
-\delta\nu_2~ \omega_2^{(1)}(x),~~~~~\delta\nu_1,\delta\nu_2\in{\bf C}
$$
the new function 
$v^{(0)}-\delta\nu_1 \omega_1^{(1)} 
-\delta\nu_2 \omega_2^{(1)}$ vanishes as $x\to 1$ (this may happen, 
since the divergence of $\omega_2^{(1)}$, which is  $\ln (1-x)$ as $x\to 1$, may cancel the divergence of $v^{(0)}(x)$ as $x\to 1$).  
Therefore, (\ref{impossibile}) is a transcendent (\ref{impossibileno}) where 
$$ 
  \nu_1^{(1)}= \nu_2^{(0)}+\delta\nu_1,
~~~~\nu_2^{(1)}= \nu_1^{(0)}+\delta\nu_2
$$
and 
$$
 v^{(1)}(x) = v^{(0)}-\delta\nu_1 ~\omega_1^{(1)} 
-\delta\nu_2~ \omega_2^{(1)}
$$

Unfortunately,  we are not able to say if the above rescaling is possible using the local 
analysis.  
And if it is, we do not know $\delta\nu_1$, $\delta\nu_2$.

\vskip 0.3 cm 

 The answer to questions 1) and 2)  can  be obtained by the method of isomonodromic deformations (which, on the other hand, has some limitations in providing the local behavior).

 \subsection{Picard Solutions}
 
 Before  answering Questions 1 and  2 in general,  
we recall that there is 
one case when the local analysis in elliptic 
representation becomes global, namely when $v(x)=0$. This happens when 
$\alpha=\beta=\gamma=1-2\delta=0$. 
This case was already well known to Picard \cite{Picard} and 
it was  studied in \cite{M}. 
The equation (\ref{fuchseq}) reduces to the 
hyper-geometric equation  $
  {\cal L}(u)=0$, and it has general solution $
  u(x)/ 2 = \nu_1 \omega^{(0)}_1(x) +\nu_2\omega^{(0)}_2(x)
$. Therefore,  
$$
  y(x) =\wp\left( \nu_1 \omega^{(0)}_1(x) +\nu_2\omega^{(0)}_2(x);\omega^{(0}_1,\omega^{(0)}_2 \right)+{1+x\over 3}
$$

\vskip 0.2 cm 
\noindent 
 Question 1)
 can be answered immediately. There is no need to prove Theorem 1 
because we have no $v(x)$! The function  $
  u(x)/ 2 = \nu_1 \omega^{(0)}_1(x) +\nu_2\omega^{(0)}_2(x)
$ in the argument of $\wp$ is defined for any $x\neq 0,1,\infty$. 

We  obtain the critical behavior by Fourier-expanding 
$$ 
\wp\left( \nu_1 \omega^{(0)}_1(x) +\nu_2\omega^{(0)}_2(x);\omega^{(0}_1,\omega^{(0)}_2 \right)=\wp\left( \nu_1 \omega^{(0)}_1(x) +[\nu_2+2N]\omega^{(0)}_2(x);\omega^{(0}_1,\omega^{(0)}_2 \right)
$$
 For any fixed $N$, the expansion is performed if $|x|<1$ and  
$$ 
\left|\Im \left[{\nu_1\over 2} +\left({\nu_2\over 2} +N\right)\tau(x)  \right] \right|< \Im \tau(x)
$$ 
Namely: 
\be
   (\Re \nu_2 +2+2N) \ln{|x|\over 16} +O(x)< \Im \nu_2~\arg x +\pi ~\Im \nu_1< 
 (\Re \nu_2 -2+2N) \ln{|x|\over 16}+O(x)
\label{nuovadef}
\ee
Note that 
this defines a domain ${\cal D}(r;\nu_1,\nu_2+2N)$ which  contains the 
 union of 
 ${\cal D}_1(r<1;\nu_1,\nu_2+2N)$,  ${\cal D}_2(r<1;\nu_1,\nu_2+2N)$,  
${\cal D}_1(r<1;\nu_1,\nu_2+
2[N+1])$,  ${\cal D}_2(r<1;\nu_1,\nu_2+2[N+1])$ introduced in the proof of
 Theorem 1. 

 The critical behavior along 
$$ 
\arg x =\arg x_0 +{\Re \nu_2 +2N - {\cal V}\over \Im \nu_2} \ln |x|,
~~~ -2 \leq {\cal V} \leq 2,~~~\Im \nu_2\neq 0   
$$ 
is: 
\vskip 0.2 cm 
\noindent  
For $0<{\cal V} <1$ 
 $$ 
y(x) = -{1\over 4} \left[ 
{e^{i\pi \nu_1} \over 16^{\nu_2+2N-1}}
\right]
x^{\nu_2+2N} ~(1+O(x^{\nu_2+2N}))
$$
\vskip 0.2 cm 
\noindent  
For $1<{\cal V} <2$ 
 $$ 
y(x) = -{1\over 4} \left[ 
{e^{i\pi \nu_1} \over 16^{\nu_2+2N-1}}
\right]^{-1}
x^{2-\nu_2-2N} ~(1+O(x^{2-\nu_2-2N}))
$$
\vskip 0.2 cm 
\noindent  
For ${\cal V}=1$ 
$$
y(x) = x~ \sin^2\left(i{1-\nu_2-2N\over 2} \ln{x\over 16} +{\pi \nu_1
\over 2}\right)(1+O(x))
$$

\vskip 0.2 cm 
\noindent  
For ${\cal V}=0$
$$
  y(x)=\left[ 
{x\over 2} +\sin^{-2} \left( 
-i{\nu_2+2N\over 2} \ln{x\over 16} +{\pi \nu_1\over 2} - i {\nu_2+2N\over 2} 
{F_1(x) \over F(x)}
\right)
\right](1+O(x))
$$
 
\vskip 0.2 cm 
\noindent  
For ${\cal V}=2$
$$
  y(x)=\left[ 
{x\over 2} +\sin^{-2} \left( 
i{2-\nu_2-2N\over 2} \ln{x\over 16} +{\pi \nu_1\over 2} +i {2-\nu_2-2N\over 2} 
{F_1(x) \over F(x)}
\right)
\right](1+O(x))
$$

\vskip 0.2 cm 
\noindent  
For $-1<{\cal V}<0$: it is like $1<{\cal V}<2$ with  $N \mapsto N+1$.

\vskip 0.2 cm 
\noindent  
For $-2<{\cal V}<-1$: it is like $0<{\cal V}<1$ with  $N \mapsto N+1$.

\vskip 0.2 cm 
\noindent  
For ${\cal V}=-1$: it is like ${\cal V}=1$ with $N \mapsto N+1$.

\vskip 0.2 cm 
\noindent  
For ${\cal V}=-2$: it is like ${\cal V}=0$ with $N \mapsto N+1$. 

 \vskip 0.2 cm 
We observe that the choice of $N$ is arbitrary, therefore the {\it same } transcendent has different critical behaviors on different domains (\ref{nuovadef}) 
specified by different values of $N$. This answers Question 1.

\vskip 0.3 cm 
\noindent
{\it Remark:} Note that in the cases ${\cal V}=-2,0,2$, the $\sin^2(...)$ may vanish, therefore there may be (movable) poles. Actually, this happens because the domain now is bigger than 
 that of the generic case since we did not impose 
that $\sin^2(...)\neq 0$, as we did in the proof of Theorem 1. This is also 
the reason why we have
 to keep $F_1(x)\over F(x)$ in the argument of $\sin^2(...)$ in the 
Fourier expansion. 

\vskip 0.3 cm

\vskip 0.2 cm 
If $\Im \nu_2=0$, we choose the convention $ 
 0 \leq \nu_i <2
$. 
The critical behavior for $0<\nu_2<1$ is the same of the case $\Im \nu_2
\neq 0$ with $N=0$ and  $0<{\cal V}<1$; for $1<\nu_2<2$ it 
is the same of the case $\Im \nu_2\neq 0$ with $N=0$ and  $1<{\cal V}<2$. Finally, 
$$
y(x) = (1+O(x))~x~\sin^2\left({\pi \nu_1\over 2}\right),~~~
\hbox{ if } \nu_2=1
$$
$$
    y(x) = \left[{x\over 2} + \sin^{-2}\left({\pi \nu_1 \over 2}
\right) \right](1+O(x)),~~~\hbox{ if }\nu_2=0,~~~\nu_1\neq 0
$$

\vskip 0.3 cm 

Question 2) can be answered. Actually 
$$
  y(x) =\wp\left( \nu_1 \omega^{(0)}_1(x) +\nu_2\omega^{(0)}_2(x);\omega^{(0}_1,\omega^{(0)}_2 \right)+{1+x\over 3}
$$
$$
=\wp\left( \nu_2 \omega^{(1)}_1(x) +\nu_1\omega^{(1)}_2(x);\omega^{(0}_1,\omega^{(0)}_2 \right)+{1+x\over 3}
$$
$$
=\wp\left( \nu_1 \omega^{(\infty)}_1(x) +(\nu_2-\nu_1)
\omega^{(\infty)}_2(x);\omega^{(0}_1,\omega^{(0)}_2 \right)+{1+x\over 3}
$$
Namely, the couples 
$(\nu_1^{(i)},\nu_2^{(i)})$, $i=0,1,\infty$, are: 
$$
 (\nu_1,\nu_2),~~~(\nu_2,\nu_1),~~~(\nu_1,\nu_2-\nu_1)
$$


\subsection{Solution of the Connection Problem 
 in a Non-Generic Case}\label{so2qng}

 We first give  
  the solutions of Questions 1) and 2) in the special 
case $\beta=\gamma=1-2\delta=0$ 
 studied in \cite{DM} and in \cite{guz1},  because of 
its connection to 2-D topological field theory and to Frobenius 
Manifolds \cite{Dub1} \cite{Dub2} and quantum cohomology \cite{KM} \cite{Manin} \cite{guz2}. 
 This sub-section mainly reviews the results of \cite{guz1} and connects them to the framework of the present paper. 
In the next sub-section we will consider the generic case. 

 The Painlev\'e VI equation is the 
isomonodromic
 deformation equation of the fuchsian system \cite{JM1}
\be
{dY\over dz} = \left( 
{A_0(x)\over z}+{A_1(x)\over z-1}+{A_x(x) \over z-x} 
\right)~ Y \equiv A(z,x)~Y
\label{LAB}
\ee
where  $A_i(x)$ ($i=0,1,x$) are $2 \times 2$ matrices.
 They depend on $x$ in 
such a way that the monodromy matrices at $z=0,1,x$, do not 
change for small deformations of $x$. Moreover,   
\be
A_0(x)+A_1(x)+A_x(x) = -{1\over 2} \pmatrix{ \theta_{\infty} & 0 \cr 
                                                 0       & -\theta_{\infty}\cr
}  , 
 ~~~~ \hbox{eigenvalues of } A_i(x)=\pm {1\over 2} \theta_i,~~~~i=0,1,x
\label{CONDIZIONI}
\ee
The constants $\theta_{\nu}$, $\nu=0,1,x,\infty$ are linked to the parameters of PVI: 
$$ 
    \alpha= {1\over 2} (\theta_{\infty} -1)^2,
~~\beta=-{1\over 2} \theta_0^2, 
~~ \gamma={1\over 2} \theta_1^2,
~~\delta={1\over 2} (1-\theta_x^2) 
$$
 A Painlev\'e transcendent is the solution $z=z(x)$ of 
$A(z,x)_{12}=0$: actually 
 $A_{12}(z,x)= \kappa(x) (z-y(x))/z(z-1)(z-x)$, where the function $\kappa(x)$ does not interests us now (see \cite{JM1}). The monodromy 
matrices give a representation of the fundamental group of ${\bf P}^1 \backslash \{0,1,x,\infty\}$. They are usually denoted by $M_{\nu}$, $\nu=0,1,x,\infty$ and they correspond 
to  counterclockwise loops around $z=0,1,x,\infty$ respectively. We order the 
 loops  as in figure \ref{figure3} (an arbitrary base point is chosen).
The monodromy matrices  are not independent, because 
$$ 
 M_1M_xM_0= M_{\infty}
$$
Moreover
$$
M_{\infty}~ \hbox{ is similar to }~
 e^{2\pi i\pmatrix{ -{\theta_{\infty}\over 2} & 0 \cr 
0 & {\theta_{\infty}\over 2}\cr}} ~ e^{2\pi i R},~~~~~ 
   R=\left\{ \matrix{ 0 ~~~~~\hbox{ if } \theta_{\infty}\not \in {\bf Z} 
\cr
\cr
             \pmatrix{ 0 & b \cr 0 & 0} ~~~~~  \hbox{ if } 
0<\theta_{\infty} \in {\bf Z} 
\cr
\cr
 \pmatrix{ 0 & 0\cr b & 0}  ~~~~~ \hbox{ if } 
0>\theta_{\infty} \in {\bf Z} \cr
}\right.
$$ 
where $b$ is a complex number. 
\begin{figure}
\epsfxsize=12cm
\centerline{\epsffile{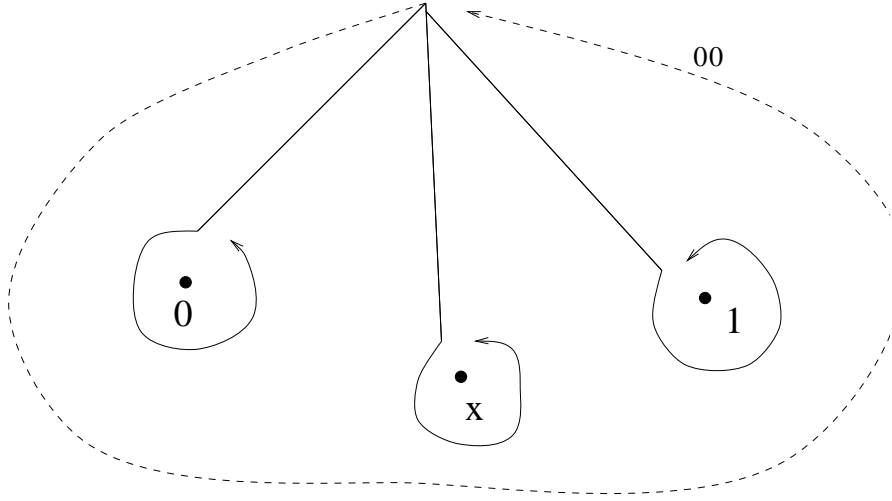}}
\caption{The order of the basis of loops of the Fuchsian system.}
\label{figure3}
\end{figure}
Following  \cite{DM}, let  $\alpha=:{(2\mu-1)^2\over 2}$ (namely 
$\theta_{\infty} = 2\mu$, $\theta_0=\theta_x=\theta_1=0$) and:   
$$ 
   2-x^2_0=\hbox{tr} ~M_0M_x,~~~2-x^2_1=\hbox{tr}~ M_1 M_x,~~~2-x^2_{\infty}=\hbox{tr}~ M_0 M_1.
$$
They satisfy 
\be
x_0^2+x_1^2+x_{\infty}^2-x_0x_1x_{\infty}=4\sin^2(\pi \mu).
\label{ziononsa}
\ee
 For any fixed $\mu$, 
it is proved in \cite{DM} 
 that there  is a one to one correspondence between   
 Painlev\'e 
transcendents and   triples
 of monodromy data $(x_0,x_1,x_{\infty})$, 
defined up to the change of two signs,  satisfying $
x_i \neq \pm 2$, $i=0,1,\infty$ 
and at most one $x_i=0$. 
 A transcendent in one to one correspondence to $(x_0,x_1,x_{\infty})$  
will  be denoted  $y=y(x;x_0,x_1,x_{\infty})$.

\vskip 0.2 cm
Let $\sigma$ be a complex number such that  $2\cos(\pi \sigma)=\hbox{tr}(M_0 
M_x)$, namely: 
$$ 
   \cos \pi \sigma = 1- {x_0^2\over 2},~~~~
0\leq \Re \sigma \leq 1, ~~~~\sigma\neq 1
  $$
We define a domain in the universal covering  $\widetilde{{\bf C}_0}$ of ${\bf C}\backslash \{0\}$, for small $|x|$, as follows:  let $\epsilon$ be a  small positive number, $0<\tilde{\sigma}<1$ a real 
number 
 arbitrarily close to 1 and $\vartheta_1$, $\vartheta_2$ two real parameters. We define
\be
   D(\epsilon;\sigma):= 
\{{ x}\in \widetilde{{\bf C}_0}
   \hbox{ s.t. } |x|<\epsilon ,~~ e^{-\vartheta_1 \Im \sigma}
   |x|^{\tilde{\sigma}} \leq |{x}^{\sigma}| 
\leq e^{-\vartheta_2 \Im
   \sigma} ,~~~0<\tilde{\sigma}<1
   \}. 
\label{dodoalmuseodivienna}
\ee 
If $0\leq \sigma <1$ we simply let $D(\epsilon;\sigma):= 
\{{ x}\in \widetilde{{\bf C}_0}
   \hbox{ s.t. } |x|<\epsilon \}$.  The critical behavior of  $y=y(x;x_0,x_1,x_{\infty})$ 
  was obtained  in \cite{guz1} in the domains $D(\epsilon_n;\pm\sigma+2n)$ for any integer $n$. We note that such domains 
can be rewritten as
$$\Re (\pm \sigma +2n) \ln |x| +\vartheta_2~ \Im(\pm \sigma) \leq \Im (\pm \sigma) \arg x < 
[\Re( \pm \sigma +2n)-\tilde{\sigma}]~\ln|x| + \vartheta_1~\Im(\pm \sigma),
~~~~|x|<\epsilon_n
$$
As it is proved in \cite{guz1},  $\epsilon_n<1$ is small (it decreases as $n$ increases) 
 and  $\Im(\pm \sigma)~\vartheta_1$ decreases  
 as $-\ln| n|$. 
The critical behavior in the domains $D(\epsilon_n;\pm\sigma+2n)$ is:
\be 
  y(x;x_0,x_1,x_{\infty})= a(\pm\sigma+2n;x_0,x_1,x_{\infty})
 ~x^{1-[\pm\sigma+2n]} (1+O(x^{\pm\sigma+2n}
+x^{1-[\pm\sigma+2n]}))
\label{guzzo}
\ee
where 
$$
 a(\sigma;x_0,x_1,x_{\infty})={ i 16^{\sigma} \Gamma\left({\sigma+1\over 2}\right)^4 \over 8 
\sin(\pi \sigma) \Gamma\left(1-\mu+{\sigma\over 2}\right)^2\Gamma\left(
\mu+{\sigma\over 2}
\right)^2}\left[2(1+e^{-i\pi \sigma})- \right.
$$
\be \left. - f(x_0,x_1,x_{\infty})(x_{\infty}^2+
e^{-i\pi\sigma}x_1^2)\right]f(x_0,x_1,x_{\infty})
\label{AHAHAH}
\ee
and
$$
f(x_0,x_1,x_{\infty}):= {4-x_0^2\over x_1^2+x_{\infty}^2-x_0x_1x_{\infty}}
$$
The above formula for $a$ has some limit cases when $\sigma\to 0, \pm 2\mu+2m$, $m$ integer. Namely: 

\vskip 0.2 cm 
\noindent
I) For $x_0=0$. 
           $$ \sigma=0,~~~~
    a= {x_{\infty}^2\over x_1^2+x_{\infty}^2}.
$$
provided that $x_1\neq 0$ and $x_{\infty}\neq 0$.   
\vskip 0.15 cm 
\noindent
II) For  $x_0^2=4 \sin^2(\pi\mu)$ (and then  
 $x_{\infty}^2=-x_1^2 ~\exp(\pm 2\pi i \mu)$ as it follows from (\ref{ziononsa})),  
 there are four sub-cases:

\vskip 0.15 cm
   II1)  $x_{\infty}^2=-x_1^2 e^{- 2
\pi i \mu}$
$$\sigma=  2\mu + 2m,~~~~ a=-{1\over 4 x_1^2} 
               { 16^{2\mu+2m}
 \Gamma(\mu+m+{1\over 2})^4\over \Gamma(m+1)^2 \Gamma(2\mu+m)^2  } ,~~~
m=0,1,2,...
$$

\vskip 0.15 cm
  II2) $x_{\infty}^2=-x_1^2 e^{2\pi i \mu}$
$$\sigma=2\mu+2m,~~~~ m=-1,-2,-3,...$$
$$
            a=-{\cos^4(\pi\mu)\over 4 \pi^4}
                 16^{2\mu+2m} \Gamma(\mu+m+{1\over 2})^4
\Gamma(-2\mu-m+1)^2 \Gamma(-m)^2~ x_1^2 
$$

\vskip 0.15 cm
  II3) $x_{\infty}^2=-x_1^2e^{2\pi i\mu}$
$$\sigma=-2\mu+2m,~~~~a=-{1\over 4 x_1^2}
         {
  16^{-2\mu+2m} \Gamma(-\mu+m+{1\over 2})^4\over \Gamma(m-2\mu+1)^2
\Gamma(m)^2  },,~~~~m=1,2,3,...
$$

\vskip 0.15 cm
II4)  $x_{\infty}^2=-x_1^2 e^{-2\pi i \mu}$
$$
\sigma=-2\mu+2m,~~~~m=0,-1,-2,-3,...
$$
$$
  a= -{\cos^4(\pi\mu)\over 4 \pi^4} 16^{-2\mu+2m}
\Gamma(-\mu+m+{1\over 2})^4\Gamma(2\mu-m)^2 \Gamma(1-m)^2 ~ x_1^2
$$

\vskip 0.2 cm 
Conversely, let be given a transcendent  with 
behavior $y(x)=a x^{1-\sigma}(1+~\hbox{higher orders})$ as $x\to 0$, for given $\sigma$ and $a$. We define a triple $(x_0,x_1,x_{\infty})$ by the formulae below: 
$$
  x_0= 2 \sin \left({\pi \over 2} \sigma \right)
$$
$$
  x_1= i\left( {\sqrt{a}\over f(\sigma,\mu) G(\sigma,\mu)}-{G(\sigma,\mu)
\over \sqrt{a}} 
\right)
$$
$$
  x_{\infty}= {e^{i\pi \sigma\over 2}\sqrt{a}\over f(\sigma,\mu) G(\sigma,\mu)}+{G(\sigma,\mu)
\over e^{i\pi \sigma\over 2}\sqrt{a}}
$$
where
$$ 
f(\sigma,\mu):= {2 \cos^2\left({\pi \over 2} \sigma\right) \over \cos(\pi \sigma) 
- \cos(2\pi \mu)}\equiv f(x_0,x_1,x_{\infty}),~~~~G(\sigma,\mu):= {4^{\sigma}\Gamma\left({\sigma+1\over 2}\right)^2 \over 2 \Gamma\left(1-\mu+{\sigma\over 2}\right)\Gamma\left(
\mu+{\sigma\over 2}
\right)}
$$
Again, there are limiting cases:

 \vskip 0.15 cm
\noindent
 i) $\sigma =0$
$$\left\{ 
              \matrix{
                         x_0=0  \cr\cr
                        x_1^2=2\sin(\pi\mu)~\sqrt{1-a} \cr\cr
                        x_{\infty}^2=2\sin(\pi\mu)~\sqrt{a}
                        \cr
                       }
 \right.  
$$
 
\vskip 0.15 cm
\noindent
 ii) $\sigma= \pm 2\mu +2 m$.  

\vskip 0.15 cm
  ii1) $\sigma=2\mu+2m$, $m=0,1,2,...$
 $$ 
         \left\{ 
   \matrix{ 
              x_0=2\sin(\pi\mu)\cr \cr
              x_1= -{i\over 2}{16^{\mu+m} \Gamma(\mu+m+{1\over 2})^2 \over
   \Gamma(m+1) \Gamma(2\mu+m)}~{1\over\sqrt{a}} 
                                           \cr \cr
               x_{\infty}=i~x_1~e^{-i\pi\mu}\cr
          }
                  \right.
$$

 ii2) $\sigma=2\mu+2m$, $m=-1,-2,-3,...$
$$
  \left\{ 
\matrix{
            x_0=2\sin(\pi\mu)\cr \cr
            x_1=2i{\pi^2\over \cos^2(\pi\mu)}
                {1\over 16^{\mu+m} \Gamma(\mu+m+{1\over 2})^2
\Gamma(-2\mu-m+1) \Gamma(-m)}~\sqrt{a}
\cr\cr
x_{\infty}=-ix_1 e^{i\pi\mu}
}
\right.
$$

ii3) $\sigma=-2\mu +2m$, $m=1,2,3,...$
$$
   \left\{ 
\matrix{
           x_0=-2\sin(\pi\mu) \cr\cr
         x_1= -{i\over 2}{16^{-\mu+m} \Gamma(-\mu+m+{1\over 2})^2\over
\Gamma(m-2\mu+1) \Gamma(m)}~{1\over \sqrt{a}} \cr\cr
        x_{\infty}= ix_1e^{i\pi\mu}\cr
}
\right.
$$

ii4) $\sigma= -2\mu+2m$, $m=0,-1,-2,-3,...$
$$
  \left\{
\matrix{
            x_0=-2\sin(\pi\mu) \cr\cr
            x_1= 2i {\pi^2\over \cos^2(\pi\mu)}{1\over 16^{-\mu+m}
\Gamma(-\mu+m+{1\over 2})^2\Gamma(2\mu-m) \Gamma(1-m)}  ~\sqrt{a}
\cr\cr
x_{\infty}=-ix_1 e^{-i\pi \mu}\cr
}
 \right.
$$
In all the above formulae the sign of the square roots is arbitrary (a triple is define up to the change of two signs). The relation
$x_0^2+x_1^2+x_{\infty}^2-x_0x_1x_{\infty}=4 \sin^2(\pi \mu)$  is
automatically satisfied. Note that $\sigma \neq 1$ implies $x_0\neq \pm
2$. 

 It is proved in \cite{guz1} that if the triple $(x_0,x_1,x_{\infty})$ defined above is such that 
$x_i\neq \pm 2$ ($i=0,1,\infty$) and at most one $x_i=0$, then there exists only one transcendent with behavior $y(x)=a x^{1-\sigma}(1+~\hbox{higher orders})$ as $x\to 0$ specified by 
the given $a$ and $\sigma$.  It coincides with 
$y(x;x_0,x_1,x_{\infty})$.

 \vskip 0.3 cm

 Let us return to the elliptic representation. 
We observe that in the special case $\beta=\gamma=1-2\delta=0$ there is only  
${\partial \over \partial u}\wp(u/2)$ 
 in the r.h.s. of (\ref{fuchseq}). 
Therefore, $\varepsilon_1=\varepsilon_2=0$ in the proof of Theorems 1 and 3. 
As a consequence, not only at $x=0$, but also at $x=1,\infty$ the domain ${\cal D}$ is 
larger.  More precisely, at $x=0$ 
the result of section \ref{special} applies: 
{\it For any $N\in {\bf Z}$, and for any complex $\nu_1^{(0)},
\nu_2^{(0)}$, such that $0<\nu_2^{(0)}<2$ if $\nu_2^{(0)}$ is real, there 
exists a  transcendent $y=\wp \bigl(\nu_1^{(0)}\omega_1^{(0)}
+\nu_2^{(0)} \omega_2^{(0)} 
+ v^{(0)} \bigr) +{1+x\over 3}$ such that $v^{(0)}(x)$ is holomorphic in 
$$
{\cal D}(r_N;\nu^{(0)}_1,\nu^{(0)}_2+2N):= \left\{ 
x\in \tilde{{\bf C}_0} ~|~ |x|<r_N,~\left| e^{-i\pi \nu^{(0)}_1} \left({x\over 16} 
\right)^{2-\nu^{(0)}_2-2N}\right|<r,~\left| e^{i\pi \nu^{(0)}_1} \left({x\over 16} 
\right)^{\nu^{(0)}_2+N}\right|<r
\right\}
$$ 
where it has the convergent expansion
$$
 v(x;\nu^{(0)}_1,\nu^{(0)}_2+2N) = 
\sum_{n\geq 1} a_n x^n + \sum_{n\geq 0,m\geq 1} b_{nm}
x^n\left[e^{-i\pi \nu^{(0)}_1}\left({x\over 16} \right)^{2-\nu^{(0)}_2-2N} \right]^m
+
$$
$$
+
 \sum_{n\geq 0,m\geq 1} c_{nm}
x^n\left[e^{i\pi \nu^{(0)}_1}\left({x\over 16} \right)^{\nu^{(0)}_2+2N} \right]^m
$$
}

At $x=1$. {\it For any $N\in {\bf Z}$, and for any complex $\nu_1^{(1)},
\nu_2^{(1)}$, such that $0<\nu_2^{(1)}<2$ if $\nu_2^{(1)}$ is real, there 
exists a  transcendent $y=\wp \bigl(\nu_1^{(1)}\omega_1^{(1)}
+\nu_2^{(1)} \omega_2^{(1)} 
+ v^{(1)} \bigr) +{1+x\over 3}$ such that $v^{(1)}(x)$ is holomorphic in  
$$
{\cal D}(r_N;\nu_1^{(1)},\nu_2^{(1)}+2N) := 
\left\{
x\in \widetilde{{\bf P}^1_1}~\hbox{ such that }~|1-x|<r_N, ~\left|e^{i\pi \nu_1^{(1)}} \left( {1-x\over 16}\right)^{2-\nu^{(1)}_2
-2N}\right|<r_N,\right.
$$
$$ \left. \left|e^{-i\pi \nu_1^{(1)}} \left( {1-x\over 16}\right)^{\nu^{(1)}_2+2N} \right|<r_N 
\right\}
$$
where it has a  convergent expansion
$$
  v^{(1)}(x;\nu_1^{(1)},\nu_2^{(1)}+2N) = 
\sum_{n\geq 1} a_n (1-x)^n+\sum_{n\geq 0,m\geq 1} 
 b_{nm} (1-x)^n \left[e^{i\pi \nu_1^{(1)}} \left( {1-x\over 16}\right)^{2-\nu^{(1)}_2
-2N} \right]^m +
$$
$$
 +\sum_{n\geq 0,m\geq 1} 
 c_{nm} (1-x)^n \left[e^{-i\pi \nu_1^{(1)}} \left( {1-x\over 16}\right)^{\nu^{(1)}_2+2N} \right]^m 
$$
}

At $x=\infty$. {\it For any $N\in {\bf Z}$, and for any complex $\nu_1^{(\infty)},
\nu_2^{(\infty)}$, such that $0<\nu_2^{(\infty)}<2$ if $\nu_2^{(\infty)}$ is real, there exists a  transcendent $y=\wp \bigl(\nu_1^{(\infty)}\omega_1^{(\infty)}
+\nu_2^{(\infty)} \omega_2^{(\infty)} 
+ v^{(\infty)} \bigr) +{1+x\over 3}$ such that $v^{(\infty)}(x)$ 
is holomorphic in  
$$
{\cal D}(r_N;\nu_1^{(\infty)},\nu_2^{(\infty)}+2N) := 
\left\{
x\in \widetilde{{\bf P}_{\infty}}~\hbox{ such that }~|x^{-1}|<r_N, ~\left|
e^{-i\pi \nu_1^{(\infty)}} \left( {16\over x}\right)^{2-\nu^{(\infty)}_2
-2N}\right|<r_N,\right.$$
$$
\left.
 \left|e^{i\pi \nu_1^{(\infty)}} \left( {16\over 
x}\right)^{\nu^{(\infty)}_2+2N} \right|<r_N 
\right\}
$$
where it has a convergent expansion 
$$
 x^{1\over 2} 
v^{(\infty)}(x;\nu_1^{(\infty)},\nu_2^{(\infty)}+2N) = \sum_{n\geq 1} 
a^{(1)} \left({1\over x}\right)^n
+\sum_{n\geq0,m\geq 1} 
 b_{nm} \left({1\over x}\right)^n \left[e^{-i\pi \nu_1^{(\infty)}} 
\left( {16\over x}\right)^{2-\nu^{(\infty)}_2
-2N} \right]^m +$$
$$
+\sum_{n\geq 0,m\geq 1} 
 c^{(1)}
 \left({1\over x}\right)^n \left[e^{i\pi \nu_1^{(\infty)}} 
\left( {16\over x}\right)^{\nu^{(\infty)}_2+2N} \right]^m 
$$
}
\vskip 0.15 cm 
 Note that we have used the same notations for the coefficients $a_n$, $b_{nm}$, $c_{nm}$ both at $x=0,1$ and $x=\infty$. This is for simplicity of 
notations, but they are  different! 

\vskip 0.3 cm 
 We are ready to identify the transcendents $y(x;x_0,x_1,x_{\infty})$ with the transcendents in elliptic representation. We know that for $x\to 0$ the transcendent  
$y(x)= \wp\bigl(\nu_1\omega_1^{(0)}
+\nu_2 \omega_2^{(0)} 
+ v^{(0)}(x;\nu_1,\nu_2+2N) \bigr) +{1+x\over 3}$   has behavior 
$$
  y(x)=
    - { 1 \over 4}\left[ 
{e^{i\pi \nu_1} \over 16^{\nu_2+2N-1}}
\right]
x^{\nu_2+2N}
(1+\hbox{ higher orders}),~~~0<{\cal V}<1 $$
$$=-{1\over 4}\left[
 {e^{i\pi \nu_1} \over 16^{\nu_2+2N-1}}
\right]^{-1}
x^{2-\nu_2-2N}(1+\hbox{ higher orders}),~~~1<{\cal V}<2 
$$
On the other hand, a transcendent $y(x;x_0,x_1,x_{\infty})$ has behavior
(\ref{guzzo}):
$$
  y(x;x_0,x_1,x_{\infty})= a x^{1-[\pm\sigma+2n]}(1+\hbox{ higher orders})
$$
As we mentioned above, the critical  behavior uniquely determines the transcendent. Therefore, 
 the transcendent  
in the
 elliptic representation, for given $\nu_1,\nu_2,N$,  
must coincide with some $y(x;x_0,x_1,x_{\infty})$ having the same critical behavior. 
Let us choose the convention 
$$0
\leq \Re \sigma \leq 1, ~~~
0\leq \Re \nu_2 <2
$$
We take $N=n=0$.  
A necessary condition for the transcendents to coincide  is that 
$$ 
  D(\epsilon_0;\sigma)\cap {\cal D}(r_0;\nu_1,\nu_2) \neq \emptyset
$$   
They actually coincide when they have the same critical behavior, namely if: 
$$
1):~~~  \sigma=1-\nu_2,~~~~ a= -{1\over 4} \left[{e^{i\pi \nu_1} \over 16^{\nu_2-1}} \right] ,~~~
\hbox{ if }0\leq \Re \nu_2 \leq 1
$$
$$
2):~~~\sigma=\nu_2-1,~~~~ a= -{1\over 4} \left[{e^{i\pi \nu_1} \over 16^{\nu_2-1}} \right]^{-1}
,~~~\hbox{ if }1\leq \Re \nu_2 < 2
$$
(recall that $\sigma\neq 1$ and $\nu_2\neq 0,2$). Equivalently: 
 $$ 
1):~~  \nu_2=1-\sigma,~~~e^{i\pi \nu_1} = - 4 a(\sigma) ~16^{-\sigma}~~\hbox{ if } 
0\leq \Re \nu_2 \leq 1
$$
$$
2):~~\nu_2=1+\sigma,~~~e^{-i\pi \nu_1} = - 4 a(\sigma) ~16^{-\sigma}~~\hbox{ if } 
1\leq \Re \nu_2 \leq 2
$$

\vskip 0.3 cm
We are now ready  to answer Question 1) and Question 2). Actually, this is possible 
because the identifications 1) and 2) above allow
 us to express the parameters $\nu_1,\nu_2$ in terms 
of the triple of monodromy data $(x_0,x_1,x_{\infty})$ and $\mu$.

 We observe that  transformation $\sigma\mapsto -\sigma$ is $\sigma=1- \nu_2 \mapsto 
\nu_2-1$ in case 1), and  $\sigma=\nu_2-1 \mapsto 
1-\nu_2$ in case 2). Namely, in case 1) the behavior 
of $y(x;x_0,x_1,x_{\infty})$ is 
$$ 
y(x;x_0,x_1,x_{\infty}) = a(\sigma;x_0,x_1,x_{\infty}) x^{1-\sigma} (1 + \hbox{ higher orders}),~~\hbox{ in } 
D(\epsilon_0;\sigma)
$$
and 
$$ 
 y(x;x_0,x_1,x_{\infty})= a(-\sigma;x_0,x_1,x_{\infty}) x^{1-[-\sigma]} (1 + \hbox{ higher orders}),~~\hbox{ in } 
D(\epsilon_0;-\sigma)
$$
This corresponds to the  behaviors obtained in elliptic representation: 
 $$
  y(x)=
    - { 1 \over 4}\left[ 
{e^{i\pi \nu_1} \over 16^{\nu_2-1}}
\right]
x^{\nu_2}
(1+\hbox{ higher orders}),~~~0<{\cal V}<1 $$ 
and 
$$
y(x)=-{1\over 4}\left[
 {e^{i\pi \nu_1} \over 16^{\nu_2-1}}
\right]^{-1}
x^{2-\nu_2}(1+\hbox{ higher orders}),~~~1<{\cal V}<2 
$$
respectively, 
with  $\cos \pi \nu_2={x_0^2\over 2}-1$ and 
$$
  e^{i\pi \nu_1}= 
-{i \Gamma^4\left( 1-{\nu_2\over 2}\right)\over 2 \sin(\pi \nu_2)\Gamma^2\left({3\over 2} 
-\mu - {\nu_2\over 2}\right) \Gamma^2\left({1\over 2} +\mu -{\nu_2\over 2}
\right)}
\left[2(1-e^{i\pi \nu_2}) -\right.
$$
\be
\left.  - f(x_0,x_1,x_{\infty})(x_{\infty}^2- e^{i\pi \nu_2}
x_1^2) \right]f(x_0,x_1,x_{\infty}) 
\label{cambiolavoro1}
\ee
 In case 2), the behavior of    $y(x;x_0,x_1,x_{\infty})$ is still 
$$ 
y(x;x_0,x_1,x_{\infty}) = a(\sigma;x_0,x_1,x_{\infty}) x^{1-\sigma} (1 + \hbox{ higher orders}),~~\hbox{ in } 
D(\epsilon_0;\sigma)
$$
and  
$$ 
y(x;x_0,x_1,x_{\infty}) = a(-\sigma;x_0,x_1,x_{\infty}) x^{1-[-\sigma]} (1 + 
\hbox{ higher orders}),~~\hbox{ in } 
D(\epsilon_0;-\sigma)
$$
This corresponds to 
$$y(x)=-{1\over 4}\left[
 {e^{i\pi \nu_1} \over 16^{\nu_2-1}}
\right]^{-1}
x^{2-\nu_2}(1+\hbox{ higher orders}),~~~1<{\cal V}<2 
$$
and 
$$
  y(x)=
    - { 1 \over 4}\left[ 
{e^{i\pi \nu_1} \over 16^{\nu_2-1}}
\right]
x^{\nu_2}
(1+\hbox{ higher orders}),~~~0<{\cal V}<1 
$$
respectively,  
with  $\cos \pi \nu_2={x_0^2\over 2}-1$ and 
$$
  e^{-i\pi \nu_1}= 
{i \Gamma^4\left( {\nu_2\over 2}\right)\over 2 \sin(\pi \nu_2)\Gamma^2\left({1\over 2} 
-\mu + {\nu_2\over 2}\right) \Gamma^2\left(-{1\over 2} +\mu +{\nu_2\over 2}
\right)}
\left[2(1-e^{-i\pi \nu_2})-\right.
$$
\be
\left. 
 - f(x_0,x_1,x_{\infty})(x_{\infty}^2- e^{-i\pi \nu_2}
x_1^2) \right]f(x_0,x_1,x_{\infty})
\label{cambiolavoro2}
\ee

\vskip 0.3 cm
\noindent
{\bf Answer to Question 1)}.  
 Let us  study the effect of the shift 
$\sigma\mapsto \sigma+2n$ on the elliptic representation 
of $y(x;x_0,x_1,x_{\infty})$.  
$y(x;x_0,x_1,x_{\infty})$ has critical  behaviors on  the domains 
$D(\epsilon_n;\sigma+2n)$, with exponents $\sigma+2n$ and coefficients 
$a=a(\sigma+2n;x_0,x_1,x_{\infty})$. The shift corresponds to 
 $ \nu_2 \mapsto \nu_2-2n$ in case 1), and $ \nu_2 \mapsto \nu_2
+2n$ in case 2).  Using the properties of the  $\Gamma$-function and 
(\ref{cambiolavoro1}),(\ref{cambiolavoro2})  we obtain, both for case 1) and 
 case 2):
$$
\bigl.e^{i\pi \nu_1}\bigr|_{\nu_2+2n} =\bigl.e^{i\pi \nu_1}\bigr|_{\nu_2}~
K(\nu_2,n)
$$
where 
$$
K(\nu_2,n):=
\left\{ \matrix{ 1,~~~~n=0 
\cr
\cr
 \Pi_{k=1}^{n} {[\nu_2-1+2\mu+2(k-1)]^2
[\nu_2-1-2\mu+2k]^2 \over (\nu_2+2k-1)^4},~~~~n>0
\cr
\cr
 \Pi_{k=1}^{|n|}{(\nu_2-2k)^4 \over
[\nu_2-1+2\mu-2k]^2[\nu_2-1-2\mu-2(k-1)]^2},~~~~n<0 
}
 \right.
$$
Thus
$$ 
  \bigl.\nu_1\bigr|_{\nu_2+2n}= \bigl.
\nu_1\bigr|_{\nu_2}-{i\over \pi} \ln K(\nu_2,n)
$$

 In this way we have answered to Question 1): the representations 
$$ 
y(x)= \wp\left(\bigl[\nu_1-{i\over \pi}\ln K(\nu_2,N)\bigr]~ \omega_1
+\nu_2~\omega_2+v\bigl(x;\nu_1-{i\over \pi}\ln K(\nu_2,N),\nu_2+2N\bigl)
\right)+{1+x\over 3}
$$
define the same  transcendent  on different domains ${\cal D}
\left(r_N; \nu_1-{i\over \pi}\ln K(\nu_2,N),\nu_2+2N\right)$.

\vskip 0.3 cm 
\noindent
{\it Remark:} We now know the elliptic representation and the critical behavior of 
the transcendent $y(x;x_0,x_1,x_{\infty})$ 
 on the union of the domains  ${\cal D}
(r_N; \nu_1-{i\over \pi}\ln K(\nu_2,N),\nu_2+2N)$, $n\in {\bf Z}$ (for $\Im \nu_2\neq 0$).  
However, some regions in 
the $(\ln |x|,\Im\nu_2 \arg x)$-plane, for $x$ close to 0, 
  are not included in the union. 
In these regions the transcendent  may have movable poles.

\vskip 0.3 cm 
\noindent
{\bf Answer to  Question 2)}. 
 In \cite{DM},\cite{guz1} the connection problem was solved by showing that a given transcendent $y(x;x_0,x_1,x_{\infty})$ has critical behaviors:  
\be
y(x;x_0,x_1,x_{\infty}) = 
 a^{(0)} x^{1-\sigma^{(0)}}(1+ \hbox{ higher orders in $x$}), ~~~x\to 0  
\label{BIEVO1}
\ee
\be=
 1- a^{(1)} (1-x)^{1-\sigma^{(1)}}(1+\hbox{ higher orders in $(1-x)$}),
~~~x\to 1
\label{BIEVO2}
\ee
\be = 
a^{(\infty)} x^{\sigma^{(\infty)}}(1+\hbox{ higher orders in $x^{-1}$}),~~~x\to 
\infty 
\label{BIEVO3}
\ee
where 
$$ 
  \cos(\pi \sigma^{(i)})= 1-{x_i^2\over 2}, ~~~0\leq \Re \sigma^{(i)}\leq 1
,~~~\sigma^{(i)}\neq 1,~~~i=0,1,\infty 
$$   
The parameters $a^{(0)}$, $\sigma^{(0)}$ stand for  $a,\sigma$ used before.  
 $a^{(1)}, a^{(\infty)}$ 
are obtained by the  formula (\ref{AHAHAH}) (or limit cases)    
 with the substitutions $(x_0,x_1,x_{\infty}) 
\mapsto (x_1,x_0,x_0x_1-x_{\infty})$, $\sigma^{(0)}\mapsto \sigma^{(1)}$ and 
 $(x_0,x_1,x_{\infty}) 
\mapsto (x_{\infty}, -x_1,x_0-x_1x_{\infty})$,  $\sigma^{(0)}\mapsto 
\sigma^{(\infty)}$ respectively. The behavior (\ref{BIEVO3}) holds in a domain 
\be 
    D(M; \sigma^{(\infty)}):=\{ 
 {x}\in 
\widetilde{ {\bf P}^1 \backslash \{\infty\} } 
\hbox{ s.t. } 
|x|>M,~ 
e^{ -\vartheta_1\Im\sigma^{(\infty)} }
 |x|^{-\tilde{\sigma}}
\leq |{x}^{ -\sigma^{(\infty)}}|
\leq  e^{-\vartheta_2\Im\sigma^{(\infty)}}
 \}
\label{dodoalmuseodiviennainf}
\ee
where $M>0$ is sufficiently big; the behavior (\ref{BIEVO2}) holds in 
\be
  D(\epsilon;\sigma^{(1)}):= \{ {x}\in
   \widetilde{{\bf C}\backslash \{1\}} 
   \hbox{ s.t. } |1-x|<\epsilon ,~~ e^{-\vartheta_1 \Im \sigma}
   |1-x|^{\tilde{\sigma}} \leq |(1-{x})^{\sigma^{(1)}}| 
 \leq e^{-\vartheta_2 \Im
   \sigma} \}  
\label{dodoalmuseodivienna1}
\ee

 At $x=1$, $x=\infty$ we repeat the analysis we did above for $x=0$, 
 identifying $y(x;x_0,x_1,x_{\infty})$ with a transcendent in elliptic 
representation with the same critical behaviors (by uniqueness of the critical behavior).  
Thus, we conclude that the transcendent $y(x;x_0,x_1,x_{\infty})$ has 
a representation 
$$
  y(x)= \wp (\nu_1^{(0)}\omega_1^{(0)}+\nu_2^{(0)}\omega_2^{(0)}+v^{(0)})
+{1+x\over 3} 
$$
 at $x=0$; it has a representation 
$$
 y(x)= 
   \wp (\nu_1^{(1)}\omega_1^{(1)}+\nu_2^{(1)}\omega_2^{(1)}+v^{(1)})
+{1+x\over 3}
$$
 at $x=1$ and it has a  representation 
$$
y(x)= \wp (\nu_1^{(\infty)}\omega_1^{(\infty)}+\nu_2^{(\infty)}
\omega_2^{(\infty)}+v^{(\infty)})+{1+x\over 3}, 
$$
where the  parameters $\nu_2^{(i)}$ are obtained from 
$$ 
\cos \pi \nu_2^{(i)}= {x_i^2\over 2} -1,~~0\leq \Re \nu_2^{(i)}\leq 1~~(\nu_2^{(i)} \neq 0), ~~~~~\hbox{ {\rm or } }   1\leq \Re \nu_2^{(i)}<2
~~~~i=0,1,\infty
$$ 
and the  parameter $\nu_1^{(0)}$ is 
obtained from (\ref{cambiolavoro1}) for the choice $
0\leq \Re \nu_2^{(0)}\leq 1$ and from  (\ref{cambiolavoro2}) for the choice $
1\leq \Re \nu_2^{(0)}<2$. Moreover,   
$\exp\{ \mp i \pi \nu_1^{(1)}\}$, $\exp\{\pm i \pi \nu_1^{(\infty)}\}$ 
are given by formulae analogous to 
(\ref{cambiolavoro1}), (\ref{cambiolavoro2}) with the substitutions 
$(x_0,x_1,x_{\infty}) 
\mapsto (x_1,x_0,x_0x_1-x_{\infty})$, $\nu_2^{(0)}\mapsto \nu_2^{(1)}$ and 
 $(x_0,x_1,x_{\infty}) 
\mapsto (x_{\infty}, -x_1,x_0-x_1x_{\infty})$,  $\nu_2^{(0)}\mapsto 
\nu_2^{(\infty)}$ respectively.

\vskip 0.3 cm 
 Conversely, given 
\be
y(x)= \wp \bigl({\nu}_1\omega_1^{(0)}(x)+
{\nu}_2\omega_2^{(0)}(x)+v^{(0)}(x;{\nu}_1,
{\nu}_2)\bigr) +
{1+x\over 3},~~~ \hbox{ at } x=0
\label{startingpoint1}
\ee
it coincides 
 with $y(x;~x_0,x_1,x_{\infty})$, with the following monodromy data.

\vskip 0.2 cm 
\noindent 
If $0\leq \Re \nu_2 \leq 1$:  
$$
  x_0=2\cos\left({\pi \over 2}\nu_2\right)
$$
$$
x_1=  \left[{ 4^{-\nu_2}~2~e^{i {\pi\over 2} \nu_1 } \over f(\nu_2,\mu) G(\nu_2,\mu)}+{G(\nu_2,\mu) \over  4^{-\nu_2}~2~e^{i{ \pi\over 2} \nu_1 }} \right]
$$
$$
x_{\infty} = \left[ { 4^{-\nu_2}~2~e^{i{\pi \over 2} (\nu_1-\nu_2)}\over 
  f(\nu_2,\mu) G(\nu_2,\mu)}+ {G(\nu_2,\mu)\over  4^{-\nu_2}~2~e^{i{\pi \over 2} (\nu_1-\nu_2)}} \right]
$$
where 
$$ 
  f(\nu_2,\mu)= - {2\sin^2\left({\pi \over 2} \nu_2\right) \over 
\cos(\pi \nu_2) +\cos(2\pi \mu)}
,~~~G(\nu_2,\mu)=  4^{-\nu_2}~2~ { \Gamma\left(1-{\nu_2\over 2}\right)^2 \over
\Gamma\left({3\over 2} - \mu - {\nu_2 \over 2} \right) 
\Gamma\left({1\over 2} +\mu - {\nu_2\over 2} \right) 
}
$$

\vskip 0.2 cm 
\noindent 
If $1\leq \Re \nu_2 <2$:  
$$
x_0=  2 \cos \left({\pi\over 2} \nu_2 \right)
$$
$$
x_1= \left[ {e^{-i{\pi \over 2} \nu_1} \over  4^{1-\nu_2}~ 2~f(\nu_2,\mu) G_1(\nu_2,\mu)} + {  4^{1-\nu_2}~  2~G_1(\nu_2,\mu) \over e^{-i{\pi \over 2} \nu_1}}
 \right] 
$$
$$ 
x_{\infty}
= \left[{ e^{i {\pi \over 2} (\nu_2-\nu_1)} \over  4^{1-\nu_2}~2~ f(\nu_2,\mu) 
G_1(\nu_2,\mu) } +{   4^{1-\nu_2}~2~ G_1(\nu_2,\mu)\over e^{i {\pi \over 2} (\nu_2-\nu_1)}} \right]
$$
where 
$$ 
  G_1(\nu_2,\mu)= {1\over 4^{1-\nu_2}~2} {\Gamma\left({\nu_2\over 2} \right)^2
\over 
\Gamma\left( {1\over 2}- \mu +{\nu_2\over 2}\right)\Gamma\left(
-{1\over 2} +\mu+{\nu_2\over 2}\right)}
$$
(The limit cases are left as an exercise for the reader). 
After computing the monodromy data, we can write the  elliptic 
representations of  $y(x;x_0,x_1,x_{\infty})$ at $x=1$ and $x=\infty$. 
These are the elliptic representations at $x=1$, $x=\infty$  of 
(\ref{startingpoint1}). Thus,  
we have solved the connection problem for (\ref{startingpoint1}).

\vskip 0.3 cm
 We observed that there  is a one to one correspondence between   
 Painlev\'e 
transcendents and   triples
 of monodromy data $(x_0,x_1,x_{\infty})$, 
defined up to the change of two signs,  satisfying $
x_i \neq \pm 2$, $i=0,1,\infty$ (i.e. $\sigma^{(i)}\neq 1$)  
and at most one $x_i=0$. The cases when  these conditions are not satisfied are studied in \cite{M}. 
However, if  $x_i=\pm2$ (namely the  trace is $ -2$) the problem of 
finding the critical behavior at the corresponding critical point $x=i$ is still 
open (except when  {\it all the 
three} $x_i$  are $\pm 2$: in this case  there is a one-parameter class of 
 solutions called {\it Chazy solutions } in \cite{M}). We conclude that the results of our paper 
(together with \cite{guz1}), 
plus the results of \cite{M} cover all the possible transcendents, except the special case when 
 one or two 
$x_i$ are $\pm 2$. We plan to cover this last case soon.



\subsection{Solution of the Connection Problem in the Generic Case}
\label{Solution of the Connection Problem in the Generic Case}

 We solve the connection problem for the elliptic representation in the 
{\it generic case}: 
$$ 
    \nu_2,~\theta_0,~\theta_x,~\theta_1,~\theta_{\infty}~\not \in {\bf Z};
~~~~
{\pm 1\pm\nu_2 \pm \theta_1 \pm \theta_{\infty}\over 2},~
{\pm 1\pm\nu_2 \pm \theta_0 \pm \theta_x \over 2} ~\not \in {\bf Z}
$$ 
 We show  in the Appendix, following \cite{Jimbo}, that if $\theta_0,\theta_x,\theta_1,\theta_{\infty}\not \in {\bf Z}$ 
there is a one to one correspondence between transcendents and the 
 monodromy data   $\theta_0, \theta_x,\theta_1,\theta_{\infty}, \hbox{tr}(M_0M_x),  \hbox{tr}(M_0M_1),
 \hbox{tr}(M_1M_x)$ of the  fuchsian system (\ref{LAB}). For this reason 
we write     
 \be
y(x)=
y\bigl(x;~ \theta_0, \theta_x,\theta_1,\theta_{\infty}, \hbox{tr}(M_0M_x),  \hbox{tr}(M_0M_1),
 \hbox{tr}(M_1M_x)~ \bigr)
\label{PAALLAA}
\ee

  In the Appendix it is  shown that for {\it any} value of $\alpha$, $\beta$, $\gamma$, $\delta$, for any $a\neq 0$, for any complex 
$\sigma \not \in (-\infty,0]
\cup [1,+\infty)$ and for any additional real parameters $0<\tilde{\sigma}<1$, $\vartheta_1$,$\vartheta_2$ there exists a sufficiently small $\epsilon$ and  a transcendent $
y(x;\sigma,a)$ in the domain  $D(\epsilon;\sigma)$ defined by (\ref{dodoalmuseodivienna}). It has critical behavior 
$$
y(x;\sigma,a)= 
 a x^{1-\sigma}(1+~\hbox{higher orders})
$$
for $x\to 0 $ along a regular path in  $ D(\epsilon;\sigma)$.    

In the Appendix it is also proved that in  the generic case
  a  transcendent (\ref{PAALLAA}) 
associated to the {\it monodromy data}    
 $\theta_0$,$\theta_x$,$\theta_1$,$\theta_{\infty}$, $\hbox{tr}(M_0M_x)$,$\hbox{tr}(M_0M_1)$,
$\hbox{tr}(M_1M_x)$,  
  coincides with a transcendent $y(x;\sigma,a)$ when $x$ lies in the domain$D(\epsilon;\sigma)$.  $\sigma$ and $a$ are obtained from the monodromy data by explicit formulae 
$$
  2 \cos(\pi \sigma) =\hbox{tr}(M_0M_x)
$$
$$
  a=a\bigl(\sigma;\theta_0, \theta_x,\theta_1,\theta_{\infty}, \hbox{tr}(M_0M_x),  \hbox{tr}(M_0M_1),
 \hbox{tr}(M_1M_x) \bigr)
$$
The  formula for $a$ is rather long and complicated, and we do not need its explicit form here; 
therefore we refer to the Appendix for its computation. 
 Note that $\sigma$ is defined up to $ \sigma\mapsto \pm \sigma+2n$, $n $ integer. We can 
fix $
    0\leq \Re \sigma \leq 1
$; thus,  a  transcendent (\ref{PAALLAA}) 
associated to the monodromy data coincides with the transcendents 
$$
y\Bigl(x;\pm \sigma+2n,a\bigl(\pm\sigma+2n;\theta_0, \theta_x,\theta_1,\theta_{\infty}, \hbox{tr}(M_0M_x),  \hbox{tr}(M_0M_1),
 \hbox{tr}(M_1M_x) \bigr)\Bigr)
$$ 
in the 
 domains $D(\epsilon_n; \pm\sigma+2n)$ defined for sufficiently small $\epsilon_n$  (note that only $0<\sigma<1$ is allowed if $\sigma$ is real, and in this 
 case no translation $ \pm\sigma+2n$ is possible). We observe that $a$ depends on the monodromy data, but also explicitly on $\sigma$. Namely, it 
changes if we do $\sigma \mapsto \pm \sigma +2n$ (being the monodromy data 
fixed). 
 In the following, for brevity we write $a=a(\sigma)$, understanding 
 the dependence on the monodromy data. 

 We finally remark that  
in the Appendix the formula 
$ a=a\bigl(\sigma;\theta_0, \theta_x,\theta_1,\theta_{\infty}, \hbox{tr}(M_0M_x),  \hbox{tr}(M_0M_1),
 \hbox{tr}(M_1M_x) \bigr)$ is derived for any $\sigma$ such that $\sigma \not \in (-\infty,0]
\cup [1,+\infty)$, but in \cite{Jimbo} also the case $\sigma=0$ is calculated. On the other 
hand the cases 
$\Re \sigma<0$ and $\Re \sigma\geq 1$ are not considered in \cite{Jimbo}, so we have included them in 
the Appendix. 

\vskip 0.3 cm 
 Conversely, we can start from $y(x;\sigma,a)$. Namely, 
for given $\theta_0$, $\theta_x$, $\theta_1$, 
$\theta_{\infty}$ and  $\sigma,a$,   we can compute tr($M_0M_x$), tr($M_1M_x$), tr($M_0M_1$) as functions of $a,\sigma,\theta_{\nu}$ ($\nu=0,x,1,\infty$) 
 by the formulae of the Appendix (in particular   $\hbox{tr}(M_0Mx)= 2 \cos(\pi \sigma)$). So the 
transcendent 
 $y(x;\sigma,a)$ associated to $\sigma,a$ 
  will coincide with the transcendent (\ref{PAALLAA}) associated to 
 the monodromy data.

\vskip 0.3 cm 
 We  claim the following:

\vskip 0.2 cm 
\noindent
{\it Let  $\sigma \not \in 
 (-\infty,0]\cup [1,+\infty)$ and $a\neq 0$.  Let $y(x)$ be a solution 
of $PVI$ in the generic case such that 
 $y(x)= a x^{1-\sigma}$(1+higher order terms)  as $x\to 0 $ in an open  
 domain contained in $D(\epsilon;\sigma)$.  
  Then, $y(x)$ coincides with $y(x;\sigma,a)$.} 

\vskip 0.2 cm

 The claim is proved in Proposition A1 in 
 the Appendix. As a consequence, a transcendents $y(x;\sigma,a)$  coincides with  a transcendents in elliptic representation with the same critical behavior.

 We fix the range of $\nu_2$ and $\sigma$ by  $ 
0\leq \Re \nu_2 < 2$, $0\leq \Re \sigma \leq 1
$. 
If $\Im \nu_2 \neq 0$, we consider a transcendent 
$$
 y(x) = \wp \Bigl(\nu_1\omega_1(x)+\nu_2\omega_2(x) +v(x;\nu_1,\nu_2);~ \omega_1(x),\omega_2(x)
 \Bigr)+{1+x\over 3}
$$
 on the domain  $ 
{\cal D}(r;\nu_1,\nu_2)  $. For $x\to 0$ it has behavior  
$$
y(x)= -{1\over 4} \left[{e^{i\pi \nu_1} \over 16^{\nu_2-1}}\right] ~
          x^{\nu_2} ~\left(1+~\hbox{higher orders}~\right).
$$
Thus it coincides  with 
\be
 y(x;\sigma, a(\sigma))=a(\sigma) x^{1-\sigma}  ~\left(1+~\hbox{higher orders}~\right)~~~ \hbox{ if } 0\leq \Re \nu_2 \leq 1 
\label{UNO11}
\ee
or
\be
  y(x;-\sigma, a(-\sigma))= a(-\sigma) x^{1+\sigma} 
 ~\left(1+~\hbox{higher orders}~\right) ~~~\hbox{ if } 1\leq \Re \nu_2 \leq 2
\label{DUE22}
\ee
where 
$$
\nu_2=1-\sigma,~~~e^{i\pi \nu_1} =-4a(\sigma) 16^{-\sigma}\equiv
 - 4a(1-\nu_2)~16^{\nu_2-1} ~~~\hbox{ in case (\ref{UNO11})} 
$$
$$
\nu_2=1+\sigma,~~~e^{i\pi \nu_1} = -4a(-\sigma) 16^{\sigma}\equiv
- 4a(1-\nu_2)~16^{\nu_2-1} ~~~\hbox{ in case (\ref{DUE22})} 
$$
We note that the corresponding domains ${\cal D}(r;\nu_1,\nu_2)$ and $D(\epsilon;\sigma)$ 
have non-empty 
intersection, therefore Proposition A1 -- namely the above claim -- holds. 
If $\nu_2$ is real, we consider a transcendent 
$$
y(x) =\wp \Bigl(\nu_1\omega_1(x)+\nu_2\omega_2(x) +v(x;\nu_1,\nu_2);~ \omega_1(x),\omega_2(x)
 \Bigr)+{1+x\over 3}
$$
\be
= -{1\over 4} \left[
{e^{i\pi \nu_1} \over 16^{\nu_2-1}} \right]x^{\nu_2} ~\left(1+~\hbox{higher orders}~\right), ~~\hbox{ if } 0<\nu_2<1~~~~ 
\label{UNOUNO1}
\ee
or 
$$
y(x)= \wp \Bigl(\nu_1\omega_1(x)+\nu_2\omega_2(x) +v(x;-\nu_1,2-\nu_2);~ \omega_1(x),\omega_2(x)
 \Bigr)+{1+x\over 3} 
$$
\be
=-{1\over 4} \left[
{e^{i\pi \nu_1} \over 16^{\nu_2-1}} \right]^{-1}
x^{2-\nu_2} ~\left(1+~\hbox{higher orders}~\right), ~~\hbox{ if } 1<\nu_2<2 ~
\label{DUEDUE2}
\ee
In both cases, they coincide with 
$$ 
   y(x;\sigma, a(\sigma))=a(\sigma) x^{1-\sigma}  ~\left(1+~\hbox{higher orders}~\right), ~~~ 0<\sigma<1
$$
where
$$
\nu_2=1-\sigma, ~~~~e^{i\pi \nu_1} = -4a(\sigma) 16^{-\sigma} 
\equiv -4 a(1-\nu_2) 16^{\nu_2-1} ~~~\hbox{ in case (\ref{UNOUNO1})} 
$$
$$ 
\nu_2=1+\sigma,~~~~e^{-i\pi \nu_1}= -4 a(\sigma) 16^{-\sigma}
\equiv -4 a(\nu_2-1) 16^{1-\nu_2} ~~~\hbox{ in case (\ref{DUEDUE2})}
$$
 The above identifications give $\nu_1$ and $\nu_2$ in terms of the 
monodromy data associated to $a,~\sigma$.

\vskip 0.3 cm 
 Conversely, we can start from $y(x;\sigma,a)$ ($0\leq 
\Re\sigma \leq 1$, $\sigma\neq 0,1$) and find its 
elliptic representation. This will be $
 y(x) = \wp \Bigl(\nu_1\omega_1+\nu_2\omega_2+v(x;\nu_1,\nu_2)
 \Bigr)+{1+x\over 3}$, $0\leq \Re \nu_2 \leq 1$ ($\nu_2\neq 0,1$),   with parameters  
$$\nu_2=1-\sigma,~~~~ e^{i\pi \nu_1} = 
-4 a(\sigma) ~16^{-\sigma}
$$
Again $(\epsilon;\sigma)$ and  $ 
{\cal D}(r;\nu_1,\nu_2)$ (or in  ${\cal D}_0(r)$ when
 $\Im \sigma=0$) have non empty intersection. 
 Equivalently, according to the Observation 2 in the Introduction,   $y(x;\sigma,a)$ coincides with  $
 y(x) = \wp \Bigl(\tilde{\nu}_1\omega_1+\tilde{\nu_2}\omega_2(x) +
v(x;-\tilde{\nu_1},2-\tilde{\nu_2})\Bigr)+{1+x\over 3}$,  $1\leq \Re 
\tilde{\nu}_2 \leq 2$ ($\tilde{\nu}_2\neq 1,2$), with parameters
$$
\tilde{\nu}_2=1+\sigma,~~~~
e^{-i\pi \tilde{\nu}_1} =- 4 a(\sigma) ~16^{-\sigma}
$$
 The domain is now 
 $ 
{\cal D}(r;-\tilde{\nu_1},2-\tilde{\nu_2})  
$ which has non-empty intersection with $D(\epsilon;\sigma)$.

\vskip 0.3 cm 
\noindent
{\bf Answer to Question 2)}. 
 We are ready to solve the connection problem.  In the Appendix we show, following \cite{Jimbo},  that 
 $$ 
y\bigl(x;~ \theta_0, \theta_x,\theta_1,\theta_{\infty}, \hbox{tr}(M_0M_x),  \hbox{tr}(M_0M_1),
 \hbox{tr}(M_1M_x)~ \bigr)=
$$
$$
  =  a^{(0)} x^{1-\sigma^{(0)}}~(1+~\hbox{higher orders}),~~~x\to 0
$$
   \be
= 1-a^{(1)} (1-x)^{1-\sigma^{(1)}}~(1+~\hbox{higher orders}),~~~x\to 1
\label{periduno}
\ee
\be
=
 a^{(\infty)} x^{\sigma^{(\infty)}}~(1+~\hbox{higher orders}),~~~x\to \infty
\label{periddue}
\ee
as $x$ converges to the critical points along regular paths contained in   domains $D(\epsilon;\sigma^{(0)}) $ defined in (\ref{dodoalmuseodivienna}),  $D(\epsilon;\sigma^{(1)}) $ defined in (\ref{dodoalmuseodivienna1}), $D(M;\sigma^{(\infty)})$ defined in 
 (\ref{dodoalmuseodiviennainf}) respectively.  
 The parameters are computed in terms of the monodromy data by 
 $$
2\cos(\pi \sigma^{(0)})= \hbox{tr}M_0M_x,~~
2\cos(\pi \sigma^{(1)})= \hbox{tr}M_1M_x,~~2\cos(\pi \sigma^{(\infty)})= \hbox{tr}M_0M_1
$$
$$
   0\leq \Re \sigma^{(i)} \leq 1,~~~\sigma^{(i)} \neq 0,1,~~~(i=0,1,\infty)
$$
 The coefficients $a^{(0)}$, $a^{(1)}$, $a^{(\infty)}$ are also  given in terms of 
the monodromy data by explicit 
formulae which are rather long to write here. The procedure  for their computation is explained in the Appendix. Note that here we  denote 
 $\sigma$, $a$ introduced  before  with  $\sigma^{(0)}$, $a^{(0)}$. 
By uniqueness of critical behavior (claim above -- or Proposition A1),    
 the elliptic representation at $x=1$ of  (\ref{periduno}) is   $y(x)= 
\wp ( \nu_1^{(1)} \omega_1^{(1)}+ \nu_2^{(1)} \omega_2^{(1)}+
v^{(1)}(x;\nu_1^{(1)},\nu_2^{(1)}))+{1+x\over 3}$, with 
\be
\nu^{(1)}_2=1-\sigma^{(1)},~~~e^{-i\pi \nu^{(1)}_1} = - 4a^{(1)}(1-\nu^{(1)}_2)~16^{\nu^{(1)}_2-1}
\label{ARIGATONI1}
\ee
At  $x=\infty$ 
 the representation of  (\ref{periddue}) is $y(x)= 
\wp ( \nu_1^{(\infty)} \omega_1^{(\infty)}+ \nu_2^{(\infty)} 
\omega_2^{(\infty)}+v^{(\infty)}(x;\nu_1^{(\infty)},\nu_2^{(\infty)}))+{1+x\over 3}$, with  
\be
\nu^{(\infty)}_2=1-\sigma^{(\infty)},~~~e^{i\pi \nu^{(\infty)}_1} = - 4a^{(\infty)}(1-\nu^{(\infty)}_2)~16^{\nu^{(\infty)}_2-1}
\label{ARIGATONI2}
\ee

Therefore, the solution of the connection problem 
 is a 
follows. 
           In the case  $\Im \nu_1^{(0)}\neq 0$, let us  start from an  elliptic representation  in ${\cal D}(r;\nu_1^{(0)},\nu_2^{(0)})$ close to $x=0$: 
\be
   y(x) = \wp ( \nu_1^{(0)} \omega_1^{(0)}+ \nu_2^{(0)} 
\omega_2^{(0)}+v^{(0)}(x;\nu_1^{(0)},\nu_2^{(0)}))+{1+x\over 3}.
\label{ritiro}
\ee  
 It must coincide  with $y(x;\sigma^{(0)},a^{(0)}(\sigma^{(0)}))$ or
 $y(x;-\sigma^{(0)},a^{(0)}(-\sigma^{(0)}))$, according to the value of 
$\Re \nu^{(0)}_2$, with parameters identified by the formulae which follow   (\ref{UNO11}) and (\ref{DUE22}). In the case 
 $\Im \nu_2^{(0)}= 0$,  
 let us  start from 
\be
   y(x) = \wp ( \nu_1^{(0)} \omega_1^{(0)}+ \nu_2^{(0)} 
\omega_2^{(0)}+v^{(0)}(x;\nu_1^{(0)},\nu_2^{(0)}))+{1+x\over 3}~~ ~~~0<\nu_2^{(0)}<1 
\label{ritiro1}
\ee
or 
\be
   y(x) = \wp ( \nu_1^{(0)} \omega_1^{(0)}+ \nu_2^{(0)} 
\omega_2^{(0)}+v^{(0)}(x;-\nu_1^{(0)},2-\nu_2^{(0)}))+{1+x\over 3}~~ ~~~1<\nu_2^{(0)}<2 
\label{ritiro2}
\ee
They must coincide  with $y(x;\sigma,a)$ with parameters identified by the formulae which follow  (\ref{UNOUNO1}) and (\ref{DUEDUE2})  respectively. 
 After these identifications, we compute monodromy data corresponding to  $\sigma^{(0)}$, $a^{(0)}$, and therefore we compute $\sigma^{(1)}$, 
 $a^{(1)}$ and $\sigma^{(\infty)}$, $a^{(\infty)}$ from the monodromy data.   
Finally, we write the 
elliptic representation of (\ref{periduno}), (\ref{periddue}) with parameters (\ref{ARIGATONI1}) and (\ref{ARIGATONI2}) respectively. These are the 
elliptic representations of (\ref{ritiro}) (or (\ref{ritiro1}), (\ref{ritiro2})) 
  at $x=1$ and $x=\infty$. This answers
 Question 2.

\vskip 0.3 cm
\noindent
{\bf Answer to Question 1)}.  
 As for Question 1, we recall that $
y\bigl(x;~ \theta_0, \theta_x,\theta_1,\theta_{\infty}, \hbox{tr}(M_0M_x),  \hbox{tr}(M_0M_1)
 \hbox{tr}(M_1M_x)~ \bigr)
$
 has representations $y(x; \pm\sigma+2n,a(\pm\sigma +
 2n))$, $n\in {\bf Z}$, on different domains. 
 Let us consider the family of  transcendents in elliptic representation (at $x=0$) 
$$ 
 y(x) =
\wp \bigl(\nu_1(N) \omega_1 +[\nu_2 +2N]~\omega_2+v(x;\nu_1(N)
,\nu_2+2N) \bigr)+{1+x\over 3} 
$$
$$
\equiv \wp \bigl(\bigl.\nu_1(N) \omega_1 +\nu_2 \omega_2+v(x;\nu_1(N)
,\nu_2+2N) \bigr)+{1+x\over 3}
$$ 
where $\nu_1=\nu_1(N)$ means that $\nu_1$ changes with $N$. 
 They are the elliptic representations
 of the same transcendent in different domains 
if and only if they coincide with 
$$
  y(x; \sigma-2N,a(\sigma-2N)) ~~~~\hbox{ if } 0 \leq  \Re \nu_2 \leq 1,
$$
or 
$$
  y(x; -\sigma-2N,a(-\sigma-2N)) ~~~~\hbox{ if } 1 \leq  \Re \nu_2 \leq 2,
$$
where  $y(x; \sigma-2N,a(\sigma-2N))$, or $y(x; -\sigma-2N,a(-\sigma-2N))$, correspond to the same monodromy data (i.e. they are  representations of  
$
y\bigl(x;~ \theta_0, \theta_x,\theta_1,\theta_{\infty}, \hbox{tr}(M_0M_x),  \hbox{tr}(M_0M_1)
 \hbox{tr}(M_1M_x)~ \bigr)
$). We can explicitly compute ( through the formulae of the Appendix) 
$\nu_1(N)=\hbox{function}(N,\nu_2,\nu_1(N=0))$ from 
$$ 
  e^{i\pi \nu_1(N)} = -4 a(1-\nu_2-2N) 16^{\nu_2+2N-1}
$$
  as we did in (\ref{UNO11}) and (\ref{DUE22}). 

\vskip 0.3 cm
\noindent
{\it Remark:} We have found the elliptic representations and the critical behaviors of 
$$
y\bigl(x;~ \theta_0, \theta_x,\theta_1,\theta_{\infty}, \hbox{tr}(M_0M_x),  \hbox{tr}(M_0M_1)
 \hbox{tr}(M_1M_x)~ \bigr)
$$
in the union of the domains 
$ {\cal D}(r_N;\nu_1(N),\nu_2+2N)$ ($\Im \nu_2\neq 0$). However, some regions in the $(\ln|x|,\Im \nu_2 \arg x)$-plane, for $x$ close to 0, are not included in the union. In these regions there may be movable poles (see figure \ref{figur1}).

\vskip 0.3 cm 
\noindent
In the generic case, we prove in the Appendix  the  
 one-to-one correspondence between transcendents and monodromy data. 
However, the condition $\sigma^{(i)} \neq 1$, namely $\nu_2^{(i)} = 0,2$,   
implies that we can not give the critical 
behaviors (and the elliptic representation) of 
$$
y\bigl(x;~ \theta_0, \theta_x,\theta_1,\theta_{\infty}, \hbox{tr}(M_0M_x),  \hbox{tr}(M_0M_1)
 \hbox{tr}(M_1M_x)~ \bigr)
$$ 
at $x=0$ for 
$\hbox{tr}(M_0M_x)=- 2$, at $x=1$ for 
$\hbox{tr}(M_1M_x)=-2$, at $x=\infty$ for 
$\hbox{tr}(M_0M_1)=- 2$. These cases have still to be studied.

 In the generic case we have also assumed $\sigma\not\in {\bf Z}$; nevertheless, in the case $\sigma^{(i)}=0$ (tr($M_iM_j)=2$) the critical behavior and the  solution of the connection problem  can be found in the paper \cite{Jimbo} by Jimbo.  We note that 
 the corresponding $\nu_2^{(i)}$ should be equal to 1, but  unfortunately the condition 
$\nu_2^{(i)}\neq 1$  which we had to impose to study the elliptic representation (except for special 
cases like $\beta=\gamma=1-2\delta=0$) did not allow us to know  the analytic properties and the critical behavior of the  elliptic representation in this case.  We expect that the properties of $u(x)$ are such to exactly produce the critical behavior found by Jimbo for $\sigma^{(i)}=0$, but we still have to cover this case.


\section{ APPENDIX}

 We give a brief account of the solution of the connection problem for the 
 generic PVI, following \cite{Jimbo}, with the extension of the values of  monodromy data for 
which the results apply. In \cite{Jimbo} the case $|\hbox{tr}(M_iM_j)|> 2$ is not 
considered, so we have to do it now.  This extension is necessary to identify the 
transcendent with the elliptic representation. Such a generalization is proved 
exactly as in  \cite{guz1}, to which we refer for a detailed analysis of the 
non-generic case $\beta=\gamma=1-2\delta=0$ (already reviewed in section \ref{so2qng}). 

\vskip 0.2 cm

 The Painlev\'e VI equation is the isomonodromic
 deformation equation of the fuchsian system (\ref{LAB}).  
 The  Fuchsian system is obtained from the following Riemann-Hilbert 
 problem.  We fix 
$$
\theta_0,~\theta_x,~\theta_1,~\theta_{\infty}~\not \in {\bf Z}
$$ and 
  monodromy matrices 
 $M_0$, $M_x$, $M_1$ corresponding to the loops in the basis of figure 
\ref{figure3}, with eigenvalues $\exp\{\pm i \pi 
 \theta_i \}$, $i=0,x,1$ 
respectively.  
 At infinity, the monodromy is  $M_1 M_x M_0$, and we require that this has 
eigenvalues  $\exp\{\pm i \pi 
 \theta_{\infty} \}$. Having assigned the monodromy matrices and their eigenvalues, there exist $2\times 2$ invertible matrices 
 $C_0$, $C_x$, $C_1$, $C_{\infty}$ such that 
\be
  M_i= C_i^{-1}~ \exp\left\{2\pi i \pmatrix{ {\theta_i\over 2} &  0
 \cr 
0  &    - {\theta_i\over 2} \cr} \right\}~ C_i,~~~i=0,x,1
\label{CONDIZIONE1}
\ee
$$
   M_1 M_x M_0 = C_{\infty}^{-1}~ \exp\left\{2\pi i \pmatrix{- {\theta_{\infty}\over 2} &  0
 \cr 
0  &    {\theta_{\infty}\over 2} \cr} \right\} C_{\infty}
$$

In order to 
construct the fuchsian system, we have to find a $2\times 2$ invertible 
matrix $Y(z;x)$ holomorphic in ${\bf P}^1\backslash \{0,x,1,\infty\}$ 
with the assigned monodromy at $z=0,1,x,\infty$. We also fix the {\it indices} of 
the Riemann-Hilbert problem, namely we require that $Y(z;x)$ be normalized 
as follows: 
\be
Y(z;x) = 
 \left\{ \matrix{  \left(I+O(\left({1\over z}\right) \right)~
z^{\pmatrix{- {\theta_{\infty}\over 2} &  0
 \cr 
0  &    {\theta_{\infty}\over 2} \cr}}~C_{\infty},~~~z\to \infty  
\cr 
\cr
G_i (I+O(z-i)) ~(z-i)^{\pmatrix{ {\theta_i\over 2} &  0
 \cr 
0  &    - {\theta_i\over 2} \cr} }~C_i,~~~z\to i=0,x,1
\cr
}
\right.
\label{SOLUZRH}
\ee
 Here  $G_i$ are invertible matrices. The fuchsian system is 
 obtained from 
\be
   A(z;x):={dY(z;x)\over dz}~Y(z;x)^{-1} ~~ \Longrightarrow~~
 A_i=G_i \pmatrix{ {\theta_i\over 2} &  0
 \cr 
\cr
0  &    - {\theta_i\over 2} \cr}~G_i^{-1}
\label{NONHAMAIFINE}
\ee
The fuchsian system satisfies (\ref{CONDIZIONI}). 
We remark that we have to assume $
  \theta_0,\theta_x,\theta_1,\theta_{\infty}\not \in {\bf Z}
$, 
otherwise the general form of $Y(z;x)$ would be different from (\ref{SOLUZRH}).
\vskip 0.2 cm 

 A $2\times 2$  Riemann-Hilbert always has solution \cite{AB}. In our case, 
the solution is unique, up to diagonal conjugation $A_i \mapsto D A_i D^{-1}$,
 where $D$ is any diagonal matrix. To prove this fact, let  
 $C_{\nu}$, $\tilde{C}_{\nu}$, $\nu=0,x,1,\infty$, be such that 
$$ 
   C_\nu^{-1}~ \exp\left\{2\pi i \pmatrix{ {\theta_\nu\over 2} &  0
 \cr 
0  &    - {\theta_\nu\over 2} \cr} \right\}~ C_\nu = 
\tilde{C}_\nu^{-1}~ \exp\left\{ 2\pi i \pmatrix{ {\theta_\nu\over 2} &  0
 \cr 
0  &    - {\theta_\nu\over 2} \cr} \right\}~ \tilde{C}_\nu
$$
 If follows that $C_{\nu} \tilde{C}_{\nu}^{-1}$ is any diagonal matrix. We denote $Y$ and $\tilde{Y}$ the solutions (\ref{SOLUZRH}) 
with $C_\nu$ and $\tilde{C}_\nu$ 
respectively. Since  $C_{\nu} \tilde{C}_{\nu}^{-1}$ is diagonal we conclude
 that $Y~\tilde{Y}^{-1}$ is holomorphic at $z=0,1,x,\infty$. Therefore, being 
holomorphic on ${\bf P}^1$,  it 
 is a constant diagonal matrix 
 $D:= C_{\infty}\tilde{C}_{\infty}^{-1}$. This proves uniqueness up to 
 diagonal conjugation. We remark that if some $\theta_{\nu}$ ($\nu=0,x,1,\infty$) is integer, the uniqueness may fail. The cases when this happens  if  $\beta=\gamma=1-2\delta=0$ are discussed in 
 \cite{M}. 

\vskip 0.2 cm 
 Once the Riemann-Hilbert problem is solved, we compute $y(x)$ from 
$A(z;x)_{12}=0$. It is clear that $y(x)$ 
 depends on the monodromy data $M_0$, $M_x$, $M_1$, $\theta_{\nu}$ 
($\nu=0,1,x,\infty$). 
Note that for any invertible matrix $C$ we obtain the  fuchsian system  
(\ref{NONHAMAIFINE})  
from any $Y(z;u)C$. The solution  $Y(z;u)C$ corresponds to  
   the monodromy  matrices $
C^{-1} M_i C$. This implies that $y(x)$ depends on the invariants 
of $M_i$ ($i=0,x,1$) with respect to  the conjugation $ M_i \mapsto 
C^{-1} M_i C$. They are traces and determinants of the products on the $M_i$'s. 
The determinants are 1, the traces of $M_i$ are specified by the eigenvalues
 and the traces of the products $M_{i_1} M_{i_2}...M_{i_n}$ for $n>2$ are 
functions of the traces of $M_{i_1} M_{i_2}$. Hence: 
$$
y(x)= y\bigl(x;~ \theta_0, \theta_x,\theta_1,\theta_{\infty}, \hbox{tr}(M_0M_x),  \hbox{tr}(M_0M_1),
 \hbox{tr}(M_1M_x)~ \bigr)
$$

\vskip 0.2 cm 
 PVI is equivalent to the  Schlesinger equations which ensure  the isomonodromicity of the fuchsian system \cite{JMU}: 
\be 
\left.\matrix{
               {d A_0\over dx}= {[A_x,A_0]\over x} \cr
\cr               
               {dA_1\over dx}={[A_1,A_x]\over 1-x} \cr
\cr  
               {d A_{x} \over d x} = {[A_x,A_0]\over x}+{[A_1,A_x]\over 1-x}\cr
}\right.
\label{sch}
\ee
 The system (\ref{sch}) is a particular case of  
\be 
\left. \matrix{
{d A_{\mu}\over dx}=\sum_{\nu=1}^{n_2} [A_{\mu},B_{\nu}]~f_{\mu \nu}(x)
\cr\cr
{dB_{\nu}\over dx}= -{1\over x} \sum_{ \nu^{\prime}=1 }^{ n_2 }[B_{\nu},
B_{\nu^{\prime}}]+
\sum_{\mu=1}^{n_1}[B_{\nu},A_{\mu}]~g_{\mu \nu}(x) +     
\sum_{\nu^{\prime}=1}^{n_2} [B_{\nu},B_{\nu^{\prime}}] ~h_{\nu
\nu^{\prime} }(x) 
}
\right.
\label{sch1}
\ee
where the functions $f_{\mu\nu}$, $g_{\mu \nu}$, $h_{\mu \nu}$ are
meromorphic with poles at $x=1,\infty$ and $$
\sum_{\nu} B_{\nu} +
\sum_{\mu} A_{\mu}=-A_{\infty},~~~~ A_{\infty}:= 
\pmatrix{ {\theta_{\infty} \over 2} & 0 \cr 0 & -{\theta_{\infty} \over 2} }
$$
 System (\ref{sch}) is obtained
for $f_{\mu \nu}=g_{\mu \nu}=b_{\nu}/(a_{\mu}-xb_{\nu})$, $h_{\mu
\nu}=0$, $n_1=1$, $n_2=2$, $a_1=b_2=1$, $b_1=0$ 
and $B_1=A_0$, $B_2=A_x$, $A_1=A_1$.

\vskip 0.2 cm 

 Let $\sigma$ be a complex number such that 
$$ 
  \sigma \not\in (-\infty,0)\cup[1,+\infty)
$$
and let us consider the domain   (\ref{dodoalmuseodivienna}) with  additional parameters 
$\vartheta_1$, $\vartheta_2 \in {\bf R}$, $0<\tilde{\sigma} 
<1$.
 In \cite{guz1} the reader can find the proof of the following lemma, which is 
a generalization in $ D(\epsilon;\sigma)$ 
of \cite{SMJ}, page 262

\vskip 0.3 cm
\noindent
 {\bf Lemma A1: } {\it  
                  Consider  matrices $B_{\nu}^0$ ($\nu=1,..,n_2$), 
$A_{\mu}^0$ ($\mu=1,..,n_1$) and   $\Lambda$, independent of $x$ and such
                  that 
$$
   \sum_{\nu} B_{\nu}^0 + \sum_{\mu} A_{\mu}^0 = -A_{\infty}
$$
$$
   \sum_{\nu} B_{\nu}^0=\Lambda, ~~~~\hbox{ eigenvalues}(\Lambda)=
   {\sigma \over 2},~-{\sigma\over 2},~~~
\sigma \not\in(-\infty,0)\cup [1,+\infty). 
$$
  Suppose that  $f_{\mu \nu}$, $g_{\mu \nu}$, $h_{\mu \nu}$ are holomorphic
if  $|x|<\epsilon^{\prime}$, for some small $\epsilon^{\prime}<1 $.

For any $0<\tilde{\sigma}<1$ and  $\vartheta_1,\vartheta_2$ real there exists a
sufficiently small $0<\epsilon<\epsilon^{\prime}$ such that 
the system  (\ref{sch1}) has holomorphic 
solutions 
 $A_{\mu}({x})$, $B_{\nu}({x})$   in
$D(\epsilon;\sigma)$   satisfying:
$$
   || A_{\mu}({x})-A_{\mu}^0||\leq C~ |x|^{1-\sigma_1}$$
$$
           || {x}^{-\Lambda} B_{\nu}({x})~ {x}^{\Lambda} 
-B_{\nu}^0||\leq C~ |x|^{1-\sigma_1}
$$          
            Here $C$ is a positive  constant and $\tilde{\sigma}<\sigma_1<1$
}

\vskip 0.3 cm

\noindent
{\bf Lemma A2:} {\it Let $\theta_{\infty}\neq 0$. 
Let $r$, $s$ be two complex numbers not equal to zero. Let 
$\hat{G}_0$ be such that 
$$ 
     \hat{G}_0^{-1} \Lambda \hat{G}_0 = \hbox{\rm diag}\left({\sigma\over 2}, 
-{\sigma \over 2} \right),~~~~\sigma\neq 0
$$
The general solution of 
$$A_0^0+A_x^0+A_1^0=\hbox{\rm diag}\left(-{\theta_{\infty}\over 2},
{\theta_{\infty}\over 2}\right),~~~ A_0^0+A_x^0=\Lambda, ~~~\hbox{ eigenvalues } A_i 
=\pm {\theta_i \over 2}~~~(i=0,1,x)
$$
is 
$$
   A_0^0=  \hat{G}_0 \pmatrix{{\alpha(\beta+1-\gamma)\over \alpha-\beta} & 
                     s \cr \cr
{\alpha\beta(\beta+1-\gamma)(\gamma-1-\alpha)\over (\alpha-\beta)^2 }~{1\over
 s}
&
{\beta(\gamma-\alpha-1)\over\alpha-\beta} \cr}
\hat{G}_0^{-1} +{\theta_0\over 2},
$$
\vskip 0.2 cm 
$$
A_x^0= \hat{G}_0\pmatrix{ {\alpha(\gamma-\alpha-1)\over \alpha-\beta} & 
                - s  \cr\cr
-(A_0^0)_{21}   &  {\beta(\beta+1-\gamma)\over \alpha-\beta} \cr
}\hat{G}_0^{-1} +{\theta_x\over 2}
$$
where 
$$
\alpha:={-\sigma+\theta_0+\theta_x \over 2},~~~\beta:={\sigma+\theta_0+\theta_x \over 2},~~~\gamma:=1+\theta_0
$$
and 
$$
   \Lambda=  \pmatrix{{\alpha(\beta+1-\gamma)\over \alpha-\beta} & 
                     r \cr \cr
{\alpha\beta(\beta+1-\gamma)(\gamma-1-\alpha)\over (\alpha-\beta)^2 }~{1\over
 r}
&
{\beta(\gamma-\alpha-1)\over\alpha-\beta} \cr}
 +{\sigma\over 2},~~~
A_1^0= \pmatrix{ {\alpha(\gamma-\alpha-1)\over \alpha-\beta} & 
                - r  \cr \cr
-(\Lambda)_{21}   &  {\beta(\beta+1-\gamma)\over \alpha-\beta} \cr
} +{\theta_1\over 2}
$$
\vskip 0.2 cm 
$$
\hat{G}_0 = \pmatrix{ 1 & 1 \cr \cr
{\alpha(1+\beta-\gamma)\over \beta-\alpha} {1\over r} & \left\{ 
{\alpha(1+\beta-\gamma)\over \beta-\alpha}+1-\gamma 
\right\}{1\over r} \cr}
$$
where  
$$
\alpha:= {\theta_{\infty}+\theta_1+\sigma\over 2},~~~\beta := 
 {-\theta_{\infty}+\theta_1+\sigma\over 2},~~~\gamma:=1+\sigma
$$
}(note that  $\alpha,\beta,\gamma$ are not 
 the coefficients of PVI. We apologize for using the same symbols.)

\vskip 0.3 cm 
 Setting $A(z;x)_{12}=0$ we obtain from Lemma A1 and Lemma A2: 
$$ 
  y(x)= -x~{(A_0)_{12}\over (A_1)_{12}} ~(1+O(|x|^{1-\tilde{\sigma}}))
$$
\be 
=   {1\over 16 \sigma^3}\left[ \theta_0+\theta_x+\sigma \right]
\left[ -\theta_0+\theta_x+\sigma \right]
\left[ \theta_0+\theta_x-\sigma \right]
\left[ \theta_0-\theta_x+\sigma \right]~{1\over s} ~ x^{1-\sigma} (1 +
O(|x|^{\Delta}))
\label{F}
\ee
where $\Delta$ is a small positive  number and $x\to 0$ in 
$D(\epsilon;\sigma)$.  
 Namely, we have the following theorem (for a detailed proof see \cite{guz1}):

\vskip 0.3 cm 
\noindent
{\bf Theorem A1: } {\it Let $\theta_{\infty} \neq 0$.   
For any complex 
$\sigma \not \in (-\infty,0]\cup [1,+\infty)$, for
any complex 
 $a\neq 0$, for any $\vartheta_1,\vartheta_2  \in {\bf R}$ and for any
$0<\tilde{\sigma }<1$, there exists a sufficiently small positive $\epsilon$ 
and a transcendent 
$y(x;\sigma,a)$ with behavior  
\be
    y(x;\sigma,a)=a x^{1-\sigma} \left(
1+O(|x|^{\Delta})
\right),~~~~0<\Delta<1,
\label{asy0}
\ee
as $x\to 0$ along a regular path in  $D(\epsilon;\sigma)$.  

The above local behavior is valid along any regular path, with  the  exception, which occurs if 
 $\Im \sigma\neq 0$, of the paths 
 $\Im \sigma \arg(x) = \Re \sigma \log|x| + b
$, where $b$ is a constant such that the path is 
 contained in $D(\epsilon;\sigma)$: in this case the
behavior is 
\be
   y({x};\sigma,a)= \sin^2\left({i\sigma \over 2} \ln x -{i \over 2} \ln (4a) 
-{\pi \over 2} \right)~x~ (1+O(|x|^{\Delta}))
\label{cecilia}
\ee
}
\vskip 0.3 cm 
Note that we have excluded $\sigma=0$ from Lemma A2 and therefore from Theorem A1. Such a case is however computed in \cite{Jimbo} and so  Theorem A1  holds also 
for $\sigma=0$. On the other hand, in \cite{Jimbo} only $0\leq \Re \sigma<1$ is considered, 
thus Theorem A1 is a generalization of the result of \cite{Jimbo}.  

\vskip 0.2 cm
  
We assume in the following that $ 
  \theta_{\nu} \not \in {\bf Z}$, $\nu=0,x,1,\infty
$. We choose a solution of  the fuchsian system  
\be
  {dY\over dz} = \left[{A_0({x})\over z} +{A_x({x})\over z-x}+
{A_1({x})\over z-1}
   \right] ~Y
\label{stofucs}
\ee
normalized as follows 
\be
Y(z,{x})= \left(I+O\left({1\over z}\right)\right) 
~z^{-A_{\infty}} ,~~~z\to \infty,~~~~A_{\infty}=\hbox{diag} \left(
{\theta_{\infty}\over 2}, -{\theta_{\infty}\over 2}\right).~~
\label{solutionzx}
\ee 
As a corollary of Lemma A1, 
the limits 
$$
   \hat{Y}(z):=\lim_{{x}\to 0} Y(z,{x}),~~~~
 \tilde{Y}(z):=\lim_{{x}\to 0} ~{x}^{-\Lambda} 
Y({x}z,{x})
$$
exist when  ${x}\to 0$ in $D(\epsilon;\sigma)$. They satisfy
\be
{d \hat{Y}\over dz} = \left[{A_1^0 \over z-1}+{\Lambda \over z}
   \right] ~\hat{Y}
\label{systemhat}
\ee
\be
  {d\tilde{Y}\over dz} = \left[{A_0^0\over z}+{A_x^0\over z-1}   \right]
    ~\tilde{Y}
\label{systemtilde}
\ee
For the system (\ref{systemhat}) we choose a fundamental
matrix solution   normalized as follows
\be
   \hat{Y}_N(z)= \left(I+O\left({1\over z}\right)\right) ~z^{-A_{\infty}},
~~~~z\to \infty
\label{solutionhat}
\ee
$$
  =   (I+O(z))~z^{\Lambda}~\hat{C}_0,~~~~z\to 0 
$$
$$
   = \hat{G}_1 (I+O(z-1))~(z-1)^{\hbox{diag}\left({\theta_1\over 2},-{\theta_1\over 2} \right)} ~\hat{C}_1,~~~~z \to 1 
$$
Where $\hat{G}_1^{-1} A_1^0 \hat{G}_1=\hbox{diag}\left({\theta_1\over 2},-{\theta_1\over 2} \right)$.  
$\hat{C}_0$,  $\hat{C}_1$ are { connection matrices}. For  
(\ref{systemtilde}) we choose a fundamental matrix solution normalized
as follows
 \be
   \tilde{Y}_N(z)= \left(I+O\left({1\over z}\right)\right) ~z^{\Lambda},
~~~~z\to \infty
\label{solutiontilde}
\ee
$$
  =  \tilde{G}_0 (I+O(z))~
z^{\hbox{diag}\left({\theta_0\over 2},-{\theta_0\over 2} \right)}~\tilde{C}_0,~~~~z\to 0 
$$
$$
   = \tilde{G}_1 (I+O(z-1))~
(z-1)^{\hbox{diag}\left({\theta_x\over 2},-{\theta_x\over 2} \right)} ~\tilde{C}_1,~~~~z \to 1. 
$$
Here $\tilde{G}_i^{-1}A_i^0 \tilde{G}_i =\hbox{diag}\left({\theta_i\over 2},-{\theta_i\over 2} \right)$ and  $\tilde{C}_0$,  $\tilde{C}_1$ are { connection matrices}.
 As it is proved in \cite{Jimbo} and in \cite{guz1},   we have:    
$$
    \hat{Y}(z) = \hat{Y}_N(z)
$$ 
\be     
 \tilde{Y}(z)=
 \tilde{Y}_N(z)~\hat{C}_0
 \label{lim}
\ee
and, as a  consequence of isomonodromicity, we have  
\be
  M_1= \hat{C}_1^{-1} e^{2\pi i 
~\hbox{diag}\left({\theta_1\over 2},-{\theta_1\over 2} \right)
 } \hat{C}_1
\label{*}
\ee
\be
  M_0=\hat{C}_0^{-1}\tilde{C}_0^{-1} e^{2\pi i
~\hbox{diag}\left({\theta_0\over 2},-{\theta_0\over 2} \right)
 } \tilde{C}_0 \hat{C}_0
\label{**}
\ee
\be
  M_x=\hat{C}_0^{-1}\tilde{C}_1^{-1} e^{2\pi i
~\hbox{diag}\left({\theta_x\over 2},-{\theta_x\over 2} \right)
 } \tilde{C}_1 \hat{C}_0
\label{***}
\ee

\vskip 0.2 cm 
 The connection matrices $\hat{C}_0$, $\tilde{C}_0$, $\tilde{C}_1$ can be 
computed explicitly because the  $2\times 2$ fuchsian systems  
(\ref{systemhat}) (\ref{systemtilde}) 
can be reduced to the
hyper-geometric equation.

\vskip 0.3 cm
\noindent 
{\bf Lemma A3: } {\it 
The Gauss hyper-geometric equation 
\be
   z(1-z)~ {d^2 y \over dz^2} +[\gamma_0-z(\alpha_0+\beta_0+1)]~{dy\over dz}
   -\alpha_0 \beta_0 ~y=0
\label{hyper1}
\ee
 is equivalent to the system 
\be
{d\Psi\over dz}= \left[{1\over z}\pmatrix{0 & 0 \cr
                                       -\alpha_0 \beta_0 & -\gamma_0 \cr}
                                +{1\over z-1}\pmatrix{0&1\cr
                                                     0 & \gamma_0
                                       -\alpha_0-\beta_0
                                                     \cr} \right]~\Psi
\label{hyper2}
\ee
where $\Psi=\pmatrix{y \cr (z-1) {dy\over dz}\cr}$. 
}

\vskip 0.3 cm 
\noindent
{\bf  Lemma A4: } { \it Let $B_0$ and $B_1$ be matrices  of eigenvalues $0,
1-\gamma$,  and  $0, \gamma-\alpha-\beta-1$ respectively, such
that }
$$
  B_0+B_1= \hbox{ diag}(-\alpha, -\beta),~~~\alpha\neq \beta$$
{\it Then
$$
   B_0=   \pmatrix{{\alpha(\beta+1-\gamma)\over \alpha-\beta} & 
                     r_1 \cr 
{\alpha\beta(\beta+1-\gamma)(\gamma-1-\alpha)\over (\alpha-\beta)^2 }~{1\over r_1}
&
{\beta(\gamma-\alpha-1)\over\alpha-\beta} \cr}
$$
\vskip 0.2 cm
$$
B_1= \pmatrix{ {\alpha(\gamma-\alpha-1)\over \alpha-\beta} & 
                - r_1  \cr
-(B_0)_{21}   &  {\beta(\beta+1-\gamma)\over \alpha-\beta} \cr
}
$$
for any  $r_1\neq 0$. 
}
\vskip 0.3 cm 
\noindent
We leave the proof as an exercise (note that 
 $\alpha,\beta,\gamma$ are not the coefficients of PVI. We apologize   
for using the same symbols).       
The following lemma connects Lemmas A3 and A4:
\vskip 0.3 cm

\noindent
{\bf Lemma A5: } {\it The system (\ref{hyper2}) with 
$$\alpha_0=\alpha,~~~\beta_0=\beta+1,~~~\gamma_0=\gamma,~~~\alpha\neq \beta$$ 
is gauge-equivalent to
the system 
\be
 {d X\over dz}= \left[{B_0\over z}+{B_1\over z-1}\right] ~X
\label{hyper3}
\ee
where $B_0$, $B_1$ are given in Lemma A4. This means that there exists
a matrix 
$$
   G(z)= \pmatrix{ 1  &   0  \cr  
{\alpha((\alpha-\beta)z+\beta+1-\gamma)\over \beta-\alpha}  ~{1\over r_1} & 
   - {z\over r_1} \cr
}
$$
such that $X(z)=G(z) ~\Psi(z)$. 
 It follows that (\ref{hyper3}) and the corresponding hyper-geometric
equation (\ref{hyper1}) have the same fuchsian singularities
0,1,$\infty$ and the same monodromy group.
}
\vskip 0.3 cm 
\noindent
{\it Proof:}  By direct computation. $\Box$

\vskip 0.3 cm
\noindent
 Note that the form of $G(z)$ ensures that if
$y_1$, $y_2$ are independent solutions of the hyper-geometric
equation, then a fundamental matrix of (\ref{hyper3}) may be chosen to
be 
$X(z)=
\pmatrix{ y_1(z) & y_2(z) \cr
           *  & *\cr}
$.

\vskip 0.3 cm
Now we compute the monodromy matrices for the systems
(\ref{systemhat}), (\ref{systemtilde}) by reduction to an
hyper-geometric equation. We assume 
$$\sigma \not \in {\bf Z}$$
Let us start with (\ref{systemhat}). With the
gauge  $ 
  Y^{(1)}(z):= z^{-{\sigma\over 2}} (z-1)^{-{\theta_1\over 2}} \hat{Y}_N(z)
$ 
  we transform (\ref{systemhat}) into 
\be 
      {d Y^{(1)}\over dz}=\left[ {A_1^0-{\theta_1\over 2}\over
      z-1}+{\Lambda-{\sigma\over2 }I\over z} \right]~Y^{(1)}
\label{novosis}
\ee
We identify the matrices $B_0$, $B_1$ of Lemma A4 with    
$\Lambda-{\sigma\over2 }$ and  $A_1^0-{\theta_1\over 2}$. Therefore 
$$ 
 \alpha= {\theta_{\infty}+\theta_1+\sigma\over 2},~~~\beta  =
 {-\theta_{\infty}+\theta_1+\sigma\over 2}, ~~~\gamma=1+\sigma
$$ 
(note that $\alpha-\beta =\theta_{\infty}\neq 0$. 
 Also note that Lemma A2  for $\Lambda$ and $A^0_1$ 
follows from Lemma A4 applied to (\ref{novosis})).  
The hyper-geometric equation connected to the present system through 
 Lemma A3  has coefficients  
$$ 
 \left\{ \matrix{                                         
\alpha_0={\theta_{\infty}+\theta_1+\sigma\over 2} \cr
\beta_0= 1+ {-\theta_{\infty}+\theta_1+\sigma\over 2}\cr
      \gamma_0=\sigma+1 \cr
}\right.
$$
Therefore, according to the standard theory of hyper-geometric equations, 
 the generic case occurs when  
 $$ 
\sigma, ~\theta_1,~ \theta_{\infty} \not \in {\bf Z}
$$
 Let $F(a,b,c,;z)$ be the hyper-geometric function. For $\sigma\not \in {\bf Z}$ we have two independent solutions of the hyper-geometric equation
$$ 
  y_1^{(0)}(x)= F(\alpha_0,\beta_0,\gamma_0;z)
$$
$$ 
y_2^{(0)}(x)= z^{1-\gamma_0} F(\alpha_0-\gamma_0+1,\beta_0-\gamma_0+1,2-\gamma_0;z)
$$
 For $\theta_1\not \in {\bf Z}$ we have also solutions 
$$ 
y_1^{(1)}(z) = F(\alpha_0,\beta_0,\alpha_0+\beta_0+1-\gamma_0;1-z) 
$$
$$
 y_2^{(1)}(z) = (1-z)^{\gamma_0-\alpha_0-\beta_0}F(\gamma_0-\beta_0, \gamma_0-\alpha_0,1+\gamma_0-\alpha_0-\beta_0;1-z)
$$
They are connected by 
$$ 
[y_1^{(0)},y_2^{(0)}]=[y_1^{(1)},y_2^{(1)}] C_{01}
$$
\vskip 0.2 cm 
$$ 
 C_{01}:= \pmatrix{ { \Gamma(\gamma_0-\alpha_0-\beta_0)\Gamma(\gamma_0)
     \over
\Gamma(\gamma_0-\alpha_0) \Gamma(\gamma_0-\beta_0) }  & 
{\Gamma(\gamma_0-\alpha_0-\beta_0)\Gamma(2-\gamma_0) \over 
 \Gamma(1-\alpha_0) \Gamma(1-\beta_0)} 
\cr 
\cr
{\Gamma(\alpha_0+\beta_0-\gamma_0)\Gamma(\gamma_0)\over \Gamma(\alpha_0)
\Gamma(\beta_0) }
&
{\Gamma(\alpha_0+\beta_0-\gamma_0) \Gamma(2-\gamma_0)\over 
\Gamma(\alpha_0+1-\gamma_0)\Gamma(\beta_0+1-\gamma_0)}
\cr
}
$$
Note that 
$$ 
Y^{(1)}(z)= \hat{G}_0 [I+O(z)]~z^{\hbox{diag}\left(0,
-\sigma\right)} ~\hat{G}_0^{-1} \hat{C}_0,~~~z\to 0
$$
$$
= \hat{G}_1 [I+O(z-1)]~(z-1)^{\hbox{diag}\left(0,
-\theta_1\right)} ~ \hat{C}_1,~~~z\to 1
$$
 From the behaviors of the hyper-geometric functions  $y_i^{(0)}$ 
for $z\to 0 $ and  $y_i^{(1)}$ for 
$z\to 1$ and for a suitable choice of  $\hat{G}_1$ we obtain  
$$ 
Y^{(1)}=\pmatrix{ y_1^{(0)} & y_2^{(0)} \cr \cr 
                     *    &  *      \cr } ~\hat{G}_0^{-1} \hat{C}_0
=\pmatrix{ y_1^{(1)} & y_2^{(1)} \cr \cr 
                     *    &  *      \cr }~ \hat{C}_1
$$
Namely
$$ 
  \hat{C}_1 (\hat{G}_0^{-1}\hat{C}_0)^{-1} = C_{01}
$$
 This will be enough for our purposes (note that we did not use 
the hypothesis $\theta_{\infty} \not \in {\bf Z}$).

\vskip 0.3 cm 
We turn to the system (\ref{systemtilde}).  With the gauge $ 
Y^{(2)}(z):=z^{-{\theta_0\over 2}} (z-1)^{-{\theta_x\over 2}}
~\hat{G}_0^{-1}\tilde{Y}_N(z)\hat{G}_0 
$ we have 
$$
  {d Y^{(2)} \over d z} = \left[ {B_0 \over
  z}+{B_1\over z-1} \right] Y^{(2)}
$$
where 
  $$
 B_0=\hat{G}_0^{-1} A_0^0 \hat{G}_0 -{\theta_0\over 2}  , ~~~~
B_1
=\hat{G}_0^{-1} A_x^0 \hat{G}_0 -{\theta_x\over 2}
  $$
$B_0$ and $B_1$ are as in Lemma A4, with 
$$ 
\alpha= {-\sigma+\theta_0+\theta_x\over 2},
~~~\beta={\sigma+\theta_0+\theta_x\over 2},
~~~\gamma= 1+\theta_0
$$
(Lemma A2  for $A^0_x$ and $A^0_0$ 
follows from Lemma A4 applied to the present case). 
The corresponding hyper-geometric equation has coefficients
$$
\left\{\matrix{\alpha_0={-\sigma+\theta_0+\theta_x\over 2} \cr
                \beta_0=1+{\sigma+\theta_0+\theta_x\over 2} \cr
\gamma_0=1+\theta_0
}\right.  
$$
Therefore, according to the standard theory of hyper-geometric equations, 
 the generic case occurs for 
$$ 
\theta_0,~\theta_x,~\sigma \not \in {\bf Z}
$$
and we can choose three couples of independent solutions 
$$ \left\{\matrix{
  y_1^{(0)}= F(\alpha_0,\beta_0,\gamma_0;z)
\cr
  y_2^{(0)}= z^{1-\gamma_0}F(1+\alpha_0-\gamma_0,1+\beta_0-\gamma_0,2-\gamma_0;z) }\right.
$$
$$ \left\{\matrix{
  y_1^{(1)}= F(\alpha_0,\beta_0,\alpha_0+\beta_0+1-\gamma_0;1-z)
\cr
  y_2^{(1)}= (1-z)^{\gamma_0-\alpha_0-\beta_0}
F(\gamma_0-\alpha_0,\gamma_0-\beta_0,\gamma_0+1-\alpha_0-\beta_0;1-z) }\right. 
$$
$$ \left\{\matrix{
  y_1^{(\infty)}=z^{-\alpha_0}
 F(\alpha_0,\alpha_0-\gamma_0+1,\alpha_0+1-\beta_0;{1\over z})
\cr
  y_2^{(\infty)}=z^{-\beta_0}
F(\beta_0,\beta_0+1-\gamma_0,\beta_0+1-\alpha_0;{1\over z} }\right. 
$$
They are connected by 
$$ 
  [y_1^{(\infty)},y_2^{(\infty)}]=[y_1^{(0)},y_2^{(0)}] C_{\infty0}^{\prime}
$$
\vskip 0.2 cm 
$$ 
  C_{\infty 0}^{\prime}= \pmatrix{ e^{-i\pi \alpha_0}{\Gamma(1+\alpha_0-\beta_0)
\Gamma(1-\gamma_0) \over \Gamma(1-\beta_0)\Gamma(1+\alpha_0-\gamma_0)} 
&
 e^{-i\pi \beta_0} {\Gamma(1+\beta_0-\alpha_0) \Gamma(1-\gamma_0) \over 
\Gamma(1-\alpha_0)\Gamma(1+\beta_0-\gamma_0)} 
\cr 
\cr
e^{i\pi(\gamma_0-1-\alpha_0)} {\Gamma( 1+\alpha_0-\beta_0)\Gamma(\gamma_0-1) 
\over \Gamma(\alpha_0) \Gamma(\gamma_0-\beta_0) } &
e^{i\pi (\gamma_0-1-\beta_0)} {\Gamma(1+\beta_0-\alpha_0) \Gamma(\gamma_0-1) 
\over \Gamma(\beta_0) \Gamma(\gamma_0-\alpha_0)} 
\cr
}
$$
\vskip 0.2 cm 
$$
  [y_1^{(\infty)}, y_2^{(\infty)}]= [y_1^{(1)}, y_2^{(1)}]C_{\infty 1}^{\prime}
$$
\vskip 0.2 cm 
$$
C_{\infty 1}^{\prime}
=
\pmatrix{
{\Gamma(\gamma_0-\alpha_0-\beta_0) \Gamma(1+\alpha_0-\beta_0) \over 
\Gamma(1-\beta_0) \Gamma(\gamma_0-\beta_0) } & 
{\Gamma(\gamma_0-\alpha_0-\beta_0) \Gamma(1+\beta_0-\alpha_0) \over 
\Gamma(1-\alpha_0) \Gamma(\gamma_0-\alpha_0) } 
\cr
\cr
e^{i\pi ( \gamma_0-\alpha_0-\beta_0)} {\Gamma(\alpha_0+\beta_0-\gamma_0) 
\Gamma(1+\alpha_0-\beta_0) \over \Gamma(1+\alpha_0-\gamma_0) \Gamma(\alpha_0)}
&
e^{i\pi(\gamma_0-\alpha_0-\beta_0)} {\Gamma(\alpha_0+\beta_0-\gamma_0) 
\Gamma(1+\beta_0-\alpha_0) \over \Gamma(1+\beta_0-\gamma_0) \Gamma(\beta_0)}
\cr
}
$$
 Considering the asymptotic behavior of the hyper-geometric equations and of $Y^{(2)}(z)$ at $z=0,1,\infty$ we obtain
$$ 
Y^{(2)}= \pmatrix{ y_1^{(\infty)} & {s\over 1+\sigma} y_2^{(\infty)} 
\cr
\cr
* & * \cr
}
=  \pmatrix{ y_1^{(0)} &  y_2^{(0)} 
\cr
\cr
* & * \cr
}\tilde{C}_0\hat{G}_0
=  \pmatrix{ y_1^{(1)} &  y_2^{(1)} 
\cr
\cr
* & * \cr
}\tilde{C}_1\hat{G}_0
$$
Therefore
$$ 
   \tilde{C}_0 \hat{G}_0 \pmatrix{ 1 & 0 \cr \cr 0 & {1+\sigma\over s} 
\cr 
}
= C_{\infty0}^{\prime},~~~
\tilde{C}_1 \hat{G}_0 \pmatrix{ 1 & 0 \cr \cr 0 & {1+\sigma \over s} \cr}
= C_{\infty 1}^{\prime}
$$

\vskip 0.3 cm 

 We introduce the following matrices: 
\be
 m_0:=  \left[
\left\{C_{\infty0}^{\prime} \pmatrix{ 1 & 0 \cr \cr 0 & {s\over 1+\sigma}}  
 \right\}^{-1} ~ e^{2\pi i~\hbox{diag}\left(
{\theta_0\over 2},-{\theta_0\over 2}
\right)}~ 
\left\{C_{\infty0}^{\prime} \pmatrix{ 1 & 0 \cr \cr 0 & {s\over 1+\sigma}}  
 \right\}
\right]
\label{m000}
\ee
\be
m_x:= \left[
\left\{C_{\infty1}^{\prime} \pmatrix{ 1 & 0 \cr \cr 0 & {s\over 1+\sigma} } 
 \right\}^{-1} ~ e^{2\pi i~\hbox{diag}\left(
{\theta_x\over 2},-{\theta_x\over 2}
\right)}~ 
\left\{C_{\infty1}^{\prime} \pmatrix{ 1 & 0 \cr \cr 0 & {s\over 1+\sigma}  }
 \right\}
\right]
\label{mxxx}
\ee
\be 
m_1:= C_{01}^{-1} ~ e^{2\pi i~\hbox{diag}\left(
{\theta_1\over 2},-{\theta_1\over 2}
\right)} ~C_{01}
\label{m111}
\ee
From the above results, from the formulae  (\ref{*}), (\ref{**}), (\ref{***}),
 and from  the cyclic properties of the trace we obtain 
$$ \hbox{tr}(M_0M_x) \equiv \hbox{tr}(m_0m_x),~~~
 \hbox{tr}(M_1M_x)\equiv \hbox{tr}(m_1m_x),
~~~\hbox{tr}(M_0M_1)\equiv \hbox{tr}(m_0m_1)
$$

\vskip 0.3 cm 
 We note that in order for the $\Gamma$-functions to 
be non-singular we have to require $(\pm \sigma \pm \theta_1 
\pm \theta_{\infty})/2 \not \in {\bf Z} $,  
 $(\pm \sigma \pm \theta_0 \pm \theta_x)/2 \not \in {\bf Z} $. Therefore, to  
summarize we define the {\it generic case}  to be  
\be
   \sigma, ~\theta_0,~\theta_x,~\theta_1,~\theta_{\infty}~\not \in {\bf Z};
~~~~
{\pm \sigma \pm \theta_1 \pm \theta_{\infty}\over 2},~
{\pm \sigma \pm \theta_0 \pm \theta_x \over 2} ~\not \in {\bf Z}
\label{GENERICCASE}
\ee

\vskip 0.2 cm
\noindent
 Computing the traces, we find
\be 
  \hbox{tr}( M_0M_x)= 2 \cos (\pi \sigma) 
\label{overrama0}
\ee
\be 
   \hbox{tr}( M_1M_0)= F_1(\sigma,\theta_0,\theta_x,\theta_1,\theta_{\infty})
~{1\over s} +  F_2(\sigma,\theta_0,\theta_x,\theta_1,\theta_{\infty})+ 
 F_3(\sigma,\theta_0,\theta_x,\theta_1,\theta_{\infty})~ s
\label{overrama}
\ee
\be
 \hbox{tr}( M_xM_1)= -e^{-i\pi \sigma} 
F_1(\sigma,\theta_0,\theta_x,\theta_1,\theta_{\infty})
~{1\over s} +  F_4(\sigma,\theta_0,\theta_x,\theta_1,\theta_{\infty})
- e^{i\pi\sigma}F_3(\sigma,\theta_0,\theta_x,\theta_1,\theta_{\infty})~ s
\label{overrama1}
\ee
where $F_i(\sigma,\theta_0,\theta_x,\theta_1,\theta_{\infty})$ are long 
expressions which can be explicitly computed  from (\ref{m000}), 
 (\ref{mxxx}), (\ref{m111}) in terms of $\Gamma$ functions (
and trigonometric functions). We omit the explicit formulae because they keep 
 some 
space; in any case, 
they are easily computed with a computer from (\ref{m000}), 
 (\ref{mxxx}), (\ref{m111}). The first equation  determines
 $\sigma$ up to $\sigma \mapsto \pm \sigma +2 n$, 
$n\in {\bf Z}$. Once $\sigma$ is chosen, the last two equations determine 
\be 
  s= {e^{i\pi \sigma}  \hbox{tr}( M_1M_0) +\hbox{tr}( M_xM_1) - F_4 - 
e^{-i\pi \sigma} F_2 \over (e^{i\pi\sigma}-e^{-i\pi\sigma})F_1}
\label{TROPPODIFFICILE}
\ee
 If we substitute $s$ in the coefficient $a$ given by (\ref{F}) 
  we get 
$$
a=a(\sigma; \theta_0, \theta_x,\theta_1,\theta_{\infty}, 
\hbox{tr}(M_0M_x),  \hbox{tr}(M_0M_1),
 \hbox{tr}(M_1M_x))
$$
 namely, we get  the parameterization $^2$ 
\be
  y(x;\sigma,a)= 
y\bigl(x;~ \theta_0, \theta_x,\theta_1,\theta_{\infty}, \hbox{tr}(M_0M_x),  \hbox{tr}(M_0M_1),
 \hbox{tr}(M_1M_x)~ \bigr)
\label{P}
\ee

\vskip 0.3 cm

\noindent
{\it Remark: } Formula (\ref{overrama0}) determines $\sigma$ up to $\sigma\mapsto 
\pm\sigma+2n$, $n\in {\bf Z}$. We can choose $
0 \leq\Re \sigma\leq 1
$ 
and therefore all solutions of (\ref{overrama0}) are  $\pm\sigma+2n$, $n\in {\bf Z}$. 
For the given monodromy data, the transcendent (\ref{P}) 
 has different
 representations  
$$
y\bigr(~x;~\pm \sigma+2n,~a(\pm\sigma+2n,~ \theta_0, \theta_x,\theta_1,\theta_{\infty}, \hbox{tr}(M_0M_x),  \hbox{tr}(M_0M_1),
 \hbox{tr}(M_1M_x)) ~\bigl)
$$ 
on different domains $D(\epsilon_n; \pm\sigma+2n)$ of theorem A1. 

\vskip 0.3 cm

\vskip 0.2 cm
\noindent
{\bf Proposition A1:} {\it Let us consider the generic case 
(\ref{GENERICCASE}) and let 
 $\sigma \not \in 
 (-\infty,0]\cup [1,+\infty)$ and $a\neq 0$.  Let $y(x)$ be a solution 
of $PVI$ such that 
 $y(x)= a x^{1-\sigma}$(1+higher order terms)  as $x\to 0 $ in an open   
 domain contained in  $D(\epsilon;\sigma)$ of theorem A1. 
  Then, $y(x)$ coincides with $y(x;\sigma,a)$ of
Theorem A1} 
\vskip 0.2 cm
\noindent
{\it  Proof:}  Observe that both $y(x)$ and $y(x;\sigma,a)$ have the same
asymptotic behavior for $x\to 0$ in $D(\sigma)$. Let  
 $A_0(x)$, $A_1(x)$, $A_x(x)$ be the matrices constructed from
 $y(x)$ and $A^*_0(x)$, $A^*_1(x)$, $A^*_x(x)$ constructed from
 $y(x;\sigma,a)$ by means of   
 the formulae of the Appendix C of \cite{JM1}. 
  It follows that $A_i(x)$ and $A^*_i(x)$, $i=0,1,x$, 
 have the same asymptotic behavior as
 $x\to 0$. This is the
behavior of Lemma A1. Therefore, according to the construction above (reduction to the hyper-geometric equation)   $A_0(x)$,
$A_1(x)$, $A_x(x)$  and $A^*_0(x)$, $A^*_1(x)$, $A^*_x(x)$ give 
 the same monodromy.  
The solution of the Riemann-Hilbert problem for such monodromy is unique, up 
to diagonal conjugation of the fuchsian systems.  
Therefore $A_i(x)$ and $ A^*_i(x)$, $i=0,1,x$ are at most diagonally 
conjugated and  
  $[A(z;x)]_{12}=0$ and  $ [A^*(z;x)]_{12}=0$ give  
$y(x)\equiv y(x;\sigma,a)$. 
\rightline{$\Box$}

\vskip 0.5 cm 
\noindent 
{\bf CONNECTION PROBLEM:} ~~We explain how
 to solve the connection problem for the
 transcendent $y(x;\sigma,a)$.

\vskip 0.2 cm 
\noindent
{\bf [$x=\infty$]}: Let 
$$
 t:={1\over x},~~~  y(x) =: {1\over  t}~\tilde{y}(t)
$$
Then $y(x)$ solves PVI with parameters $\theta_0$, $\theta_x$, $\theta_1$, $\theta_{\infty}$  if and only if $\tilde{y}(t)$ solves PVI (with  
independent  variable $t$) with parameters
\be
  (\tilde{\theta}_0,\tilde{\theta}_t, \tilde{\theta}_1, 
\tilde{\theta}_{\infty} ) := (\theta_0,\theta_1,\theta_x,\theta_{\infty})
\label{itetainf}
\ee
Namely, $\theta_1$ and $\theta_x$ are exchanged.  From  Theorem A1 if follows that there exists a solution $\tilde{y}(t;\sigma^{(\infty)},a^{(\infty)})= 
a^{(\infty)} t^{1-\sigma^{(\infty)}}(1+O(t^{\Delta}))$  
for $t\to 0$. Therefore, PVI in variable $x$ has a solution 
$$ 
  y(x) = a^{(\infty)} x^{\sigma^{(\infty)}}(1+ O(x^{-\Delta})),~~~x\to \infty
$$ 
in   $$ 
    D(M; \sigma^{(\infty)};\vartheta_1,\vartheta_2,\tilde{\sigma}):=\{ 
 {x}\in 
\widetilde{ {\bf C} \backslash \{\infty\} } 
\hbox{ s.t. } 
|x|>M,~ 
e^{ -\vartheta_1\Im\sigma^{(\infty)} }
 |x|^{-\tilde{\sigma}}
\leq |{x}^{ -\sigma^{(\infty)}}|
\leq  e^{-\vartheta_2\Im\sigma^{(\infty)}}
 $$
$$
0<
   \tilde{\sigma} < 1
\}
$$
where $M>0$ is sufficiently big and $0<\Delta<1$ is small. 

\vskip 0.2 cm
\noindent 
{\bf [$x=1$]}: Now let  
$$
   t:=1-x,~~~~y(x)=:1-\tilde{y}(t)
$$
 Then $y(x)$ satisfies PVI  with parameters $\theta_0$, $\theta_x$, $\theta_1$, $\theta_{\infty}$  if and only if $\tilde{y}(t)$ solves PVI 
(with 
independent  variable $t$) with parameters
\be
(\tilde{\theta}_0,\tilde{\theta}_t, \tilde{\theta}_1, 
\tilde{\theta}_{\infty} ) := (\theta_1,\theta_x,\theta_0,\theta_{\infty})
\label{iteta1}
\ee
Namely, $\theta_0$ and $\theta_1$ are exchanged.  From  Theorem A1 if follows that there exists a solution $\tilde{y}(t;\sigma^{(1)},a^{(1)})= 
a^{(1)} t^{1-\sigma^{(1)}}(1+O(t^{\Delta}))$  
for $t\to 0$. Therefore, PVI in variable $x$ has a solution
$$
y(x)=1-a^{(1)}(1-x)^{1-\sigma^{(1)}}~\bigl(1+O((1-x)^{\Delta})\bigr)
$$
in 
$$
  D(\epsilon;\sigma^{(1)};\vartheta_1,\vartheta_2,\tilde{\sigma}):= \{ {x}\in
   \widetilde{{\bf C}\backslash \{1\}} 
   \hbox{ s.t. } |1-x|<\epsilon ,~~ e^{-\vartheta_1 \Im \sigma}
   |1-x|^{\tilde{\sigma}} \leq |(1-{x})^{\sigma^{(1)}}| 
 \leq e^{-\vartheta_2 \Im
   \sigma} , 
$$ 
$$0 < \tilde{\sigma} < 1 \}$$  

\vskip 0.3 cm 
 Let us examine how the monodromy data change. 

\vskip 0.2 cm 
\noindent
{\bf [$x=\infty$]}: 
In the fuchsian system 
(\ref{stofucs}) we put $x={1\over t}$ and $\tilde{z}:= t z$. We get 
$$ 
    {dY \over d \tilde{z}} = \left[{A_0\over \tilde{z}}+{A_x\over \tilde{z}-1} 
+ {A_1 \over \tilde{z}-t} \right]~Y
$$ 
 Therefore $z=1, x$ have been exchanged and correspond to
 $\tilde{z} = t,1$. The monodromy matrices $\tilde{M}_0$,  $\tilde{M}_t$,  $\tilde{M}_1$ 
corresponding to the loops  
 ordered as in figure \ref{figure2} with $t$ instead of $x$ are
 related to $M_0$, $M_1$, $M_x$ of 
 (\ref{stofucs}) in the basis of figure \ref{figure2}:
$$ 
   M_0=\tilde{M}_0,~~M_x=\tilde{M}_t^{-1} \tilde{M}_1 \tilde{M}_t,~~
M_1=\tilde{M}_t 
$$
\begin{figure}
\epsfxsize=12cm
\centerline{\epsffile{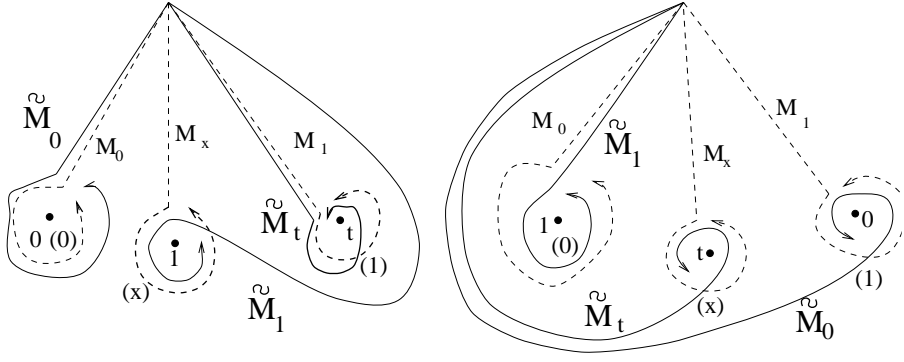}}
\caption{The basis of loops for the system in $z$  and for the system in  
$\tilde{z}$. In parenthesis are  the poles of $z$, out of parenthesis the poles of $\tilde{z}$. The left figure corresponds to $t=1/x$, $\tilde{z}=tz$. The right figure corresponds to $t=1-x$, $\tilde{z}=1-z$. The dotted loops are for the system in $z$, the others for the system in $\tilde{z}$}
\label{figure2}
\end{figure}
If we  introduce the convenient notation 
$$ 
  T_0:=\hbox{tr}(M_0M_x),~~~T_1:=\hbox{tr}(M_1M_x),~~~
T_{\infty}:=\hbox{tr}(M_0M_1)
$$
we have 
$$ 
   \left\{ \matrix{ 
T_0= \hbox{tr}(\tilde{M}_0\tilde{M}_t^{-1}\tilde{M}_1 \tilde{M}_t) \cr\cr
T_1 = \tilde{T}_1   \cr\cr
T_{\infty} = \tilde{T}_0
}
\right.
$$
where $\tilde{T}_1= \hbox{tr}(\tilde{M}_1\tilde{M}_t)$, etc. 
Now we observe that 
$$ 
 \hbox{tr}(\tilde{M}_0\tilde{M}_t^{-1}\tilde{M}_1 \tilde{M}_t)= 4\bigl[ 
\cos(\pi \tilde{\theta}_{\infty}) \cos(\pi \tilde{\theta}_t) +
 \cos(\pi \tilde{\theta}_0)\cos(\pi \tilde{\theta}_1)
\bigr]
- \hbox{tr}(\tilde{M}_1\tilde{M}_0)-\hbox{tr}(\tilde{M}_1 \tilde{M}_t) 
\hbox{tr}(\tilde{M}_0 \tilde{M}_t) 
$$
This follows from the identity 
$$ 
 \hbox{tr}(AB)= \hbox{tr}(A) \hbox{tr}(B) - \hbox{tr} (AB^{-1}), 
~~~~A,B \hbox{ $2\times 2$ matrices }, ~~~ \det(B)=1
$$ 
and from
$$
\hbox{tr}(\tilde{M}_1\tilde{M}_t\tilde{M}_0)= e^{ i \pi \tilde{\theta}_{\infty} }+ 
 e^{ -i \pi \tilde{\theta}_{\infty} },
~~~~\hbox{tr}(\tilde{M}_i) = e^{ i \pi \tilde{\theta}_i }+ 
 e^{ -i \pi \tilde{\theta}_i },~~i=0,x,1
$$
Therefore, recalling (\ref{itetainf}),  we conclude that 
$$ 
\left\{ \matrix{
\tilde{T}_0=T_{\infty}\cr
\tilde{T}_1= T_1 \cr
\tilde{T}_{\infty} =
4 \bigl[\cos(\pi \theta_{\infty}) \cos(\pi \theta_1) +
\cos(\pi \theta_x) \cos(\pi \theta_0) \bigr]
- ( T_0+T_1 T_{\infty})
}
\right.
 $$
Therefore  the parameterization of $\sigma^{(\infty)}$, $a^{(\infty)}$ 
in terms of the monodromy data $\theta_{\nu}$ ($\nu=0,x,1,\infty$), $T_0$, $T_1$, $T_{\infty}$ is obtained as  follows.  
 As for $\sigma^{(\infty)}$ we get if from  
$$ 
 2 \cos (\pi \sigma^{(\infty})) = T_{\infty}
$$
As for $a^{(\infty)}$, we compute if from the formulae 
 (\ref{TROPPODIFFICILE}) (\ref{F})   with the 
substitutions $^3$ 
$$ 
  \sigma\mapsto \sigma^{(\infty)}$$
$$
    \theta_x \mapsto \theta_1,~~~\theta_1\mapsto \theta_x
$$
$$
  T_0 \mapsto T_{\infty}, ~~~ T_{\infty} \mapsto 
4 \bigl[\cos(\pi \theta_{\infty}) \cos(\pi \theta_1) +
\cos(\pi \theta_x) \cos(\pi \theta_0) \bigr]
- ( T_0+T_1 T_{\infty})
$$

\vskip 0.3 cm 
\noindent
{\bf [$x=1$]}: If we put $x=1-t$, $z=1-\tilde{z}$, the system  (\ref{stofucs}) 
becomes 
$$ 
    {dY \over d \tilde{z}} = \left[{A_1\over \tilde{z}}+{A_0\over \tilde{z}-1} 
+ {A_x \over \tilde{z}-t} \right]~Y
$$ 
Therefore $z=0,1$ have been exchanged. 
The monodromy matrices $\tilde{M}_0$,  $\tilde{M}_t$,  $\tilde{M}_1$ are 
$$ 
  M_0= \tilde{M}_1, ~~~M_x= \tilde{M}_1 \tilde{M}_t\tilde{M}_1^{-1},
~~~
M_1=\tilde{M}_1\tilde{M}_t\tilde{M}_0 \tilde{M}_t^{-1} \tilde{M}_1^{-1}
$$
As above, this implies  that 
$$ 
\left\{ \matrix{
\tilde{T}_0=T_1 \cr
\tilde{T}_1= T_0 \cr
\tilde{T}_{\infty} =
4\left[\cos(\pi \theta_{\infty})\cos(\pi \theta_x) +
\cos(\pi \theta_0)\cos(\pi \theta_1)
 \right]-(T_{\infty} +T_0T_1)
\cr
}
\right.
$$
Therefore, we obtain $\sigma^{(1)}$ from 
$$ 
  2 \cos(\pi \sigma^{(1)})=T_1$$
and $a^{(1)} $ from the formulae   (\ref{TROPPODIFFICILE}) (\ref{F}) 
with  the substitutions $^4$:
$$
  \sigma \mapsto \sigma^{(1)}
$$
$$
   \theta_0 \mapsto \theta_1,~~~\theta_1 \mapsto \theta_0 
$$
$$
T_0 \mapsto T_1,~~~T_1\mapsto T_0,~~~T_{\infty} \mapsto 
4\left[\cos(\pi \theta_{\infty})\cos(\pi \theta_x) +
\cos(\pi \theta_0)\cos(\pi \theta_1)
 \right]-(T_{\infty} +T_0T_1)
$$

\vskip 0.2 cm 
This solves the connection problem for $y(x;\sigma,a)$, which is 
 the representation close to $x=0$ 
in $D(\epsilon;\sigma)$  of $y(x;\theta_0,\theta_x,
\theta_1,\theta_{\infty}; T_0,T_1,T_{\infty})$  for  monodromy data computed 
through (\ref{overrama0}), (\ref{overrama}), (\ref{overrama1}). 
 From the monodromy data we can 
compute $a^{(1)}, \sigma^{(1)}$,  $a^{(\infty)}, \sigma^{(\infty)}$ as we just explained above.  Therefore  $y(x;\sigma,a)$ 
has representations in some 
$D(\epsilon;\sigma^{(1)})$ and $D(M;\sigma^{(\infty)})$ at $x=1$ and 
$x=\infty$ respectively, 
 with parameters $a^{(1)}, \sigma^{(1)}$,  $a^{(\infty)}, \sigma^{(\infty)}$.

\vskip 1 cm 
{\bf NOTES}

\vskip 0.2 cm 
\noindent
NOTE 1: Frobenius manifolds 
 are the geometrical setting for  the WDVV equations
and were introduced by Dubrovin in \cite{Dub4}. They  are an 
important object in   many branches of
mathematics  like singularity theory and 
 reflection groups \cite{Saito} \cite{SYS}  \cite{Dub6}  \cite{Dub1},
algebraic and 
enumerative geometry \cite{KM} \cite{Manin}. 

\vskip 0.2 cm 
\noindent
NOTE 2: 
The limit of $a$,  for $\theta_0,\theta_1, \theta_x \to 0$, 
$\theta_{\infty}=2\mu$,  computed from 
the formula (\ref{TROPPODIFFICILE}), exists  and it 
coincides with (\ref{AHAHAH})!
 
\vskip 0.2 cm 
\noindent
NOTE 3: The substitutions we have obtained now 
 reduce to  $(x_0,x_1,x_{\infty}) 
\mapsto (x_{\infty}, -x_1,x_0-x_1x_{\infty})$ 
in the non-generic case $\beta=\gamma=1-2\delta=0$.

\vskip 0.2 cm 
\noindent
NOTE 4: The substitutions we have obtained now 
 reduce to    $(x_0,x_1,x_{\infty}) 
\mapsto (x_1,x_0,x_0x_1-x_{\infty})$   
in the non-generic case $\beta=\gamma=1-2\delta=0$.

\vskip 1 cm
{\it Acknowledgments:}  The author 
is supported by a fellowship of 
the Japan Society for the Promotion of Science (JSPS).


\end{document}